\documentclass[10pt]{article}

\usepackage{amsmath}
\usepackage{amsfonts}
\usepackage{amssymb}
\usepackage[english]{babel}
\usepackage{amsthm}
\usepackage{extarrows}
\usepackage{amscd}
\usepackage{tikz}
\usepackage{mathdots}
\usetikzlibrary{arrows,shapes}

\def\pf{\par\noindent {\bf Proof}~\par\noindent}
\def\qed{~\hfill{$\square$}\pagebreak[1]\par\medskip\par}

\newcommand{\mR}{\mathbb{R}}
\newcommand{\mC}{\mathbb{C}}

\newcommand{\mS}{\mathbb{S}}
\newcommand{\mH}{\mathbb{H}}
\newcommand{\mI}{\mathbb{I}}
\newcommand{\mJ}{\mathbb{J}}
\newcommand{\mK}{\mathbb{K}}
\newcommand{\mQ}{\mathbb{Q}}
\newcommand{\mT}{\mathbb{T}}

\newcommand{\gf}{\mathfrak{f}}
\newcommand{\gfd}{\mathfrak{f}^{\dagger}}
\newcommand{\gsl}{\mathfrak{sl}}
\newcommand{\gsp}{\mathfrak{sp}}
\newcommand{\gso}{\mathfrak{so}}
\newcommand{\gspin}{\mathfrak{spin}}
\newcommand{\gu}{\mathfrak{u}}
\newcommand{\gh}{\mathfrak{h}}
\newcommand{\geg}{\mathfrak{g}}
\newcommand{\gl}{\mathfrak{gl}}

\newcommand{\ol}{\overline}

\newcommand{\uX}{\underline{X}}

\newcommand{\uz}{\underline{z}}
\newcommand{\uzd}{\underline{z}^{\dagger}}

\newcommand{\upz}{\partial_{\uz}}
\newcommand{\upzJ}{\partial_{\uz}^J}

\newcommand{\upzd}{\partial_{\uz}^{\dagger}}
\newcommand{\upzJd}{\partial_{\uz}^{\dagger J}}

\newcommand{\p}{\partial}
\newcommand{\dirac}{\underline{\p}}

\newcommand{\eop}{\hfill$\square$}

\newtheorem{theorem}{Theorem}
\newtheorem{lemma}{Lemma}
\newtheorem{proposition}{Proposition}
\newtheorem{definition}{Definition}
\newtheorem{remark}{Remark}
\newtheorem{corollary}{Corollary}

\newtheorem{defin}{Definition}

\begin{document}

\title{Fundaments of Quaternionic Clifford Analysis I}
\author{F.\ Brackx$^\ast$, H.\ De Schepper$^\ast$, D.\ Eelbode$^{\ast\ast}$, R.\ L\'{a}vi\v{c}ka$^\ddagger$ \& V.\ Sou\v{c}ek$^\ddagger$}

\date{\small{$^\ast$ Clifford Research Group, Dept. of Math. Analysis, 
Faculty of Engineering and Architecture,\\ Ghent University,
Building S22, Galglaan 2, B--9000 Gent, Belgium\\
$\ast\ast$ University of Antwerp, Middelheimlaan 2, Antwerpen, Belgium\\
$^\ddagger$ Mathematical Institute, Faculty of Mathematics and Physics, Charles University\\
Sokolovsk\'a 83, 186 75 Praha, Czech Republic}}

\maketitle


\begin{abstract}
Introducing a quaternionic structure on Euclidean space, the fundaments for quaternionic and symplectic Clifford analysis are studied in detail from the viewpoint of invariance for the symplectic group action.
\end{abstract}


\section{Introduction}


Recently a number of papers \cite{ Ee, PSS, DES, ABBSS1, ABBSS2} appeared dealing with so--called (hermitian) quaternionic  monogenic functions. These are functions defined in Euclidean space, the dimension of which is assumed to be a fourfold, and taking values in a Clifford algebra, or subspaces thereof, which are null solutions of four first order differential operators: a quaternionic Dirac operator and three different conjugates of it. The associated function theory is termed (hermitian) quaternionic Clifford analysis. Let us situate this quaternionic Clifford analysis within the still growing but already well established domain of Clifford analysis.\\[-2mm]

Standard Clifford analysis, also called Euclidean or orthogonal Clifford analysis, is, in its most basic form, a higher dimensional generalization of holomorphic function theory in the complex plane, and a refinement of harmonic analysis, see e.g. \cite{bds, gilbert, dss, guerleb, ghs}. At the heart of this function theory lies the notion of a monogenic function, i.e.\ a Clifford algebra valued null solution of the Dirac operator $\dirac = \sum_{\alpha=1}^m e_{\alpha} \, \p_{X_{\alpha}}$, where $(e_1, \ldots, e_m)$ is an orthonormal basis of $\mR^m$, which underlies the construction of the real Clifford algebra $\mR_{0,m}$. We refer to this setting as the Euclidean case, since the fundamental group leaving the Dirac operator $\dirac$ invariant, is the orthogonal group $\mbox{O}(m)$, which is doubly covered by the Pin($m$) group in the Clifford algebra.\\[-2mm]

In the books \cite{rocha,struppa} and the series of papers \cite{sabadini,eel,partI,partII,toulouse,eel2,howe} so--called hermitian Clifford analysis emerged as a refinement of Euclidean Clifford analysis, by considering functions now taking their values in the complex Clifford algebra $\mC_{2n}$ or in complex spinor space. Hermitian Clifford analysis is based on the introduction of an additional datum, a so--called complex structure $\mI$, inducing an associated Dirac operator $\dirac_{\mI}$; it then focusses on the simultaneous null solutions of both operators $\dirac$ and $\dirac_{\mI}$, called hermitian monogenic functions. The fundamental group underlying this function theory is the unitary group U$(n)$. It is worth mentioning that the traditional holomorphic functions of several complex variables are a special case of hermitian monogenic functions.\\[-2mm]

In this paper we will show how quaternionic Clifford analysis arises in a natural way by introducing on 
$\mR^{4p}$ a so--called quaternionic structure $\mQ = (\mI, \mJ, \mK)$, where $\mJ$ is a second complex structure anti--commuting with $\mI$, and, obviously, $\mK$ is the composition of $\mI$ and $\mJ$ (Section 3). With an eye on the Fischer decomposition of homogeneous polynomials into irreducible representations of the symplectic group Sp$(p)$, the value space of the functions considered, namely the spinor space $\mS$ of the complex Clifford algebra $\mC_{4p}$, is split into  Sp$(p)$--irreducible pieces (Section 4). In Section 5 the concept of a quaternionic monogenic function is introduced in the setting based on the quaternionic structure; this is done stepwise passing from Euclidean, via hermitian, to quaternionic Clifford analysis. It should be emphasized that, while in this way laying the fundaments of quaternionic Clifford analysis, the algebra of quaternion numbers is only lurking in the background, as it is indeed our explicit aim to develop this function theory at the level of complex Clifford algebra or complex spinor space. This approach is, mutatis mutandis, similar to developing hermitian Clifford analysis at the level of real Clifford algebra making use of the above mentioned Dirac operators $\dirac$ and $\dirac_{\mI}$. In order to making this paper self--contained, the basics of Clifford algebra are recalled in Section 2.


\section{Preliminaries on Clifford algebra}


For a detailed description of the structure of Clifford algebras we refer to e.g.\ \cite{port}. Here we only recall the necessary basic notions. The real Clifford algebra $\mathbb{R}_{0,m}$ is constructed over the vector space $\mathbb{R}^{0,m}$ endowed with a non--degenerate quadratic form of signature $(0,m)$, and generated by the orthonormal basis $(e_1,\ldots,e_m)$. The non--commutative Clifford or geometric multiplication in $\mathbb{R}_{0,m}$ is governed by the rules 
\begin{equation}\label{multirules}
e_{\alpha} e_{\beta} + e_{\beta} e_{\alpha} = -2 \delta_{\alpha \beta} \ \ , \ \ \alpha,\beta = 1,\ldots ,m
\end{equation}
As a basis for $\mathbb{R}_{0,m}$ one takes for any set $A=\lbrace j_1,\ldots,j_h \rbrace \subset \lbrace 1,\ldots,m \rbrace$ the element $e_A = e_{j_1} \ldots e_{j_h}$, with $1\leq j_1<j_2<\cdots < j_h \leq m$, together with $e_{\emptyset}=1$, the identity element. The dimension of $\mR_{0,m}$ is $2^m$. Any Clifford number $a$ in $\mathbb{R}_{0,m}$ may thus be written as $a = \sum_{A} e_A a_A$, $a_A \in \mathbb{R}$, or still as $a = \sum_{k=0}^m \lbrack a \rbrack_k$, where $\lbrack a \rbrack_k = \sum_{|A|=k} e_A a_A$ is the so--called $k$--vector part of $a$.\\
Real numbers correspond with the zero--vector part of the Clifford numbers. Euclidean space $\mathbb{R}^{0,m}$ is embedded in $\mathbb{R}_{0,m}$ by identifying $(X_1,\ldots,X_m)$ with the Clifford $1$--vector $\uX = \sum_{\alpha=1}^m e_{\alpha}\, X_{\alpha}$, for which it holds that $\uX^2 = - |\uX|^2$.\\[-2mm]

When allowing for complex constants, the generators $(e_1,\ldots, e_{m})$, still satisfying (\ref{multirules}), produce the complex Clifford algebra $\mathbb{C}_{m} = \mathbb{R}_{0,m} \oplus i\, \mathbb{R}_{0,m}$. Any complex Clifford number $\lambda \in \mathbb{C}_{m}$ may thus be written as $\lambda = a + i b$, $a,b \in \mathbb{R}_{0,m}$, leading to the definition of the hermitian conjugation $\lambda^{\dagger} = (a +i b)^{\dagger} = \overline{a} - i \overline{b}$, where the bar notation stands for the Clifford conjugation in $\mathbb{R}_{0,m}$, i.e. the main anti--involution for which $\overline{e}_{\alpha} = -e_{\alpha}$, $\alpha=1, \ldots,m$. This hermitian conjugation leads to a hermitian inner product on $\mathbb{C}_{m}$ given by $(\lambda,\mu) = \lbrack \lambda^{\dagger} \mu \rbrack_0$ and its associated norm $|\lambda| = \sqrt{ \lbrack \lambda^{\dagger} \lambda \rbrack_0} = ( \sum_A |\lambda_A|^2 )^{1/2}$.\\[-2mm]

The algebra of real quaternions is denoted by $\mH$. For a quaternion
$$
q = q_0 + q_1 i + q_2 j + q_3 k = (q_0 + q_1 i )+ (q_2 + q_3 i)j = z + wj
$$
its conjugate is given by 
$$
\ol{q} = q_0 - q_1 i - q_2 j - q_3 k = (q_0 - q_1 i) - j(q_2  - q_3 i) = \ol{z} -j \ol{w}
$$
such that
$$
q \ol{q} = \ol{q} q = |q|^2 = q_0^2 + q_1^2 + q_2^3 + q_3^2 = |z|^2 + |w|^2
$$
Identifying the quaternion units $i, j$ with the respective basis vectors $e_1, e_2$, the algebra $\mH$ is isomorphic with the Clifford algebra $\mR_{0,2}$. It is also isomorphic with the even subalgebra $\mR_{0,3}^+$
of the Clifford algebra $\mR_{0,3}$, by identifying the quaternion units $i, j,k$ with the respective bivectors $e_2 e_3, e_3 e_1, e_1 e_2$.


\section{Quaternionic structure}


\subsection{A first complex structure}
\label{subsecQ1}


Let $\mR^{4p}$ be the Euclidean space of dimension $m = 2n = 4p$, with the standard orthonormal basis $(e_1, e_2, \ldots, e_{4p})$. It is isomorphic with $\mC^{2p}$ and with $\mH^p$, the natural bijections being given by
\begin{eqnarray*}
\alpha_{2p} &:& \mC^{2p} \rightarrow \mR^{4p} : (a_1+b_1 i, \ldots , a_{2p} + b_{2p} i) \mapsto (a_1, b_1, \ldots , a_{2p}, b_{2p})\\
\beta_{p} &:& \mH^{p} \rightarrow \mC^{2p} \, \, : (q_0+q_1 i + (q_2+q_3 i )j, \ldots) \mapsto (q_0+q_1 i,  q_2+q_3 i, \ldots)\\
\gamma_p = \alpha_{2p} \circ \beta_p &:& \mH^p \rightarrow \mR^{4p} \, \, : (q_0+q_1 i + q_2 j +q_3 ik, \ldots) \mapsto (q_0, q_1, q_2, q_3, \ldots)
\end{eqnarray*}
Square matrices in $M_{2p}(\mC)$ with complex entries are embedded in the space $M_{4p}(\mR)$ of real matrices in the following way: the injective homomorphism $\varphi_{2p}: M_{2p}(\mC) \rightarrow M_{4p}(\mR)$ substitutes for each complex entry $a+bi$ the $2 \times 2$ real matrix $
\begin{pmatrix} \phantom{-}a& b \\ -b & a \end{pmatrix}$.
More precisely, the following scheme, with $A \in M_{2p}(\mC)$ and $B=\varphi_{2p}(A) \in M_{4p}(\mR)$, should be commutative: 
$$
\begin{CD}
\mC^{2p} @>\hspace*{5mm} \alpha_{2p} \hspace*{5mm}>> \mR^{4p} \\
@VT_AVV @VVT_{\varphi_{2p}}(A)V \\
\mC^{2p} @>\hspace*{5mm}\alpha_{2p}\hspace*{5mm}>> \mR^{4p}
\end{CD}
$$
Here $T_A$ stands for the $\mC$--linear transformation of right multiplication by the matrix $A$, the vectors of $\mC^{2p}$ being interpreted as row matrices:
$$
T_A: \mC^{2p} \longrightarrow \mC^{2p} : (z_1,\ldots,z_{2p}) \longmapsto (z_1,\ldots,z_{2p}) A
$$
and, similarly,
$$
T_B: \mR^{4p} \longrightarrow \mR^{4p} : (X_1,\ldots,X_{4p}) \longmapsto (X_1,\ldots,X_{4p})B
$$

Square matrices in $M_p(\mH)$ with quaternion entries are first embedded in $M_{2p}(\mC)$ by the injective homomorphism $\psi_p: M_p(\mH) \rightarrow M_{2p}(\mC)$ substituting for each quaternion entry $z+wj$ the $2 \times 2$ complex matrix $\begin{pmatrix} \phantom{-}z& w \\ -\overline{w} & \overline{z} \end{pmatrix}$, making the following scheme with $A \in M_p(\mH)$ and $\psi_p(A) = B \in M_{2p}(\mC)$ commutative:
$$
\begin{CD}
\mH^p @>\hspace*{7mm} \beta_p \hspace*{7mm} >>\mC^{2p} @>\hspace*{7mm}\alpha_{2p}\hspace*{7mm}>> \mR^{4p} \\
@VT_AVV @VVT_{\psi_p}(A)V @VVT_{\varphi_{2p} \circ \psi_p}(A)V \\
\mH^p @>\hspace*{7mm}\beta_p\hspace*{7mm}>> \mC^{2p} @>\hspace*{7mm}\alpha_{2p}\hspace*{7mm}>> \mR^{4p}
\end{CD}
$$
Here
$$
T_A: \mH^{p} \longrightarrow \mH^{p} : (q_1,\ldots,q_{p}) \longmapsto (q_1,\ldots,q_{p}) A
$$
and 
$$
T_B: \mC^{2p} \longrightarrow \mC^{2p} : (z_1,\ldots,z_{2p}) \longmapsto (z_1,\ldots,z_{2p})B
$$
By the composition $\varphi_{2p} \circ \psi_p: M_p(\mH) \rightarrow M_{4p}(\mR)$, square matrices with quaternion entries are directly embedded in the space $M_{4p}(\mR)$ of real matrices. Matrices in $\varphi_{2p}(M_{2p}(\mC)) \subset M_{4p}(\mR)$ are called {\em complex linear real matrices}, matrices in $\psi_p(M_p(\mH)) \subset M_{2p}(\mC)$ are {\em quaternionic linear complex matrices}, while matrices in $\varphi_{2p} \circ \psi_p(M_p(\mH)) \subset M_{4p}(\mR)$ are called {\em quaternionic linear real matrices}.\\[-2mm]

In $\mC^{2p}$ multiplication by the imaginary unit $i$ is the $\mC$--linear transformation associated to the matrix $i E_{2p}$, $E_{2p}$ standing for the identity matrix in $M_{2p}(\mC)$.

\begin{defin}
\label{definitionQ1}
The standard complex structure $\mI_{4p}$ on $\mR^{4p}$ is the complex linear real matrix 
$$
\mI_{4p} \equiv \varphi_{2p}(i E_{2p}) = \mbox{\em{diag}} \begin{pmatrix} \phantom{-} 0 & 1 \\ -1 & 0 \end{pmatrix} = \mbox{\em{diag}} \begin{pmatrix} \phantom{-} \cos(\frac{\pi}{2}) & \sin (\frac{\pi}{2}) \\ - \sin(\frac{\pi}{2}) & \cos(\frac{\pi}{2}) \end{pmatrix}
$$
\end{defin}
\noindent As expected, there holds that $\mI_{4p}^2 = - E_{4p}$, $E_{4p}$ being the identity matrix in $M_{4p}(\mR)$. Moreover, $\mI_{4p}$ belongs to $\mbox{SO}(4p)$.\\[-2mm]

The following result is well--known and easily checked.
\begin{proposition}
\label{propositionQ1}
The following statements are equivalent:
\begin{enumerate}
\item[(i)] the matrix $B \in M_{4p}(\mR)$ is complex linear, i.e.\ $B \in \varphi_{2p}(M_{2p}(\mC))$;
\item[(ii)] the transformation $\alpha^{-1}_{2p} \circ T_B \circ \alpha_{2p}: \mC^{2p} \rightarrow \mC^{2p}$ is $\mC$--linear;
\item[(iii)] the matrix $B$ commutes with the complex structure on $\mR^{4p}$: $B \,\mI_{4p} = \mI_{4p} B$.
\end{enumerate}
\end{proposition}

The following result (see e.g. \cite{partI}) is crucial for the sequel.
\begin{proposition}
\label{propositionQ2}
The $\mbox{\em SO}(4p)$--matrices which are commuting with the complex structure $\mI_{4p}$ on $\mR^{4p}$ form a subgroup of $\mbox{\em SO}(4p)$, denoted by $\mbox{\em SO}_\mI(4p)$, which is isomorphic with the unitary group: $\varphi_{2p}(\mbox{\em U}(2p)) = \mbox{\em SO}_\mI(4p)$.\\[-5mm]
\end{proposition}

\begin{center}
\begin{tikzpicture}
\draw[fill=black, fill opacity=0.1] (0,0) ellipse (1 and 1.4);
\draw[fill=black, fill opacity=0.1] (6.3,0) ellipse (1 and 1.4);
\draw (6.5,0.3) ellipse (1.8 and 2.1);
\node at (0,1.6) {$\mbox{U}(2p)$};
\node at (6.3,1.7) {$\mbox{SO}_\mI(4p)$};
\node at (6.3,2.6) {$\mbox{SO}(4p)$};
\node at (0,0.3) (bul1) {$\bullet$};
\node at (0,0) {$A$};
\node at (6.3,0.3) (bul2) {$\bullet$};
\node at (6.3,0) {$\varphi_{2p}(A)$};
\node at (0,-2.2) {$A \in M_{2p}(\mC)$};
\node at (-0.1,-2.6) {$A^\ast A = E_{2p}$};
\node at (6.7,-2.2) {$\varphi_{2p}(A) \in M_{2p}(\mR)$};
\node at (6.7,-2.6) {commutes with $\mI$};
\draw[->,>=angle 60] (bul1.east) to [out=10,in=170] (bul2.west);
\node at (3.2,0.8) {$\varphi_{2p}$};
\draw[-] (0,1.4) to[out=5,in=175] (6.3,1.4);
\draw[-] (0,-1.4) to[out=5,in=175] (6.3,-1.4);
\end{tikzpicture}
\end{center}

As is well-known the group $\mbox{SO}(4p)$ has a twofold covering by the $\mbox{Spin}(4p)$ group in the Clifford algebra $\mR_{0,4p}$ (or its complexification $\mC_{4p}$); this spin group is easily depicted as the group of all products of an even number of unit vectors. More explicitly, the matrix $A \in \mbox{SO}(4p)$ corresponds to the $\mbox{Spin}(4p)$ element $s_A$ if $(Y_1,\ldots,Y_{4p}) = (X_1 \ldots X_{4p}) A$ and $\underline{Y} = s_A \underline{X} s_A^{-1}$, with $\underline{X} = \sum_{\alpha=1}^{4p} X_\alpha e_\alpha$ and $\underline{Y} = \sum_{\alpha=1}^{4p} Y_\alpha e_\alpha$. \\[-2mm] 

The subgroup $\mbox{SO}_\mI(4p)$ of $\mbox{SO}(4p)$ inherits a twofold covering by the subgroup $\mbox{Spin}_{\mI}(4p)$ of $\mbox{Spin}(4p)$, consisting of those elements of $\mbox{Spin}(4p)$ which are commuting with 
\begin{equation}
s_{\mI} = s_1 \ldots s_{2p},  \quad \mbox{where \ } s_j = \frac{\sqrt{2}}{2} ( 1 + e_{2j-1} e_{2j} ),  \; j=1,\ldots,2p
\label{sI}
\end{equation}
The element $s_{\mI}$ itself obviously belongs to $\mbox{Spin}_{\mI}(4p)$ and corresponds, under the double covering, to the complex structure $\mI_{4p} \in \mbox{SO}(4p)$. Recall that $\mbox{SO}_{\mI}(4p)$ is isomorphic with the unitary group $\mbox{U}(2p)$; up to this isomorphism, $\mbox{Spin}_{\mI}(4p)$ provides a double cover of $\mbox{U}(2p)$.\\[-2mm]

At the level of the Lie algebras we have the following picture. The Lie algebra $\gspin(4p)$ is the space $\mR^{(2)}_{0,4p}$ of all bivectors in the Cliford algebra. Indeed, given two unit vectors $u$ and $v$ enclosing the angle $\alpha$, such that $u \circ v = \cos(\alpha)$, $\|u \wedge v\| = \sin(\alpha)$, there exists the bivector $\alpha \frac{u \wedge v}{\| u \wedge v \|}$ such that
$$
\exp \bigl ( \alpha \frac{u \wedge v}{\| u \wedge v \|} \bigr ) = \cos(\alpha) + \sin(\alpha) \frac{u \wedge v}{\| u \wedge v \|} = u \circ v + u \wedge v = uv \quad \in \mbox{Spin}(4p)
$$
In particular, the bivector corresponding to $s_{\mI} \in \mbox{Spin}(4p)$ is given by
$$
\sigma_{\mI} = \frac{\pi}{4} \bigl ( e_1 e_2 + e_3 e_4 + \cdots + e_{4p-1} e_{4p} \bigr )
$$
since, clearly,
$$
\exp(\sigma_{\mI}) = \frac{\sqrt{2}}{2} (1+ e_1e_2) \frac{\sqrt{2}}{2} (1+ e_3e_4) \cdots \frac{\sqrt{2}}{2} (1+ e_{4p-1} e_{4p}) = s_1 \cdots s_{2p} = s_{\mI} 
$$
The Lie algebra $\gspin_{\mI}(4p)$ corresponding to the Lie group $\mbox{Spin}_{\mI}(4p)$, consists of all bivectors commuting with $\sigma_{\mI}$. The following $4p^2$ bivectors in $\gspin_{\mI}(4p)$:
\begin{eqnarray*}
\sigma_j & = & \frac{\pi}{4} \bigl ( e_{2j-1} e_{2j} \bigr ), \quad j=1,\ldots,2p \\
\sigma_{jk} & = & \frac{\pi}{4} \bigl ( e_{2j-1} e_{2k-1} + e_{2j} e_{2k}  \bigr ), \quad j \ne k = 1,\ldots,2p \\
\widetilde{\sigma_{jk}} & = & \frac{\pi}{4} \bigl ( e_{2j-1} e_{2k} + e_{2k-1} e_{2j}  \bigr ), \quad j \ne k = 1,\ldots,2p 
\end{eqnarray*}
give rise to the following $4p^2$ generators of $\mbox{Spin}_{\mI}(4p)$:
\begin{eqnarray*}
\exp(\sigma_j) &=& s_j, \quad j=1,\ldots,2p\\
\exp(\sigma_{jk}) & = & \frac{1}{2} \bigl ( 1 + e_{2j-1} e_{2k-1} \bigr ) \bigl ( 1 + e_{2j} e_{2k} \bigr ) , \quad j \ne k = 1,\ldots,2p \\
\exp( \widetilde{\sigma_{jk}}) & = & \frac{1}{2} \bigl ( 1 + e_{2j-1} e_{2k} \bigr ) \bigl ( 1 + e_{2k-1} e_{2j} \bigr ), \quad j \ne k = 1,\ldots,2p 
\end{eqnarray*}
The Lie algebra $\gspin(4p)$ of bivectors may be identified with the Lie algebra $\gso(4p)$ of skew--symmetric matrices. In this identification the bivector $\frac{1}{2} e_i e_j$($i<j$) corresponds with the skew--symmetric matrix with 1 as its $(i,j)$th entry (and -1 as its $(j,i)$th entry), all other entries being zero:
\begin{eqnarray*}
\frac{1}{2} e_i e_j \in \gspin(4p) & \longleftrightarrow &  \begin{array}{r} \mbox{\tiny $i$th row $\longrightarrow$ \hspace*{-3mm}} \\[6.5mm] \mbox{\tiny $j$th row $\longrightarrow$ \hspace*{-3mm}} \\[2.5mm] \end{array} 
\left ( \begin{array}{rcrc}  \hspace*{2mm} & & & \\[-2mm] \hspace*{2mm} 0 & \hspace*{2mm} \cdots & 1 & \\  \vdots \hspace*{0.5mm}  &  & \vdots \hspace*{0.5mm}& \\  -1 &\hspace*{2mm}   \cdots & 0 & \\
& & & \end{array} \right ) \in \gso(4p)\\[-3mm]
&& \hspace*{22mm} \uparrow \hspace*{14mm} \uparrow \\[-2mm]
&& \hspace*{17mm} \mbox{\tiny $i$th column } \hspace*{4mm} \mbox{\tiny $j$th column}
\end{eqnarray*}
Since $\mbox{SO}_{\mI}(4p)$ is isomorphic with $\mbox{U}(2p)$, we have the following result.
\begin{proposition}
\label{propositionQ3}
The real unitary Lie algebra $\gu(2p)$ of skew--hermitian matrices is isomorphic with the Lie algebra $\gso_{\mI}(4p)$ of skew--symmetric matrices commuting with the matrix $\log(\mI_{4p}) = \frac{\pi}{2} \mI_{4p}$, which is, in its turn, isomorphic with the Lie algebra $\gspin_{\mI}(4p)$ of bivectors commuting with $\sigma_{\mI}$, or schematically:
$$
\varphi_{2p} \bigl ( \gu(2p) \bigr ) = \gso_{\mI} (4p) \simeq \gspin_{\mI}(4p) 
$$
\end{proposition}

\begin{remark}
The above notation $\log(\mI_{4p})$ is unusual. We consider it to be an easy way to express the relation
$$
\exp \bigl ( \frac{\pi}{2} \mI_{4p} \bigr ) = \mI_{4p}
$$
In fact, this is a special case (for $t=\frac{\pi}{2}$) of the more general relation
$$
\exp \bigl ( t \, \mI_{4p}  \bigr ) = \mbox{\em{diag}} \begin{pmatrix} \phantom{-} \cos(t) & \sin (t) \\ - \sin(t) & \cos(t) \end{pmatrix}
$$
\end{remark}

\noindent {\bf Example} ($p=1$)\\
Consider the $\mbox{Spin}(4)$--element $s_{\mI} = \frac{1}{2} (1+e_1e_2)(1+e_3e_4)$. Its action on the vector $\uX = X_1 e_1 + X_2 e_2 + X_3 e_3 + X_4 e_4$ is given by $s_{\mI} \uX s_{\mI}^{-1} = -X_2 e_1 + X_1 e_2 - X_4 e_3 + X_3 e_4$.  The corresponding $\mbox{SO}(4)$--matrix is
$$
\mI_4 = \left ( \begin{array}{cccc} & 1 & & \\ -1 & & & \\  & & & 1 \\ & & -1 & \end{array} \right )
$$
leading, as expected, to the same action: $\left ( \begin{array}{lccr} \hspace*{-2mm} X_1 & X_2 & X_3 & X_4 \hspace*{-2mm} \end{array} \right ) \, \mI_4 = \left ( \begin{array}{lccr} \hspace*{-2mm} -X_2 &  X_1 & -X_4 & X_3 \hspace*{-2mm} \end{array} \right )$.  Obviously $\mI_4 \in \mbox{SO}_{\mI}(4)$ and the corresponding $\mbox{U}(2)$--matrix is given by
$$
\varphi_2 : A = \left ( \begin{array}{cc} i & \\ & i \end{array} \right ) \in \mbox{U}(2) \longmapsto \mI_4 \in \mbox{SO}_{\mI}(4)
$$
At the level of Lie algebras, the bivector $\sigma_{\mI} \in \gspin_{\mI}(4)$ corresponding to $s_{\mI}\in \mbox{Spin}_{\mI}(4)$, is given by
$$
\sigma_{\mI} = \frac{\pi}{4} ( e_1 e_2 + e_3 e_4)
$$
since, indeed, $\exp(\sigma_{\mI}) = s_{\mI}$. Under the isomorphism of $\gspin_{\mI}(4)$ and $\gso_{\mI}(4)$, the bivector $\sigma_{\mI}$ corresponds to the matrix
$$
i_4 = \frac{\pi}{2} \left ( \begin{array}{cccc} & 1 & & \\ -1 & & & \\  & & & 1 \\ & & -1 & \end{array} \right ) \in \gso_{\mI}(4)
$$
and it it readily verified that, as expected, $\exp(i_4) = \mI_4$. Under the isomorphism of $\gso_{\mI}(4)$ and $\gu(2)$, the matrix $i_4$ corresponds to 
$$
a = \frac{\pi}{2} \left ( \begin{array}{cc} i & \\ & i \end{array} \right ) \in \gu(2)
$$
and, moreover $\exp(a)=A$. Schematically, one has

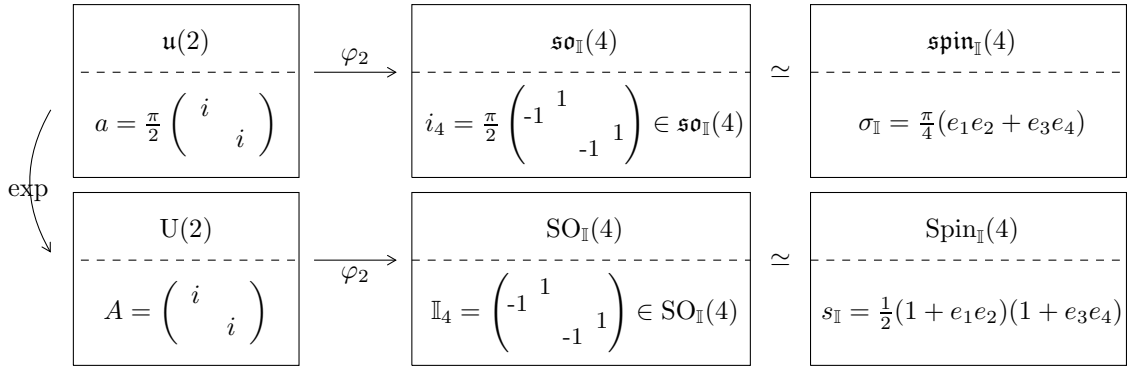
\begin{center}
\begin{tikzpicture}
\node at (0,0) {$\gu(2)$};
\node at (0,-1.1) {$a=\frac{\pi}{2} \left ( \begin{array}{cc} i & \\ & i \end{array} \right )$ };
\draw[dashed] (-1.4,-0.4) -- (1.5,-0.4);
\draw (-1.5,-1.8) -- (-1.5,0.5) -- (1.5,0.5) -- (1.5,-1.8) -- (-1.5,-1.8);
\draw[->,>=angle 60] (1.7,-0.4) to [out=0,in=180] (2.8,-0.4);
\node at (2.25,-0.2) {$\varphi_2$};
\node at (5.3,0) {$\gso_{\mI}(4)$};
\node at (5.3,-1.1) {$i_4 = \frac{\pi}{2} \left ( \begin{array}{cccc} \hspace*{-2mm} &  \hspace*{-2mm} \mbox{\small 1} & \hspace*{-2mm} & \hspace*{-2mm} \\[-2mm] \hspace*{-2mm} \mbox{\small -1} & \hspace*{-2mm} & \hspace*{-2mm} & \hspace*{-2mm} \\[-2mm]  \hspace*{-2mm} & \hspace*{-2mm} & \hspace*{-2mm} & \hspace*{-2mm} \mbox{\small 1} \hspace*{-2mm}\\[-2mm] \hspace*{-2mm} & \hspace*{-2mm} & \hspace*{-2mm} \mbox{\small -1} & \hspace*{-2mm} \end{array} \right ) \in \gso_{\mI}(4)$};
\draw[dashed] (3,-0.4) -- (7.5,-0.4);
\draw (3,-1.8) -- (3,0.5) -- (7.5,0.5) -- (7.5,-1.8) -- (3,-1.8);
\node at (7.9,-0.4) {$\simeq$};
\node at (10.45,0) {$\gspin_{\mI}(4)$};
\node at (10.45,-1.1) {$\sigma_{\mI} = \frac{\pi}{4} (e_1 e_2 + e_3 e_4)$};
\draw[dashed] (8.3,-0.4) -- (12.5,-0.4);
\draw (8.3,-1.8) -- (8.3,0.5) -- (12.5,0.5) -- (12.5,-1.8) -- (8.3,-1.8);
\draw[->,>=angle 60] (-1.8,-0.9) to [out=240,in=120]  (-1.8,-2.8);
\node at (-2.1,-2) {$\exp$};

\node at (0,-2.5) {$\mbox{U}(2)$};
\node at (0,-3.6) {$A= \left ( \begin{array}{cc} i & \\ & i \end{array} \right )$ };
\draw[dashed] (-1.4,-2.9) -- (1.5,-2.9);
\draw (-1.5,-4.3) -- (-1.5,-2) -- (1.5,-2) -- (1.5,-4.3) -- (-1.5,-4.3);
\draw[->,>=angle 60] (1.7,-2.9) to [out=0,in=180] (2.8,-2.9);
\node at (2.25,-3.1) {$\varphi_2$};
\node at (5.3,-2.5) {$\mbox{SO}_{\mI}(4)$};
\node at (5.3,-3.6) {$\mI_4 = \left ( \begin{array}{cccc} \hspace*{-2mm} &  \hspace*{-2mm} \mbox{\small 1} & \hspace*{-2mm} & \hspace*{-2mm} \\[-2mm] \hspace*{-2mm} \mbox{\small -1} & \hspace*{-2mm} & \hspace*{-2mm} & \hspace*{-2mm} \\[-2mm]  \hspace*{-2mm} & \hspace*{-2mm} & \hspace*{-2mm} & \hspace*{-2mm} \mbox{\small 1} \hspace*{-2mm}\\[-2mm] \hspace*{-2mm} & \hspace*{-2mm} & \hspace*{-2mm} \mbox{\small -1} & \hspace*{-2mm} \end{array} \right ) \in \mbox{SO}_{\mI}(4)$};
\draw[dashed] (3,-2.9) -- (7.5,-2.9);
\draw (3,-4.3) -- (3,-2) -- (7.5,-2) -- (7.5,-4.3) -- (3,-4.3);
\node at (7.9,-2.9) {$\simeq$};
\node at (10.45,-2.5) {$\mbox{Spin}_{\mI}(4)$};
\node at (10.45,-3.6) {$s_{\mI} = \frac{1}{2} (1+e_1 e_2)(1 + e_3 e_4)$};
\draw[dashed] (8.3,-2.9) -- (12.5,-2.9);
\draw (8.3,-4.3) -- (8.3,-2) -- (12.5,-2) -- (12.5,-4.3) -- (8.3,-4.3);
\end{tikzpicture}
\end{center}


\subsection{The quaternionic structure}


Consider in $\mH^p$ the left multiplication by  $i$, now regarded as a quaternionic imaginary unit; we denote this transformation by $T_i$. Note that it is a right $\mH$--linear transformation, not a left $\mH$--linear one. Nevertheless the following commutative scheme holds:
$$
\begin{CD}
\mH^{p} @>\hspace*{3mm} \gamma_{p} = \alpha_{2p} \circ \beta_p \hspace*{3mm}>> \mR^{4p} \\
@VT_iVV @VVT_1 \equiv T_{\mI_{4p}}V \\
\mH^{p} @>\hspace*{3mm}\gamma_{p} = \alpha_{2p} \circ \beta_p\hspace*{3mm}>> \mR^{4p}
\end{CD}
$$
where precisely $T_1: \mR^{4p} \longrightarrow \mR^{4p}$ is the linear transformation associated with the standard complex structure $\mI_{4p}$, so $T_1 \equiv T_{\mI_{4p}}$.\\[-2mm]

Similarly, left multiplication in $\mH^p$ with the quaternionic unit $j$ leads to the commutative scheme
$$
\begin{CD}
\mH^{p} @>\hspace*{3mm} \gamma_{p} = \alpha_{2p} \circ \beta_p \hspace*{3mm}>> \mR^{4p} \\
@VT_jVV @VVT_2 \equiv T_{\mJ_{4p}}V \\
\mH^{p} @>\hspace*{3mm}\gamma_{p} = \alpha_{2p} \circ \beta_p\hspace*{3mm}>> \mR^{4p}
\end{CD}
$$
where $\mJ_{4p}$ denotes the $M_{4p}(\mR)$--matrix associated to the linear transformation $T_2: \mR^{4p} \rightarrow \mR^{4p}$. It turns out that
$$
\mJ_{4p} = \mbox{diag} \, \left ( \begin{array}{cccc} &  & 1 & \\ & & & -1 \\  -1 & & &  \\ & 1 & & \end{array} \right )
$$
Clearly $\mJ_{4p}$ belongs to $\mbox{SO}(4p)$, with $\mJ_{4p}^2 = -E_{4p}$, and  and anti--commutes with $\mI_{4p}$. However notice that neither $\mI_{4p}$ nor $\mJ_{4p}$ are quaternionic linear and $\mJ_{4p}$  is even not complex linear.

\begin{definition}
\label{definitionQ2}
The anti--commuting $\mbox{\em SO}(4p)$--matrices $\mI_{4p}$ and $\mJ_{4p}$ form a quaternionic structure on $\mR^{4p}$.
\end{definition}

\begin{remark}
Given the quaternionic structure $(\mI_{4p},\mJ_{4p})$ on $\mR^{4p}$, there immediately arises a third $\mbox{\em SO}(4p)$--matrix 
$$
\mK_{4p} = \mI_{4p} \, \mJ_{4p} = - \mJ_{4p} \, \mI_{4p}
$$ 
for which $\mK_{4p}^2 = -E_{4p}$ and which anti--commutes with both $\mI_{4p}$ and $\mJ_{4p}$. It turns out that 
$$
\mK_{4p} = \mbox{diag} \, \left ( \begin{array}{cccc} &  &  & -1 \\ & & -1 &  \\   & 1 & &  \\  1 & & & \end{array} \right )
$$
is the $M_{4p}(\mR)$--matrix associated to the linear transformation $T_3: \mR^{4p} \rightarrow \mR^{4p}$, given in the following commutative scheme
$$
\begin{CD}
\mH^{p} @>\hspace*{3mm} \gamma_{p} = \alpha_{2p} \circ \beta_p \hspace*{3mm}>> \mR^{4p} \\
@VT_{-k}VV @VVT_3 \equiv T_{\mK_{4p}}V \\
\mH^{p} @>\hspace*{3mm}\gamma_{p} = \alpha_{2p} \circ \beta_p\hspace*{3mm}>> \mR^{4p}
\end{CD}
$$
where $T_{-k}$ stands for the right $\mH$--linear transformation of left multiplication by $-k = j\,i$. Pay attention to the fact that, under operator composition
$$
T_{\mK_{4p}} = T_{\mJ_{4p}} \circ T_{\mI_{4p}}
$$
sometimes also written as $\mK_{4p} = \mJ_{4p} \circ \mI_{4p}$. For symmetry reasons we will henceforth incorporate $\mK_{4p}$ into the quaternionic structure on $\mR^{4p}$ and denote $Q_{4p} = (\mI_{4p},\mJ_{4p},\mK_{4p}) \in \mbox{\em SO}(4p)^3$.
\end{remark}

The following result is well--known and easily checked.
\begin{proposition}
\label{propositionQ4}
The following statements are equivalent:
\begin{enumerate}
\item[(i)] the matrix $B \in M_{4p}(\mR)$ is quaternionic linear, i.e.\ $B \in \varphi_{2p} \circ \psi_p (M_{p}(\mH))$;
\item[(ii)] the transformation $\gamma^{-1}_{p} \circ T_B \circ \gamma_{p}: \mH^{p} \rightarrow \mH^{p}$, with $\gamma_p = \alpha_{2p} \circ \beta_p$,  is (left) $\mH$--linear;
\item[(iii)] the matrix $B$ commutes with the quaternionic structure on $\mR^{4p}$: 
$$
B \, \mI_{4p} = \mI_{4p} B, \quad B \, \mJ_{4p} = \mJ_{4p} B, \quad B \, \mK_{4p} = \mK_{4p} B
$$
\end{enumerate}
\end{proposition}
Also the following result is crucial for the sequel.
\begin{proposition}
\label{propositionQ5}
The $\mbox{\em SO}(4p)$--matrices which are commuting with the quaternionic structure $Q_{4p}$ on $\mR^{4p}$ form a subgroup of $\mbox{\em SO}_{\mI}(4p)$, denoted by $\mbox{\em SO}_Q(4p)$, which is isomorphic with the symplectic group $\mbox{\em Sp}(p)$:
$$
\varphi_{2p} \circ \psi_p \bigl ( \mbox{\em Sp}(p) \bigr ) = \mbox{\em SO}_Q(4p)
$$
\end{proposition}
Recall that the symplectic group $\mbox{Sp}(p)$ is the real Lie group of quaternion $p \times p$ matrices preserving the symplectic inner product, or, equivalently
$$
\mbox{Sp}(p) = \left \{ A \in \mbox{GL}_p(\mH) : AA^\ast = E_p \right \}
$$  
Notice that for $A \in \mbox{Sp}(p)$, $\psi_p(A)$ is a matrix in $\mbox{GL}_{2p}(\mC)$ for which $\varphi_{2p} ( \psi_p(A))$ is an $\mbox{SO}(4p)$--matrix commuting with the complex structure $\mI_{4p}$, hence $\psi_p(A) \in \mbox{U}(2p)$. But since $\det(A) = \det(\psi_p(A)) = +1$, it follows that $\psi_p(A) \in \mbox{SU}(2p)$, and
$$
\psi_{p} \bigl ( \mbox{Sp}(p) \bigr ) = \varphi_{2p}^{-1} \bigl ( \mbox{SO}_Q(4p) \bigr )
$$
is a subgroup of $\mbox{SU}(2p)$. Schematically, one has

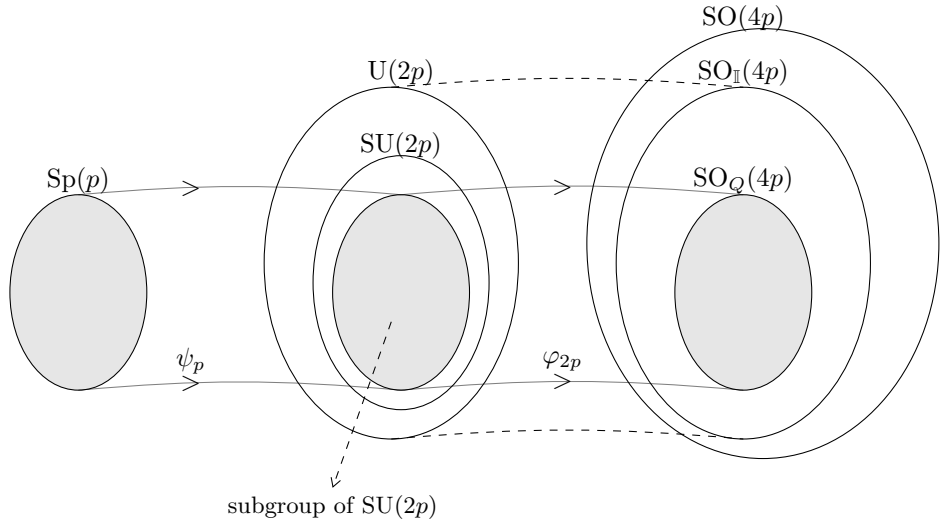
\begin{center}
\begin{tikzpicture}[scale=1.3]
\draw[fill=black,fill opacity =0.1] (-0.3,-0.2) ellipse (0.7 and 1);
\node at (-0.3,0.95) {$\mbox{Sp}(p)$};
\draw[-, draw opacity=0.5] (-0.3,0.8) to [out=5,in=175]  (3,0.8);
\draw[-, draw opacity=0.5] (3,0.8) to [out=5,in=175]  (6.5,0.8);
\draw[-, draw opacity=0.5] (-0.3,-1.2) to [out=5,in=175]  (3,-1.2);
\draw[-, draw opacity=0.5] (3,-1.2) to [out=5,in=175]  (6.5,-1.2);

\draw[fill=black,fill opacity =0.1] (3,-0.2) ellipse (0.7 and 1);
\draw (3,-0.1) ellipse (0.9 and 1.3);
\draw (2.9,0.1) ellipse (1.3 and 1.8);
\node at (3,2.05) {$\mbox{U}(2p)$};
\node at (3,1.3) {$\mbox{SU}(2p)$};
\draw[-,dashed] (2.9,1.9) to [out=5,in=175] (6.5,1.9);
\draw[-,dashed] (2.9,-1.7) to [out=5,in=175] (6.5,-1.7);

\draw[fill=black,fill opacity =0.1] (6.5,-0.2) ellipse (0.7 and 1);
\draw (6.5,0.1) ellipse (1.3 and 1.8);
\draw (6.7,0.3) ellipse (1.8 and 2.2);

\node at (6.5,2.05) {$\mbox{SO}_\mI(4p)$};
\node at (6.5,2.6) {$\mbox{SO}(4p)$};
\node at (6.5,0.95) {$\mbox{SO}_Q(4p)$};

\draw[->,>=angle 60,dashed] (2.9,-0.5) to (2.3,-2.2);
\node at (2.3,-2.4) {\small subgroup of SU$(2p)$};

\node at (0.85,0.87) {$>$};
\node at (0.85,-1.13){$>$};
\node at (0.85,-0.9){$\psi_p$};
\node at (4.65,0.88) {$>$};
\node at (4.65,-1.11){$>$};
\node at (4.65,-0.9){$\varphi_{2p}$};
\end{tikzpicture}
\end{center}

Quite naturally, the subgroup $\mbox{SO}_Q(4p)$ of $\mbox{SO}(4p)$ has a double covering by $\mbox{Spin}_Q(4p)$, the subgroup of $\mbox{Spin}(4p)$ consisting of the $\mbox{Spin}(4p)$--elements which are commuting with both $s_{\mI}$ and $s_{\mJ}$, where now $s_{\mJ}$ is the $\mbox{Spin}(4p)$--element corresponding to $\mJ_p$. Recall that $s_\mI$, corresponding to the complex structure $\mI_{4p}$, is given by $s_{\mI} = s_1 \cdots s_{2p}$, where $s_j = \frac{\sqrt{2}}{2} \bigl ( 1 + e_{2j-1} e_{2j} \bigr )$, $j=1,\ldots, 2p$, see  (\ref{sI}). Similarly, for $s_{\mJ}$ we find
\begin{equation}
s_{\mJ} = \widetilde{s_1} \cdots \widetilde{s_p}, \quad \widetilde{s_j} = \frac{1}{2} \bigl ( 1 + e_{4j-3} e_{4j-1} \bigr ) \bigl (1 - e_{4j-2} e_{4j} \bigr ), \; j=1,\ldots, p
\label{sJ}
\end{equation}
At the level of the Lie algebras we have the following picture. The real symplectic Lie algebra $\gsp(p)$ of skew--symplectic $M_p(\mH)$--matrices:
$$
\gsp(p) = \bigl \{ A \in \mbox{GL}_p(\mH) : A + A^\ast = 0 \bigr \}
$$
is isomorphic with the subalgebra $\psi_p(\gsp(p))$ of the Lie algebra $\gu(2p)$ of skew--hermitian $M_{2p}(\mC)$--matrices. Moreover, for $A \in \gsp(p)$, $\psi_p(A)$ satisfies the relation
\begin{equation}
\psi_p(A)^T \mI_{2p} + \mI_{2p} \psi_p(A) = 0
\label{vwe}
\end{equation}
On the other hand, there is the complex symplectic Lie group $\mbox{Sp}_{2p}(\mC)$ of complex linear maps preserving the standard skew--hermitian form on $\mC^{2p}$:
$$
\mbox{Sp}_{2p}(\mC) = \bigl \{ A \in \mbox{GL}_{2p}(\mC) : A^T \mI_{2p} A = \mI_{2p} \bigr \}
$$
and its corresponding complex symplectic Lie algebra $\gsp_{2p}(\mC)$ given by
$$
\gsp_{2p}(\mC) = \bigl \{ A \in \mbox{GL}_{2p}(\mC) : A^T \mI_{2p} + \mI_{2p} A = 0 \bigr \} 
$$
which can be decomposed into the direct sum of its hermitian and skew--hermitian subalgebras, both subalgebras being isomorphic under multiplication by the imaginary unit $i$:
$$
\gsp_{2p}(\mC) = \bigl ( \gsp_{2p}(\mC) \cap \gu(2p) \bigr ) \oplus i \bigl ( \gsp_{2p}(\mC) \cap \gu(2p) \bigr )
$$
The skew--hermitian subalgebra of $\gsp_{2p}(\mC)$ is mostly referred to as the {\em compact form} of $\gsp_{2p}(\mC)$. In view of (\ref{vwe}) this leads to the following result.
\begin{proposition}
\label{propositionQ6}
The real symplectic Lie algebra $\gsp(p)$ of skew--symplectic $M_p(\mH)$--matrices is isomorphic with the compact form $\gsp_{2p}(\mC) \cap \gu(2p)$ of the complex symplectic Lie algebra $\gsp_{2p}(\mC)$:
$$
\psi_p (\gsp(p)) = \gsp_{2p}(\mC) \cap \gu(2p)
$$
\end{proposition}

\noindent Schematically, one has
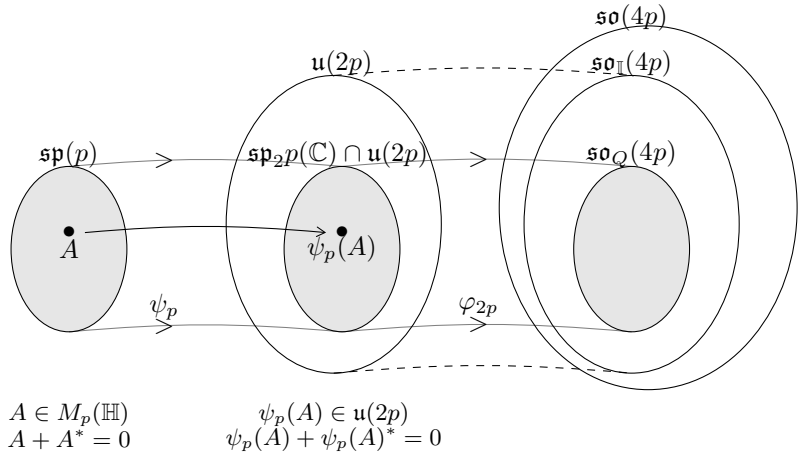
\begin{center}
\begin{tikzpicture}[scale=1.1]
\draw[fill=black,fill opacity =0.1] (-0.3,-0.2) ellipse (0.7 and 1);
\node at (-0.3,0) (bul1) {$\bullet$};
\node at (-0.3,-0.2) {$A$};
\node at (-0.3,0.95) {$\gsp(p)$};
\node at (-0.3,-2.2) {\small $A \in M_p(\mH)$};
\node at (-0.3,-2.5) {\small $A+A^\ast=0$};
\draw[-, draw opacity=0.5] (-0.3,0.8) to [out=5,in=175]  (3,0.8);
\draw[-, draw opacity=0.5] (3,0.8) to [out=5,in=175]  (6.5,0.8);
\draw[-, draw opacity=0.5] (-0.3,-1.2) to [out=5,in=175]  (3,-1.2);
\draw[-, draw opacity=0.5] (3,-1.2) to [out=5,in=175]  (6.5,-1.2);

\draw[fill=black,fill opacity =0.1] (3,-0.2) ellipse (0.7 and 1);
\node at (2.95,0.95) {$\gsp_2p(\mC) \cap \gu(2p)$};
\node at (3,0) (bul2) {$\bullet$};
\node at (3,-0.2) {$\psi_p(A)$};
\node at (2.9,-2.2) {\small $\psi_p(A) \in \gu(2p)$};
\node at (2.9,-2.5) {\small $\psi_p(A) + \psi_p(A)^\ast = 0$}; 
\draw (2.9,0.1) ellipse (1.3 and 1.8);
\node at (3,2.05) {$\gu(2p)$};
\draw[-,dashed] (2.9,1.9) to [out=5,in=175] (6.5,1.9);
\draw[-,dashed] (2.9,-1.7) to [out=5,in=175] (6.5,-1.7);
\draw[->,>=angle 60] (bul1.east) to [out=5,in=175] (bul2.west);

\draw[fill=black,fill opacity =0.1] (6.5,-0.2) ellipse (0.7 and 1);
\draw (6.5,0.1) ellipse (1.3 and 1.8);
\draw (6.7,0.3) ellipse (1.8 and 2.2);

\node at (6.5,2.05) {$\gso_\mI(4p)$};
\node at (6.5,2.6) {$\gso(4p)$};
\node at (6.5,0.95) {$\gso_Q(4p)$};

\node at (0.85,0.87) {$>$};
\node at (0.85,-1.13){$>$};
\node at (0.85,-0.9){$\psi_p$};
\node at (4.65,0.88) {$>$};
\node at (4.65,-1.11){$>$};
\node at (4.65,-0.9){$\varphi_{2p}$};
\end{tikzpicture}
\end{center}

Now let us show how to realize the Lie algebra $\gsp(p)$ inside the Clifford algebra. The Lie algebra $\gu(2p)$ is isomorphic with the subalgebra $\gso_{\mI}(4p)$ of skew--symmetric matrices in $\gso(4p)$ which are commuting with $\log(\mI_{4p}) = \frac{\pi}{2} \mI_{4p}$:
$$
\varphi_{2p} \bigl ( \gu(2p) \bigr ) = \gso_{\mI}(4p)
$$
while $\gsp(p)$ is isomorphic with the subalgebra $\gso_Q(4p)$ of skew--symmetric matrices in $\gso_{\mI}(4p)$ which are moreover also commuting with $\log(\mJ_{4p}) = \frac{\pi}{2} \mJ_{4p}$:
$$
\varphi_{2p} \circ \psi_p \bigl ( \gsp(p) \bigr ) = \gso_Q (4p)
$$
On the other hand we know that $\gu(2p)$ is isomorphic with the Lie algebra $\gspin_{\mI}(4p)$ of bivectors commuting with $\sigma_{\mI} = \frac{\pi}{4} \bigl ( e_1 e_2 + e_3 e_4 + \ldots e_{4p-1} e_{4p} \bigr )$. Let $\sigma_{\mJ}$ be the bivector
$$
\sigma_{\mJ} = \frac{\pi}{4} \bigl ( (e_1e_3 - e_2e_4) + ( e_5e_7 - e_6e_8) + \ldots + (e_{4p-3} e_{4p-1} - e_{4p-2} e_{4p} \bigr )
$$
such that $\exp ( \sigma_{\mJ}) = s_{\mJ}$, where $s_{\mJ}$ is the $\mbox{Spin}(4p)$--element corresponding to the complex structure $\mJ_{4p}$, given by (\ref{sJ}). Then the real symplectic Lie algebra $\gsp(p)$ is isomorphic with the subalgebra $\gspin_Q(4p)$ of bivectors commuting with both $\sigma_{\mI}$ and $\sigma_{\mJ}$:
$$
\varphi_{2p} \circ \psi_p \bigl (\gsp(p) \bigr ) = \gso_Q(4p) \simeq \gspin_Q(4p)
$$
The Lie algebra $\gspin_Q(4p)$, of dimension $p(2p+1)$, has the following basis:
\begin{itemize}
\item $e_{4j-3} e_{4j-2} - e_{4j-1} e_{4j}$, $j=1,\ldots,p$;\\[-7mm]
\item $e_{4j-3} e_{4j-1} + e_{4j-2} e_{4j}$, $j=1,\ldots,p$;\\[-7mm]
\item $e_{4j-3} e_{4j} - e_{4j-2} e_{4j-1}$, $j=1,\ldots,p$;\\[-7mm]
\item $e_{4j-3} e_{4k-3} + e_{4j-2} e_{4k-2} + e_{4j-1} e_{4k-1} + e_{4j} e_{4k}$ $j < k=1,\ldots,p$;\\[-7mm]
\item $e_{4j-3} e_{4k-2} - e_{4j-2} e_{4k-3} - e_{4j-1} e_{4k} + e_{4j} e_{4k-1}$ $j < k=1,\ldots,p$;\\[-7mm]
\item $e_{4j-3} e_{4k-1} + e_{4j-2} e_{4k} - e_{4j-1} e_{4k-3} - e_{4j} e_{4k-2}$ $j < k=1,\ldots,p$;\\[-7mm]
\item $e_{4j-3} e_{4k} - e_{4j-2} e_{4k-1} + e_{4j-1} e_{4k-2} - e_{4j} e_{4k-3}$ $j < k=1,\ldots,p$.
\end{itemize}
Schematically, one has
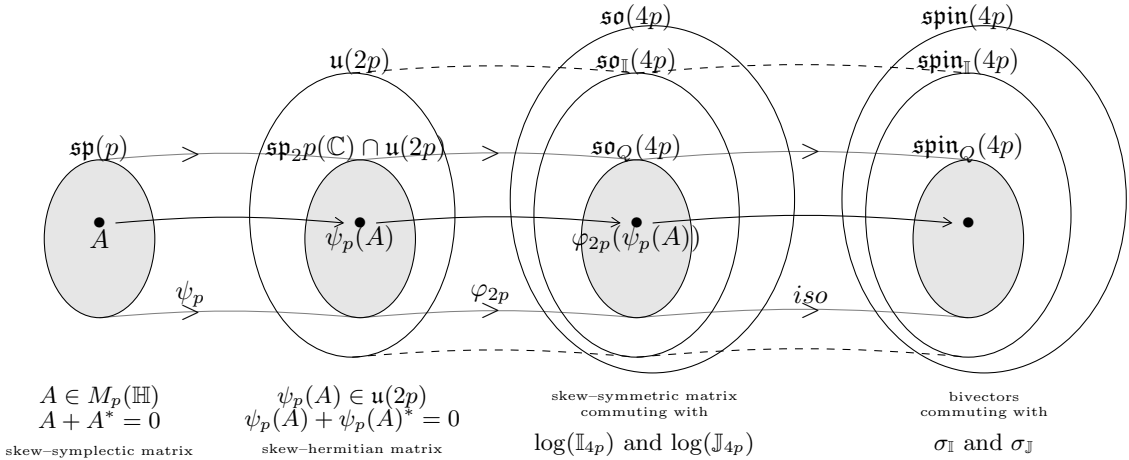
\begin{center}
\hspace*{-4mm}
\begin{tikzpicture}[scale=1.05]
\draw[fill=black,fill opacity =0.1] (-0.3,-0.2) ellipse (0.7 and 1);
\node at (-0.3,0) (bul1) {$\bullet$};
\node at (-0.3,-0.2) {$A$};
\node at (-0.3,0.95) {$\gsp(p)$};
\node at (-0.3,-2.2) {\small $A \in M_p(\mH)$};
\node at (-0.3,-2.5) {\small $A+A^\ast=0$};
\node at (-0.3,-2.9) {\tiny skew--symplectic matrix};
\draw[-, draw opacity=0.5] (-0.3,0.8) to [out=5,in=175]  (3,0.8);
\draw[-, draw opacity=0.5] (3,0.8) to [out=5,in=175]  (6.5,0.8);
\draw[-, draw opacity=0.5] (6.5,0.8) to [out=5,in=175]  (10.7,0.8);
\draw[-, draw opacity=0.5] (-0.3,-1.2) to [out=5,in=175]  (3,-1.2);
\draw[-, draw opacity=0.5] (3,-1.2) to [out=5,in=175]  (6.5,-1.2);
\draw[-, draw opacity=0.5] (6.5,-1.2) to [out=5,in=175]  (10.7,-1.2);

\draw[fill=black,fill opacity =0.1] (3,-0.2) ellipse (0.7 and 1);
\node at (2.95,0.95) {$\gsp_2p(\mC) \cap \gu(2p)$};
\node at (3,0) (bul2) {$\bullet$};
\node at (3,-0.2) {$\psi_p(A)$};
\node at (2.9,-2.2) {\small $\psi_p(A) \in \gu(2p)$};
\node at (2.9,-2.5) {\small $\psi_p(A) + \psi_p(A)^\ast = 0$}; 
\node at (2.9,-2.85) {\tiny skew--hermitian matrix};
\draw (2.9,0.1) ellipse (1.3 and 1.8);
\node at (3,2.05) {$\gu(2p)$};
\draw[-,dashed] (2.9,1.9) to [out=5,in=175] (6.5,1.9);
\draw[-,dashed] (2.9,-1.7) to [out=5,in=175] (6.5,-1.7);
\draw[->,>=angle 60] (bul1.east) to [out=5,in=175] (bul2.west);

\draw[fill=black,fill opacity =0.1] (6.5,-0.2) ellipse (0.7 and 1);
\draw (6.5,0.1) ellipse (1.3 and 1.8);
\draw (6.7,0.3) ellipse (1.8 and 2.2);
\node at (6.5,0) (bul3) {$\bullet$};
\node at (6.5,-0.2) {$\varphi_{2p}(\psi_p(A))$};
\draw[->,>=angle 60] (bul2.east) to [out=5,in=175] (bul3.west);
\node at (6.5,2.05) {$\gso_\mI(4p)$};
\node at (6.5,2.6) {$\gso(4p)$};
\node at (6.5,0.95) {$\gso_Q(4p)$};
\node at (6.6,-2.2) {\tiny skew--symmetric matrix};
\node at (6.6,-2.4) {\tiny commuting with};
\node at (6.6,-2.8) {\small $\log(\mI_{4p})$ and $\log(\mJ_{4p})$};

\draw[fill=black,fill opacity =0.1] (10.7,-0.2) ellipse (0.7 and 1);
\draw (10.7,0.1) ellipse (1.3 and 1.8);
\draw (10.9,0.3) ellipse (1.8 and 2.2);
\node at (10.7,0) (bul4) {$\bullet$};
\draw[->,>=angle 60] (bul3.east) to [out=5,in=175] (bul4.west);
\draw[-,dashed] (6.5,1.9) to [out=5,in=175] (10.7,1.9);
\draw[-,dashed] (6.5,-1.7) to [out=5,in=175] (10.7,-1.7);
\node at (10.7,2.05) {$\gspin_\mI(4p)$};
\node at (10.7,2.6) {$\gspin(4p)$};
\node at (10.7,0.95) {$\gspin_Q(4p)$};
\node at (10.9,-2.2) {\tiny bivectors};
\node at (10.9,-2.4) {\tiny commuting with};
\node at (10.9,-2.8) {\small $\sigma_{\mI}$ and $\sigma_{\mJ}$};

\node at (0.85,0.87) {$>$};
\node at (0.85,-1.13){$>$};
\node at (0.85,-0.9){$\psi_p$};
\node at (4.65,0.88) {$>$};
\node at (4.65,-1.11){$>$};
\node at (4.65,-0.9){$\varphi_{2p}$};
\node at (8.7,0.91) {$>$};
\node at (8.7,-1.11){$>$};
\node at (8.7,-0.9){$iso$};
\end{tikzpicture}
\end{center}

\noindent {\bf Example} ($p=1$) \\
Consider the $\mbox{Spin}(4)$--element
$$
s_A = \frac{1}{2} (1 + e_1 e_4)(1-e_2e_3)
$$
It can be readily verified that $s_A$ commutes with 
$$
s_{\mI} = \frac{1}{2} ( 1 + e_1 e_2) ( 1 + e_3 e_4 ) \quad \mbox{\ and \ } \quad s_{\mJ} = \frac{1}{2} ( 1+ e_1 e_3)(1 - e_2 e_4)
$$
so $s_A$ belongs to $\mbox{Spin}_Q(4)$. Its action on the vector $\uX = X_1 e_1 + X_2 e_2 + X_3 e_3 + X_4 e_4$ is given by $s_A \uX s_A^{-1} = -X_4 e_1 + X_3 e_2 - X_2 e_3 + X_1 e_4$, which is also obtained by the action of the corresponding $\mbox{SO}_Q(4)$--matrix
$$
A = \left ( \begin{array}{cccc} & & & 1 \\ & & -1 & \\  & 1 & & \\ -1 & & & \end{array} \right )
$$
which indeed commutes with $\mI_4$ and $\mJ_4$. We expect $\varphi_2^{-1}(A)$ to belong to $\mbox{SU}(2)$ and $\psi_1^{-1} \bigl ( \varphi_2^{-1}(A) \bigr )$ to belong to $\mbox{Sp}(1)$. We have indeed
\begin{eqnarray*}
B &=& \varphi_2^{-1} ( A) = \left ( \begin{array}{cc} & i \\ i & \end{array} \right ) \in \mbox{SU}(2) \\
C &=& \psi_1^{-1} \bigl ( \varphi_2^{-1}(A) \bigr ) = ( \; k \; ) \in \mbox{Sp}(1)
\end{eqnarray*}
At the level of the Lie algebras we find that for $c = \bigl ( \frac{\pi}{2} k \bigr )  \in \gsp(1)$ there holds that $\exp(c)=C$. Now we first put
$$
b = \psi_1 ( c) = \frac{\pi}{2} \left ( \begin{array}{cc} & i \\ i & \end{array} \right ) 
$$
and it may be checked that indeed $\psi_1(c) \in \gsp_2(\mC) \cap \gu(2)$, and moreover $\exp(b) = B$. Next we let
$$
a = \varphi(b) = \frac{\pi}{2}  \left ( \begin{array}{cccc} & & & 1 \\ & & -1 & \\  & 1 & & \\ -1 & & & \end{array} \right )
$$
which belongs to $\gso_Q(4)$ and for which $\exp(a)=A$. Finally, with the matrix $a$ there corresponds the bivector
$$
\sigma_A = \frac{\pi}{4} ( e_1 e_4 - e_2 e_3)
$$
which belongs to $\gspin_Q(4)$  and for which there holds $\exp(\sigma_A) = \frac{1}{2}(1+e_1e_4)(1-e_2e_3) = s_A$. This leads to the following scheme:
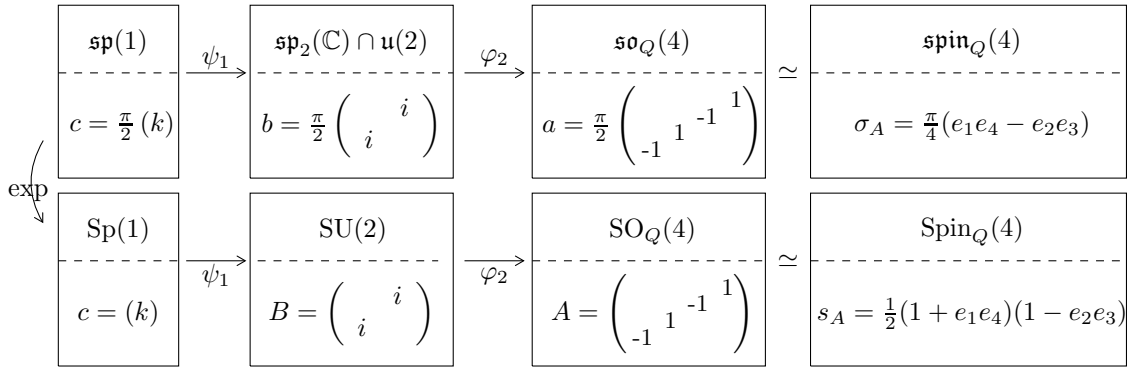
\begin{center}
\begin{tikzpicture}
\node at (-0.9,0) {$\gsp(1)$};
\node at (-0.85,-1.1) {$c=\frac{\pi}{2} \left ( k \right )$ };
\draw[dashed] (-1.7,-0.4) -- (-0.1,-0.4);
\draw (-1.7,-1.8) -- (-1.7,0.5) -- (-0.1,0.5) -- (-0.1,-1.8) -- (-1.7,-1.8);

\draw[->,>=angle 60] (0,-0.4) to [out=0,in=180] (0.8,-0.4);
\node at (0.4,-0.2) {$\psi_1$};

\node at (2.2,0) {$\gsp_2(\mC) \cap \gu(2)$};
\node at (2.2,-1.1) {$b=\frac{\pi}{2} \left ( \begin{array}{cc}  & i \\ i &  \end{array} \right )$ };
\draw[dashed] (0.95,-0.4) -- (3.55,-0.4);
\draw (0.85,-1.8) -- (0.85,0.5) -- (3.55,0.5) -- (3.55,-1.8) -- (0.85,-1.8);

\draw[->,>=angle 60] (3.7,-0.4) to [out=0,in=180] (4.5,-0.4);
\node at (4.1,-0.2) {$\varphi_2$};

\node at (6.2,0) {$\gso_Q(4)$};
\node at (6.2,-1.1) {$a = \frac{\pi}{2} \left ( \begin{array}{cccc} \hspace*{-2mm} &  \hspace*{-2mm}  & \hspace*{-2mm} & \hspace*{-2mm} \mbox{\small 1} \hspace*{-2mm} \\[-2mm] \hspace*{-2mm}  & \hspace*{-2mm} & \hspace*{-2mm} \mbox{\small -1} & \hspace*{-2mm} \\[-2mm]  \hspace*{-2mm} & \hspace*{-2mm} \mbox{\small 1} & \hspace*{-2mm}  & \hspace*{-2mm} \\[-2mm] \hspace*{-1mm} \mbox{\small -1} & \hspace*{-2mm} & \hspace*{-2mm}  & \hspace*{-2mm} \end{array} \right )$} ;
\draw[dashed] (4.6,-0.4) -- (7.7,-0.4);
\draw (4.6,-1.8) -- (4.6,0.5) -- (7.7,0.5) -- (7.7,-1.8) -- (4.6,-1.8);

\node at (8.0,-0.4) {$\simeq$};

\node at (10.45,0) {$\gspin_Q(4)$};
\node at (10.45,-1.1) {$\sigma_A = \frac{\pi}{4} (e_1 e_4 - e_2 e_3)$};
\draw[dashed] (8.3,-0.4) -- (12.5,-0.4);
\draw (8.3,-1.8) -- (8.3,0.5) -- (12.5,0.5) -- (12.5,-1.8) -- (8.3,-1.8);

\draw[->,>=angle 60] (-1.9,-1.3) to [out=240,in=120]  (-1.9,-2.4);
\node at (-2.1,-2) {$\exp$};

\node at (-0.9,-2.5) {$\mbox{Sp}(1)$};
\node at (-0.9,-3.6) {$c= \left ( k \right )$ };
\draw[dashed] (-1.7,-2.9) -- (-0.1,-2.9);
\draw (-1.7,-4.3) -- (-1.7,-2) -- (-0.1,-2) -- (-0.1,-4.3) -- (-1.7,-4.3);

\draw[->,>=angle 60] (0,-2.9) to [out=0,in=180] (0.8,-2.9);
\node at (0.4,-3.1) {$\psi_1$};

\node at (2.22,-2.5) {$ \mbox{SU}(2)$};
\node at (2.2,-3.6) {$B= \left ( \begin{array}{cc}  & i \\ i &  \end{array} \right )$ };
\draw[dashed] (0.95,-2.9) -- (3.35,-2.9);
\draw (0.85,-4.3) -- (0.85,-2) -- (3.55,-2) -- (3.55,-4.3) -- (0.85,-4.3);

\draw[->,>=angle 60] (3.7,-2.9) to [out=0,in=180] (4.5,-2.9);
\node at (4.1,-3.1) {$\varphi_2$};

\node at (6.2,-2.5) {$\mbox{SO}_Q(4)$};
\node at (6.2,-3.6) {$A = \left ( \begin{array}{cccc} \hspace*{-2mm} &  \hspace*{-2mm}  & \hspace*{-2mm} & \hspace*{-2mm} \mbox{\small 1} \hspace*{-2mm} \\[-2mm] \hspace*{-2mm}  & \hspace*{-2mm} & \hspace*{-2mm} \mbox{\small -1} & \hspace*{-2mm} \\[-2mm]  \hspace*{-2mm} & \hspace*{-2mm} \mbox{\small 1} & \hspace*{-2mm}  & \hspace*{-2mm} \\[-2mm] \hspace*{-1mm} \mbox{\small -1} & \hspace*{-2mm} & \hspace*{-2mm}  & \hspace*{-2mm} \end{array} \right ) $};
\draw[dashed] (4.6,-2.9) -- (7.7,-2.9);
\draw (4.6,-4.3) -- (4.6,-2) -- (7.7,-2) -- (7.7,-4.3) -- (4.6,-4.3);

\node at (8.0,-2.9) {$\simeq$};

\node at (10.45,-2.5) {$\mbox{Spin}_Q(4)$};
\node at (10.45,-3.6) {$s_A = \frac{1}{2} (1+e_1 e_4)(1 - e_2 e_3)$};
\draw[dashed] (8.3,-2.9) -- (12.5,-2.9);
\draw (8.3,-4.3) -- (8.3,-2) -- (12.5,-2) -- (12.5,-4.3) -- (8.3,-4.3);
\end{tikzpicture}
\end{center}


\section{Spinor space}


\subsection{Homogeneous spinor spaces}
\label{homospin}


A Clifford algebra may be decomposed as a direct sum of isomorphic copies of a spinor space $\mS$, which, abstractly, may be defined as a minimal left ideal in the Clifford algebra. A spinor space is an irreducible Spin group representation, and may be realized in the following way.\\[-2mm]

The complex structure $\mI_{4p}$ acts upon the basis vectors as follows:
\begin{align*}
e_{2k-1} =  ( 0 \cdots 1 \cdots 0 ) & \hspace*{5mm} \longmapsto & \hspace*{-15mm} \mI_{4p} [e_{2k-1} ]  =  ( 0 \cdots \phantom{(-}1\phantom{)} \cdots 0)  =  \phantom{-} e_{2k\phantom{-1}} \\[-1mm]
 \uparrow \hspace*{8.5mm} && \uparrow \hspace*{25mm} \\[-3mm]
\mbox{\tiny $(2k-1)$th place} &&  \mbox{\tiny $(2k)$th place} \hspace*{17mm}\\[2mm]
e_{2k}  =   ( 0 \cdots 1 \cdots 0 ) & \hspace*{5mm} \longmapsto & \hspace*{-15mm} \mI_{4p} [e_{2k} ]  =  ( 0 \cdots (-1) \cdots 0)  =  - e_{2k-1} \\[-1mm]
 \uparrow \hspace*{8.5mm} & & \uparrow \hspace*{25mm}\\[-3mm]
\mbox{\tiny $(2k)$th place} && \mbox{\tiny $(2k-1)$th place} \hspace*{17mm}
\end{align*}
It entails two projection operators $\frac{1}{2} ( \mathbf{1} \pm i \mI_{4p} )$ on the complexification $\mC^{4p}$ of $\mR^{4p}$, for which, with $k=1,\ldots, 2p$:
\begin{eqnarray*}
\phantom{-} \frac{1}{2} (\mathbf{1} + i \mI_{4p} ) [e_{2k-1}] &=& \phantom{-} \frac{1}{2} ( e_{2k-1} + i e_{2k} ) = \phantom{-i}\gf_k^\dagger \\
\phantom{-} \frac{1}{2} ( \mathbf{1} + i \mI_{4p} ) [e_{2k}] &=& \phantom{-} \frac{1}{2} ( e_{2k} - i e_{2k-1} ) = -i \gf_k^\dagger \\
- \frac{1}{2} ( \mathbf{1} - i \mI_{4p} ) [ e_{2k-1}] & = &  - \frac{1}{2} ( e_{2k-1} - i e_{2k} ) = \phantom{-i}\gf_k \\
 - \frac{1}{2} ( \mathbf{1} - i \mI_{4p} ) [ e_{2k}] & = &  - \frac{1}{2} ( e_{2k} + i e_{2k-1} ) = \phantom{-}i \gf_k  
\end{eqnarray*} 
leading to the so-called Witt basis $(\gf_1,\gf_1^\dagger,\ldots,\gf_{2p},\gf_{2p}^\dagger)$ of $\mC^{4p}$. The Witt basis vectors satisfy the Grassmann identities
$$
\gf_j \gf_k + \gf_k \gf_j = 0 , \quad \gf_j^\dagger \gf_k^\dagger + \gf_k^\dagger \gf_j^\dagger = 0, \qquad j,k=1,\ldots,2p
$$
including their isotropy
$$
\gf_j^2 = (\gf_j^\dagger)^2 = 0, \qquad j=1,\ldots, 2p
$$
and the duality identities
$$
\gf_j \gf_k^\dagger + \gf_k^\dagger \gf_j = \delta_{jk}, \qquad j,k=1,\ldots,2p
$$
The Witt basis vectors $(\gf_1,\ldots,\gf_{2p})$ on the one hand, and $(\gf_1^\dagger,\ldots,\gf_{2p}^\dagger)$ on the other, respectively span isotropic subspaces $W$ and $W^\dagger$ of $\mC^{4p}$, such that 
$$
\mC^{4p} = W \oplus W^\dagger
$$
those subspaces being eigenspaces of the complex structure $\mI_{4p}$ with respective eigenvalues $-i$ and $i$. They also generate the respective Grassmann algebras $\mC \bigwedge_{2p}$ and $\mC \bigwedge_{2p}^\dagger$. Note that the $\cdot^\dagger$--notation corresponds to the hermitian conjugation in the Clifford algebra $\mC_{4p}$ (see Section 2).\\[-2mm]

With the self-adjoint idempotents 
$$
I_j = \gf_j \gf_j^\dagger = \frac{1}{2} ( 1- i e_{2j-1} e_{2j} ), \qquad j=1,\ldots,2p
$$
we compose the primitive self--adjoint idempotent $I = I_1 I_2 \cdots I_{2p}$ leading to the realization of the spinor space $\mS$ as $\mS = \mC_{4p} I \simeq \mC_{2p} I$. Since $\gf_j I =0$, $j=1,\ldots,2p$, we also have $\mS \simeq \mC \bigwedge_{2p}^\dagger I$.\\[-2mm]

When decomposing the Grassmann algebra $\mC \bigwedge_{2p}^\dagger$ into its so--called homogeneous parts:
$$
\mC \bigwedge_{2p}\nolimits^\dagger = \bigoplus_{r=0}^{2p} \left ( \mC \bigwedge_{2p}\nolimits^\dagger \right )^{(r)}
$$
where $\left (\mC \bigwedge_{2p}^\dagger \right )^{(r)}$ is spanned by all products of $r$ Witt basis vectors out of $(\gf_1^\dagger,\ldots,\gf_{2p}^\dagger)$, the spinor space $\mS$ accordingly decomposes into
\begin{equation}
\mS = \bigoplus_{r=0}^{2p} \mS^r, \qquad \mbox{with \ } \mS^r \simeq \left ( \mC \bigwedge_{2p}\nolimits^\dagger \right)^{(r)} I
\label{decompspin}
\end{equation}
where the dimension of each of the homogeneous parts is given by
$$
\mbox{dim} \ \mS^r = \binom{2p}{r} = \frac{(2p)!}{r! (2p-r)!}
$$
These homogeneous parts $\mS^r$, $r=0,\ldots,2p$, of spinor space provide models for fundamental $\mbox{U}(2p)$-representations (see \cite{howe}) and for fundamental $\gsl_{2p}(\mC)$-representations (see \cite{partI}, \cite{Ee}). Our aim now is to decompose each such homogeneous spinor space $\mS^r$ into fundamental $\mbox{Sp}(p)$-representations. We have already seen (see Section 3) that the associated Lie algebra $\gsp( p)$ is isomorphic to the compact form
$\gsp_{2p}(\mC) \cap \gu(2p)$ of the Lie algebra $\gsp_{2p}(\mC)$. This complex Lie algebra $\gsp_{2p}(\mC)$ in its turn is a subalgebra of the Lie algebra $\gsl_{2p}(\mC)$ of traceless matrices, to which we turn our attention now. \\[-2mm]


\subsection{The Lie algebras $\gsl_{2p}(\mC)$, $\gsp_{2p}(\mC)$ and $\gsl_{p}(\mC)$}
\label{subsecgsp}


The Lie algebra $\gsl_{2p}(\mC)$ can be realized as a subalgebra of the Lie algebra $\mC_{4p}^{(2)}$ of bivectors in the complex Clifford algebra $\mC_{4p}$. The correspondence between bivectors in $\mC_{4p}^{(2)}$ and traceless matrices in $\gsl_{2p}(\mC)$ is the usual one (see also Section 3):
$$
b \in \mC_{4p}^{(2)} \longleftrightarrow B \in M_{4p}(\mR) \longleftrightarrow \varphi^{-1}_{2p}(B) \in \gsl_{2p}(\mC)
$$
where the respective correspondences are established by
\begin{eqnarray*}
&& \hspace*{6mm} \mbox{\small $(i,j)$ entry} \\[-2mm]
&& \hspace*{12.8mm} \downarrow \\[-5mm]
\frac{1}{2} e_{ij} & \longleftrightarrow & \left ( \begin{array}{ccc} && \hspace*{-2mm}\\ & 1 & \hspace*{-2mm}\\[-1mm] && \hspace*{-2mm}\\[-1mm]  -1 & & \hspace*{-2mm}\\ & & \hspace*{-2mm} \end{array} \right )  \qquad \mbox{and} \qquad \left ( \begin{array} {rc} a & b \\ -b & a \end{array} \right ) \longleftrightarrow a+ib \\[-5mm]
&&  \hspace*{7.5mm} \uparrow \\[-2mm]
&& \hspace*{1mm} \mbox{\small $(j,i)$ entry} 
\end{eqnarray*}
The Lie algebra $\mC_{4p}^{(2)}$ has real dimension $4p(4p-1)$. A basis is given by $e_j e_k$, $j<k=1,\ldots,4p$, or, alternatively by
\begin{itemize}
\item $\gf_j \gfd_j - \gfd_j \gf_j$, $j=1,\ldots,2p$;\\[-6mm]
\item $\gf_j \gf_k$, $\gfd_j \gf_k$, $\gf_j \gfd_k$, $\gfd_j \gfd_k$, $j<k=1,\ldots,2p$.
\end{itemize}
The Lie algebra $\gsl_{2p}(\mC)$ has real dimension $2(4 p^2-1)$. A basis is given by
\begin{itemize}
\item $e_{2j-1} e_{2j} - e_{2j+1} e_{2j+2}$, $j=1,\ldots,2p-1$, $j=1,\ldots,2p-1$;\\[-6mm]
\item $e_{2j-1} e_{2k-1} + e_{2j} e_{2k}$, $j<k=1,\ldots,2p$; $j<k=1,\ldots,2p$;\\[-6mm]
\item $e_{2j-1} e_{2k} - e_{2j} e_{2k-1}$, $j<k=1,\ldots,2p$,
\end{itemize}
or, alternatively, by
\begin{itemize}
\item $\gf_{j+1} \gfd_{j+1} - \gf_j \gfd_j$, $j=1,\ldots,2p-1$;\\[-5mm]
\item $\gfd_j \gf_k + \gf_j \gfd_k$, $j<k=1,\ldots,2p$;\\[-5mm]
\item $\gfd_j \gf_k - \gf_j \gfd_k$, $j<k=1,\ldots,2p$.
\end{itemize}
Notice that the set $\bigl \{ \gf_{j+1} \gfd_{j+1} - \gf_j \gfd_j : j=1,\ldots,2p-1 \bigr \}$ spans the Cartan subalgebra $\gh$ of $\gsl_{2p}(\mC)$, or, alternatively,
$$
\gh = \mbox{span} \bigl \{ H_j  = I_j - I_{2p} : j=1,2,\ldots,2p-1 \bigr \}
$$
To give an idea how the associated $\gsl_{2p}(\mC)$--matrices look like, we consider, for $p=2$, the bivector $\frac{1}{2} ( e_3 e_4-e_5e_6)$. Its matrix representation is
$$
\left ( \begin{array}{cccc}
\fbox{$\begin{array}{cc}   \hspace*{-2mm} \phantom{-1} & \phantom{1} \hspace*{-3mm}  \\ \hspace*{-2mm} \phantom{-1} &    \phantom{1} \hspace*{-3mm}\end{array}$} &  & & \\
& \hspace*{-2mm}  \fbox{$\begin{array}{cc}  \hspace*{-2mm} \phantom{-1} & 1 \hspace*{-3mm}  \\ \hspace*{-2mm} -1 &  \phantom{1} \hspace*{-3mm} \end{array}$ } & & \\
& & \hspace*{-2mm}  \fbox{$\begin{array}{cc}  \hspace*{-2mm} \phantom{1}  & -1 \hspace*{-2.5mm} \\ \hspace*{-2mm}  1 &  \phantom{-1} \hspace*{-3mm} \end{array}$} & \\
& & & \hspace*{-2mm} \fbox{$\begin{array}{cc}  \hspace*{-2mm} \phantom{1}  &  \phantom{-1} \hspace*{-3mm}  \\ \hspace*{-2mm} \phantom{1} &   \phantom{-1} \hspace*{-3mm}  \end{array}$}
\end{array} \right )
$$
which, under the isomorphism $\varphi_4$, corresponds with the matrix
$$
\left ( \begin{array}{cccc} & & & \\ & i & & \\ & & -i & \\ & & & \end{array} \right )
$$
indeed belonging to $\gsl_4(\mC)$. As a second example, still with $p=2$, we consider the bivector $\frac{1}{2} ( e_5 e_6-e_7e_8)$. Here we have
$$
\frac{1}{2} ( e_5 e_6-e_7e_8) \longleftrightarrow \left ( \begin{array}{cccc}
\fbox{$\begin{array}{cc}  \hspace*{-2mm} \phantom{-1} & \phantom{1} \hspace*{-3mm}  \\ \hspace*{-2mm} \phantom{-1} &  \phantom{1} \hspace*{-3mm} \end{array}$}  &  & & \\
& \hspace*{-2mm} \fbox{$\begin{array}{cc}   \hspace*{-2mm} \phantom{-1} & \phantom{1} \hspace*{-3mm}  \\ \hspace*{-2mm} \phantom{-1} &  \phantom{1} \hspace*{-3mm} \end{array}$}  & & \\
& & \hspace*{-2mm} \fbox{$\begin{array}{cc}  \hspace*{-2mm} \phantom{-1} & 1 \hspace*{-3mm}  \\ \hspace*{-2mm} -1 &  \phantom{1} \hspace*{-3mm} \end{array}$ } &\\
& & & \hspace*{-2mm}  \fbox{$\begin{array}{cc}  \hspace*{-2mm} \phantom{1}  & -1 \hspace*{-2.5mm} \\ \hspace*{-2mm}  1 &  \phantom{-1} \hspace*{-3mm} \end{array}$} 
\end{array} \right ) \longleftrightarrow \left ( \begin{array}{cccc} & & & \\ &  & & \\ & & i & \\ & & & -i \end{array} \right ) \in \gsl_4(\mC)
$$

In its turn the Lie algebra $\gsp_{2p}(\mC)$, being a subalgebra of $\gsl_{2p}(\mC)$, also may be realized as a subalgebra of bivectors in $\mC_{4p}^{(2)}$.\\

A first way of establishing this embedding is provided by the general form of an $\gsp_{2p}(\mC)$--matrix, which looks as follows (to fix the ideas we take $p=3$):
\begin{equation}
\left ( \begin{array}{rr|rr|rr} z_{11} & z_{12} & z_{13} & z_{14} & z_{15} & z_{16} \\
z_{21} & - z_{11} & z_{23} & z_{24} & z_{25} & z_{26} \\ \hline
-z_{24} & z_{14} & z_{33} & z_{34} & z_{35} & z_{36} \\
z_{23} & - z_{13} & z_{43} & - z_{33} & z_{45} & z_{46} \\ \hline
-z_{26} & z_{16} & -z_{46} & z_{36} & z_{55} & z_{56} \\
z_{25} & - z_{15} & z_{45} & - z_{35} & z_{65} & - z_{55}
 \end{array} \right ) \in \gsp_6(\mC)
\label{gsp6}
\end{equation}
So it is obvious that this matrix belongs to $\gsl_{2p}(\mC)$. Another way consists in determining a basis for $\gsp_{2p}(\mC)$ which has real dimension $2p(2p+1)$:
\begin{itemize}
\item $H_j^{{\rm sympl}} = H_{2j} - H_{2j-1} = \gf_{2j} \gfd_{2j} - \gf_{2j-1} \gfd_{2j-1}$, $j=1,\ldots,p$;\\[-5mm]
\item $\gfd_{2j-1} \gf_{2j}$, $\gfd_{2j} \gf_{2j-1}$, $j=1,\ldots,p$;\\[-5mm]
\item $\gfd_{2j} \gf_{2k} + \gf_{2j-1} \gfd_{2k-1}$, $\gfd_{2k} \gf_{2j} + \gf_{2k-1} \gfd_{2j-1}$, $j<k=1,\ldots,p$;\\[-5mm]
\item $\gf_{2j} \gfd_{2k-1} - \gfd_{2j-1} \gf_{2k}$, $\gf_{2k-1} \gfd_{2j} - \gfd_{2k} \gf_{2j-1}$, $j<k=1,\ldots,p$.
\end{itemize}
Note that
$$
\gh^{{\rm sympl}} = \bigl \{ H_j^{{\rm sympl}} = \gf_{2j} \gfd_{2j} - \gf_{2j-1} \gfd_{2j-1} : j=1,\ldots,p \bigr \} 
$$
forms the Cartan subalgebra $\gh^{{\rm sympl}}$ of $\gsp_{2p}(\mC)$.\\[-2mm]

As for the examples considered above:
\begin{itemize}
\item $\frac{1}{2} (e_3e_4-e_5e_6)$: its matrix is not of the form corresponding to (\ref{gsp6}) for $p=2$, so this bivector belongs to $\gsl_4(\mC)$ but not to $\gsp_4(\mC)$;
\item $\frac{1}{2} (e_5e_6-e_7e_8)$: its matrix does have the form (\ref{gsp6}), so this bivector belongs to $\gsp_4(\mC) \subset \gsl_4(\mC)$.\\[-2mm]
\end{itemize}

\noindent {\bf Example: explicit bases for $\gsp_2(\mC)$ and $\gsl_2(\mC)$}\\
The space of bivectors $\mC_4^{(2)}$ has complex dimension 6 and its basis reads
$$
e_1 e_2, \; e_1 e_3, \; e_1 e_4, \; e_2 e_3, \; e_2 e_4, \; e_3 e_4
$$
or, in terms of the Witt basis, 
$$
\gf_1 \gfd_1 - \gfd_1 \gf_1, \; \gf_2 \gfd_2 - \gfd_2 \gf_2, \; \gf_1 \gf_2, \; \gfd_1 \gf_2, \; \gf_1 \gfd_2, \; \gfd_1 \gfd_2 
$$
while $\gsl_2(\mC) = \gsp_2(\mC)$ has complex dimension 3 and has the bases
$$
e_1 e_2 - e_3 e_4, \; e_1 e_3 + e_2 e_4, \; e_1 e_4 - e_2 e_3 \quad 
\mbox{\ or\ } \quad \gf_2 \gfd_2 - \gf_1 \gfd_1, \; \gfd_1 \gf_2 + \gf_1 \gfd_2, \; \gfd_1 \gf_2 - \gf_1 \gfd_2
$$

\noindent {\bf Example: explicit bases for $\gsp_4(\mC)$ and $\gsl_4(\mC)$}\\
The space of bivectors $\mC_8^{(2)}$ has complex dimension 28; its basis reads $\{e_1 e_2, e_1 e_3, \ldots, e_7e_8\}$, or
\begin{itemize}
\item $\gf_1 \gfd_1 - \gfd_1 \gf_1$, $\gf_2 \gfd_2 - \gfd_2 \gf_2$, $\gf_3 \gfd_3 - \gfd_3 \gf_3$, $\gf_4 \gfd_4 - \gfd_4 \gf_4$;\\[-6mm]
\item $\gf_1 \gf_2$, $\gf_1 \gf_3$, $\gf_1 \gf_4$, $\gf_2 \gf_3$, $\gf_2 \gf_4$, $\gf_3 \gf_4$; \\[-6mm]
\item $\gfd_1 \gf_2$, $\gfd_1 \gf_3$, $\gfd_1 \gf_4$, $\gfd_2 \gf_3$, $\gfd_2 \gf_4$, $\gfd_3 \gf_4$;\\[-6mm]
\item $\gf_1 \gfd_2$, $\gf_1 \gfd_3$, $\gf_1 \gfd_4$, $\gf_2 \gfd_3$, $\gf_2 \gfd_4$, $\gf_3 \gfd_4$;\\[-6mm]
\item $\gfd_1 \gfd_2$, $\gfd_1 \gfd_3$, $\gfd_1 \gfd_4$, $\gfd_2 \gfd_3$, $\gfd_2 \gfd_4$, $\gfd_3 \gfd_4$.
\end{itemize}
The Lie algebra $\gsl_4(\mC)$ has complex dimension 15 and its basis reads
\begin{itemize}
\item $e_1e_2-e_3e_4$, $e_3e_4-e_5e_6$, $e_5e_6-e_7e_8$;\\[-7mm]
\item $e_1 e_3 + e_2 e_4$, $e_1 e_5 + e_2 e_6$ , $e_1 e_7 + e_2 e_8$, $e_3 e_5 + e_4 e_6$, $e_3 e_7 + e_4 e_8$, $e_5 e_7 + e_6 e_8$;\\[-7mm]
\item $e_1 e_4 - e_2 e_3$, $e_1 e_6 - e_2 e_5$ , $e_1 e_8 - e_2 e_7$, $e_3 e_6 - e_4 e_5$, $e_3 e_8 - e_4 e_7$, $e_5 e_8 - e_6 e_7$
\end{itemize}
or
\begin{itemize}
\item $\gf_2 \gfd_2 - \gf_1 \gfd_1$, $\gf_3 \gfd_3 - \gf_2 \gfd_2$, $\gf_4 \gfd_4 - \gf_3 \gfd_3$;\\[-5mm]
\item $\gfd_1 \gf_2 + \gf_1 \gfd_2$, $\gfd_1 \gf_3 + \gf_1 \gfd_3$, $\gfd_1 \gf_4 + \gf_1 \gfd_4$, $\gfd_2 \gf_3 + \gf_2 \gfd_3$, $\gfd_2 \gf_4 + \gf_2 \gfd_4$, $\gfd_3 \gf_4 + \gf_3 \gfd_4$;\\[-5mm]
\item $\gfd_1 \gf_2 - \gf_1 \gfd_2$, $\gfd_1 \gf_3 - \gf_1 \gfd_3$, $\gfd_1 \gf_4 - \gf_1 \gfd_4$, $\gfd_2 \gf_3 - \gf_2 \gfd_3$, $\gfd_2 \gf_4 - \gf_2 \gfd_4$, $\gfd_3 \gf_4 - \gf_3 \gfd_4$
\end{itemize}
while $\gsp_4(\mC)$ has complex dimension 10 and has the basis
\begin{itemize}
\item $e_1e_2-e_3e_4$, $e_5e_6-e_7e_8$, $e_1e_3+e_2e_4$, $e_5e_7+e_6e_8$, $e_1e_4-e_2e_3$, $e_5e_8-e_6e_7$;\\[-7mm]
\item $e_1e_5 + e_2e_6+e_3e_7+e_4e_8$, $e_1e_7+e_2e_8-e_3e_5-e_4e_6$, $e_1e_6-e_2e_5-e_3e_8+e_4e_7$, $e_1e_8-e_2e_7+e_3e_6-e_4e_5$
\end{itemize}
or
\begin{itemize}
\item $H_1^{{\rm sympl}} = \gf_2 \gfd_2 - \gf_1 \gfd_1$, $H_2^{{\rm sympl}} = \gf_4 \gfd_4 - \gf_3 \gfd_3$;\\[-6mm]
\item $\gfd_1 \gf_2$, $\gfd_2 \gf_1$, $\gfd_3 \gf_4$, $\gfd_4 \gf_3$;\\[-6mm]
\item $\gfd_2 \gf_4 + \gf_1 \gfd_3$, $\gfd_4 \gf_2 + \gf_3 \gfd_1$, $\gf_2 \gfd_3 - \gfd_1 \gf_4$, $\gf_3 \gfd_2 - \gfd_4 \gf_1$.
\end{itemize}
The explicit conversion table for $\gsp_4(\mC)$ reads
\begin{itemize}
\item $H_1^{{\rm sympl}} = \gf_2 \gfd_2 - \gf_1 \gfd_1 = \frac{1}{2} i \bigl ( e_1 e_2 - e_3 e_4 \bigr )$;\\[-4mm]
\item $H_2^{{\rm sympl}} = \gf_4 \gfd_4 - \gf_3 \gfd_3 = \frac{1}{2} i \bigl ( e_5 e_6 - e_7 e_8 \bigr )$;\\[-4mm]
\item $\gfd_1 \gf_2 = - \frac{1}{4} \bigl ( e_1 e_3 + e_2 e_4 - i e_1 e_4 + i e_2 e_3 \bigr ) $;\\[-5mm]
\item $\gfd_2 \gf_1 = - \frac{1}{4} \bigl ( -e_1 e_3 - e_2 e_4 - i e_1 e_4 + i e_2 e_3 \bigr ) $;\\[-5mm]
\item $\gfd_3 \gf_4 = - \frac{1}{4} \bigl ( e_5 e_7 + e_6 e_8- i e_5 e_8 + i e_6 e_7 \bigr ) $;\\[-5mm]
\item $\gfd_4 \gf_3 = - \frac{1}{4} \bigl ( -e_5 e_7 - e_6 e_8 - i e_5 e_8 + i e_6 e_7 \bigr ) $;\\[-5mm]
\item $\gfd_2 \gf_4 + \gf_1 \gfd_3 = - \frac{1}{4} \bigl ( e_3 e_7 + e_4 e_8 + e_1 e_5 + e_2 e_6 - i e_3 e_8 + i e_4 e_7 + i e_1 e_6 - i e_2 e_5 \bigr )$;\\[-5mm]
\item $\gfd_4 \gf_2 + \gf_3 \gfd_1 = - \frac{1}{4} \bigl ( -e_3 e_7 - e_4 e_8 - e_1 e_5 - e_2 e_6 - i e_3 e_8 + i e_4 e_7 + i e_1 e_6 - i e_2 e_5 \bigr )$;\\[-5mm]
\item $\gf_2 \gfd_3 - \gfd_1 \gf_4 = - \frac{1}{4} \bigl ( e_3 e_5 + e_4 e_6 - e_1 e_7 - e_2 e_8 - i e_4 e_5 + i e_3 e_6 + i e_1 e_8 - i e_2 e_7 \bigr )$;\\[-5mm]
\item $\gf_3 \gfd_2 - \gfd_4 \gf_1 = - \frac{1}{4} \bigl (- e_3 e_5 - e_4 e_6 + e_1 e_7 + e_2 e_8 - i e_4 e_5 + i e_3 e_6 + i e_1 e_8 - i e_2 e_7 \bigr )$.\\[-3mm]
\end{itemize}

Now it is well--known that $\gsp_{2p}(\mC)$ contains a copy of $\gsl_p(\mC)$. This subalgebra of $\gsp_{2p}(\mC)$ which is isomorphic with $\gsl_p(\mC)$ is generated by the following basis
\begin{itemize}
\item $H_j^{{\rm sl}} = H_j^{{\rm sympl}} - H_p^{{\rm sympl}}$, $j=1,\ldots,p-1$;\\[-6mm]
\item $\gfd_{2j} \gf_{2k} + \gf_{2j-1} \gfd_{2k-1}$, $ \gfd_{2k} \gf_{2j} + \gf_{2k-1} \gfd_{2j-1}$, $j<k=1,\ldots,p$.\\[-3mm]
\end{itemize}
For $p=2$ this explicitly means
\begin{itemize}
\item $H_1^{{\rm sl}} = H_1^{{\rm sympl}} - H_2^{{\rm sympl}} = \gf_2 \gfd_2 - \gf_1 \gfd_1 - \gf_4 \gfd_4 + \gf_3 \gfd_3$;\\[-6mm]
\item $\gfd_2 \gf_4 + \gf_1 \gfd_3$, $\gfd_4 \gf_2 + \gf_3 \gfd_1$
\end{itemize}
or: $e_1 e_2 - e_3 e_4 - e_5 e_6 + e_7 e_8$, $e_1 e_5 + e_2 e_6 + e_3 e_7 + e_4 e_8$, $e_1 e_6 - e_2 e_5 - e_3 e_8 + e_4 e_7$.\\[-2mm]

Schematically, we have the following overview picture showing, in addition, the {\em real} dimensions of the involved spaces:\\[-6mm]
\begin{center}
\begin{tikzpicture}[scale=0.9]
\draw[white] (0,0) edge node[rotate=90,black] {$\simeq$} (0,-2.5);
\draw[white] (4,0) edge node[rotate=90,black] {$\simeq$} (4,-2.5);
\draw[white] (8,0) edge node[rotate=90,black] {$\simeq$} (8,-2.5);
\draw[white] (0,-1.9) edge node[rotate=-45,black] {$\subset$} (2,-2.9);
\draw[white] (2,-3.06) edge node[rotate=-45,black] {$\subset$} (4,-4.66);
\draw[white] (4,-1.9) edge node[rotate=-45,black] {$\subset$} (6,-2.9);
\draw[white] (0,-4.5) edge node[rotate=45,black] {$\subset$} (2,-3.4);
\draw[white] (2,-2.9) edge node[rotate=45,black] {$\subset$} (4,-1.9);
\draw[white] (4,-4.5) edge node[rotate=45,black] {$\subset$} (6,-3.4);
\draw[white] (6,-2.9) edge node[rotate=45,black] {$\subset$} (8,-1.9);

\node at (0,0) {$\mC_p^{(2)}$};
\node at (0,-0.5) {$\# p(p-1)$};
\node at (4,0) {$\mC_{2p}^{(2)}$};
\node at (4,-0.5) {$\# 2p(2p-1)$};
\node at (8,0) {$\mC_{4p}^{(2)}$};
\node at (8,-0.5) {$\# 4p(4p-1)$};

\node at (0,-1.9) {$\gso_p(\mC)$};
\node at (4,-1.9) {$\gso_{2p}(\mC)$};
\node at (8,-1.9) {$\gso_{4p}(\mC)$};
\node at (2,-2.9) {$\gsl_p(\mC)$};
\node at (6,-2.9) {$\gsl_{2p}(\mC)$};
\node at (2,-3.4) {$\# 2(p^2-1)$};
\node at (6,-3.4) {$\# 2(4 p^2-1)$};
\node at (0,-4.5) {$\gsp_{p}(\mC)$};
\node at (4,-4.5) {$\gsp_{2p}(\mC)$};
\node at (0,-5) {$\# p(p+1)$};
\node at (4,-5) {$\# 2p(2p+1)$};
\end{tikzpicture}
\end{center}

which for $p=2$ explicitly reads\\[-4mm]
\begin{center}
\begin{tikzpicture}[scale=0.9]
\draw[white] (0,0) edge node[rotate=90,black] {$\simeq$} (0,-2.5);
\draw[white] (4,0) edge node[rotate=90,black] {$\simeq$} (4,-2.5);
\draw[white] (8,0) edge node[rotate=90,black] {$\simeq$} (8,-2.5);
\draw[white] (0,-1.9) edge node[rotate=-45,black] {$\subset$} (2,-2.9);
\draw[white] (2,-3.06) edge node[rotate=-45,black] {$\subset$} (4,-4.66);
\draw[white] (4,-1.9) edge node[rotate=-45,black] {$\subset$} (6,-2.9);
\draw[white] (0,-4.5) edge node[rotate=45,black] {$\subset$} (2,-3.4);
\draw[white] (2,-2.9) edge node[rotate=45,black] {$\subset$} (4,-1.9);
\draw[white] (4,-4.5) edge node[rotate=45,black] {$\subset$} (6,-3.4);
\draw[white] (6,-2.9) edge node[rotate=45,black] {$\subset$} (8,-1.9);

\node at (0,0) {$\mC_2^{(2)}$};
\node at (0,-0.5) {$\# 2$};
\node at (4,0) {$\mC_{4}^{(2)}$};
\node at (4,-0.5) {$\# 12$};
\node at (8,0) {$\mC_{8}^{(2)}$};
\node at (8,-0.5) {$\# 56$};

\node at (0,-1.9) {$\gso_2(\mC)$};
\node at (4,-1.9) {$\gso_{4}(\mC)$};
\node at (8,-1.9) {$\gso_{8}(\mC)$};
\node at (2,-2.9) {$\gsl_2(\mC)$};
\node at (6,-2.9) {$\gsl_{4}(\mC)$};
\node at (2,-3.4) {$\# 6$};
\node at (6,-3.4) {$\# 30$};
\node at (0,-4.5) {$\gsp_{2}(\mC)$};
\node at (4,-4.5) {$\gsp_{4}(\mC)$};
\node at (0,-5) {$\# 6$};
\node at (4,-5) {$\# 20$};
\end{tikzpicture}
\end{center}


\subsection{The symplectic cells of homogeneous spinor space}
\label{sympcells}


The fact that the Lie algebra $\gsp_{2p}(\mC)$ contains a copy of $\gsl_{p}(\mC)$ is crucial, since it allows for determining irreducible $\gsp_{2p}(\mC)$--modules in terms of irreducible representations for $\gsl_{p}(\mC)$. With this observation in mind, we start the quest for the decomposition of homogeneous spinor space $\mS^{r}$ into $\gsp_{2p}(\mC)$--irreducibles.\\[-2mm]

In a first step we introduce the left multiplication operators
\begin{eqnarray*}
P &=& \gf_2 \gf_1 + \gf_4 \gf_3 + \ldots + \gf_{2p} \gf_{2p-1} \\
Q &=& \gfd_1 \gfd_2 + \gfd_3 \gfd_4 + \ldots + \gfd_{2p-1} \gfd_{2p} \ = \ P^{\dagger}
\end{eqnarray*}
The action of the Witt basis vectors as left multiplication operators being
\begin{eqnarray*}
\gf_j &:& \mS^r \longrightarrow \mS^{r-1}, \; r=1,\ldots,n  \\
\gf_j &:& \mS^0 \longrightarrow \{0\} \\
\gfd_j &:& \mS^r \longrightarrow \mS^{r+1}, \; r=0,\ldots,n-1 \\
\gfd_j &:& \mS^n \longrightarrow \{0\}  
\end{eqnarray*}
we find that $P: \mS^r \rightarrow \mS^{r-2}$ and $Q: \mS^r \rightarrow \mS^{r+2}$. The operators $P$ and $Q$ generate an $\gsl_2(\mC)$--structure as is seen from the following relations which may be straightforwardly verified. These relations involve the so--called spin--Euler operator $\beta = \gfd_1 \gf_1 + \gfd_2 \gf_2 + \ldots + \gfd_{2p} \gf_{2p}$, measuring the homogeneity degree of a homogeneous spinor subspace.

\begin{lemma}
\label{lemmaS1}
One has
\begin{enumerate}
\item[(i)] $[P,Q] = p - \beta$;
\item[(ii)] $[P,\beta] = 2P$;
\item[(iii)] $[Q,\beta] = -2Q$
\end{enumerate}
\end{lemma}

We also have the following fundamental result.
\begin{lemma}
\label{lemmaS2}
The operators $P$ and $Q$ are $\gsp_{2p}(\mC)$--invariant.
\end{lemma}

\pf
The action of $\gsp_{2p}(\mC)$ being (left) multiplication, it suffices to prove that $P$ and $Q$ are commuting with the basis elements generating $\gsp_{2p}(\mC)$, listed in Section \ref{subsecgsp}. A straightforward computation shows that this is indeed the case.
\qed

\noindent Also the following results may be directly verified.
\begin{lemma}
\label{lemmaS3}
One has
\begin{enumerate}
\item[(i)] $\mbox{\em Ker} \, P |_{\mS^r} = \{0\}$ for $r=p+1,\ldots,2p$;
\item[(ii)] $\mbox{\em Ker} \, Q|_{\mS^r} = \{0\}$ for $r=0,\ldots,p-1$;
\item[(iii)] $\mbox{\em Ker} \, P|_{\mS^p} = \mbox{\em Ker} \, Q|_{\mS^p}$. 
\end{enumerate}
\end{lemma}

Now we define, for $r=0,\ldots,p$, the subspaces
$$
\mS_r^r = \mbox{Ker} \, P |_{\mS^r}, \qquad \mS_r^{2p-r} = \mbox{Ker} \, Q |_{\mS^{2p-r}}
$$
and for $k=0,\ldots,p-r$, the subspaces
$$
\mS_r^{r+2k} = Q^k \, \mS^r_r, \qquad \mS_r^{2p-r-2k} = P^k \, \mS^{2p-r}_r
$$
Notice that these definitions are compatible since $\mbox{Ker} \, P|_{\mS^p} = \mS_p^p = \mbox{Ker} \, Q|_{\mS^p}$ by Lemma \ref{lemmaS3}(iii), while $Q^k \, \mS_{p-2k}^{p-2k} = \mS_{p-2k}^{p}  = P^k \, \mS^{p+2k}_{p-2k} $ due to symmetry reasons.

\begin{lemma}
\label{lemmaS4}
One has
$$
\mbox{\em Ker} \, P|_{\mS} = \bigoplus_{r=0}^p \mS_r^r \qquad \mbox{and} \qquad \mbox{\em Ker} \, Q|_{\mS} = \bigoplus_{r=0}^p \mS_r^{2p-r} 
$$ 
\end{lemma}

\pf
This result follows from Lemma \ref{lemmaS3} and from the decomposition (\ref{decompspin}) of spinor space $\mS$ into its $\mbox{U}(2p)$--irreducible parts $\mS^r$, $r=0,\ldots,2p$.
\qed

\begin{lemma}
\label{lemmaS5}
One has, for $k= 0,1,\ldots,p-r-1$:
\begin{enumerate}
\item[(i)] $Q$ is an isomorphism $\mS_r^{r+2k} \longrightarrow \mS_r^{r+2k+2}$ with inverse $Q^{-1} = \displaystyle\frac{1}{\alpha_r^k} \, P$;\\[-4mm]
\item[(ii)]  $P$ is an isomorphism $\mS_r^{2p-r-2k} \longrightarrow \mS_r^{2p-r-2k-2}$ with inverse $P^{-1} = \displaystyle\frac{1}{\alpha_r^{p-r-k-1}} \, Q$
\end{enumerate}
where the coefficients are given by
$$
\alpha_r^k = (k+1)(p-r-k) = \alpha_r^{p-r-k-1}
$$
\end{lemma}

\pf
The proofs being similar, we only prove (i) explicitly. The operator $Q: \mS_r^{r+2k} \longrightarrow \,S_r^{r+2k+2}$ is surjective by definition. Let $\varphi_r^{r+2k} \in \mS_r^{r+2k}$ be such that $Q \varphi_r^{r+2k}=0$, then 
$$
\varphi_r^{r+2k} \in \mbox{Ker} \, Q|_{\mS} = \bigoplus_{j=0}^p \mS_j^{2p-j}
$$
As $r+2k\leq p-2$, this is impossible unless $\varphi_r^{r+2k}$ is trivially zero. So $Q|_{\mS_r^{r+2k}}$ is injective. Now assume that
$$
Q^{-1}|_{\mS_r^{r+2k-2}} = \frac{1}{\alpha_r^{k-1}} \, P|_{\mS_r^{r+2k-2}}
$$
then for $\varphi_r^{r+2k} \in \mS_r^{r+2k}$ it holds, also in view of Lemma \ref{lemmaS1}, that
$$
\frac{1}{\alpha_r^k} \, PQ \, \varphi_r^{r+2k} = \frac{1}{\alpha_r^k} \bigl ( QP +p - \beta \bigr ) \varphi_r^{r+2k}= \varphi_r^{r+2k}
$$
or still
$$
\frac{1}{\alpha_r^k} \bigl ( \alpha_r^{k-1} + p -\beta \bigr ) \varphi_r^{r+2k} = \varphi_r^{r+2k}
$$
from which it follows that 
$$
\alpha_r^k = \alpha_r^{k-1} + p - r - 2k
$$
With a similar reasoning we find, in particular, that $\alpha_r^0=p-r$, $r=0,\ldots,p-1$, allowing for a recursive computation resulting into  $\alpha_r^k = \alpha_r^{p-r-k-1} = (k+1)(p-r-k)$, $k=0,1,\ldots,p-r-1$.
\qed

\begin{corollary}
The composition of the multiplicative operators $P$ and $Q$ is constant on each symplectic cell; more specifically one has
\begin{itemize}
\item[(i)] $P \, Q = \alpha_r^k$ on \, $\mS_r^{r+2k}$ and on \, $\mS_r^{2p-r-2k-2}$
\item[(ii)] $Q \, P = \alpha_r^k$ on \,  $\mS_r^{r+2k+2}$ and on \, $\mS_r^{2p-r-2k}$
\end{itemize}
\end{corollary}
\vspace*{3mm}

With respect to the Fischer inner product
$$
\langle \lambda , \mu \rangle_r = [\lambda^\dagger \mu]_0  \quad \lambda, \mu \in \mS^r
$$
where, see Section 2, $[ \, . \, ]_0$ denotes the scalar part of a Clifford number, each of the homogeneous spinor subspaces can be decomposed as the direct sum
$$
\mS^r = \mS^r_r \oplus (\mS^r_r)^{\perp} \quad , \quad r=0,\ldots,p
$$
and
$$
\mS^{2p-r} = \mS^{2p-r}_r \oplus (\mS^{2p-r}_r)^{\perp} \quad , \quad r=0,\ldots,p
$$
where the orthogonal complements  $(\mS^r_r)^{\perp}$  and  $(\mS^{2p-r}_r)^{\perp}$ are isomorphic with Im$_P(\mS^r)$ and Im$_Q(\mS^{2p-r} )$ respectively. We will now determine those orthogonal complements explicitly.\\

\begin{lemma}
\label{lemmaS6}
With respect to the Fischer inner product, the operators $P$ and $Q$ are adjoint operators, i.e. for $\lambda \in \mS^r$ and $\mu \in \mS^{r-2}$ there holds
$$
\langle P \lambda, \mu \rangle_{r-2} = \langle \lambda , Q \mu \rangle_{r}
$$
\end{lemma}

\pf
We have indeed
$$
\langle P \lambda, \mu \rangle_{r-2} = [(P \lambda)^\dagger \mu]_0 = [\lambda^\dagger P^\dagger \mu]_0 =  [ \lambda^\dagger Q \mu]_0 =    \langle \lambda , Q \mu \rangle_{r}
$$
\eop

\begin{proposition}
For $r=0,\ldots,p$, the space $\mS^r$ may be decomposed as
\begin{eqnarray}
\label{step1}
\mS^r = \mS^r_r \oplus Q \, \mS^{r-2}
\end{eqnarray}
\end{proposition}

\pf
In fact we prove that, with respect to the Fischer inner product,
$
\left( Q \,  \mS^{r-2} \right)^{\perp} =  \mS^r_r.
$\\
Let $\lambda \in \mS^r_r$, then $P \lambda = 0$ and so $0 = \langle P \lambda, \mu \rangle_{r-2} = \langle \lambda, Q \mu \rangle_{r} $ for all $\mu \in  \mS^{r-2}$, which means that $\lambda$ is orthogonal to $Q \,  \mS^{r-2}$ or $\lambda \in \left(Q \,  \mS^{r-2} \right)^\perp$.\\
Conversely, let $\lambda \in \left( Q \,  \mS^{r-2} \right)^{\perp}$. Then $P \lambda \in \mS^{r-2}$ and 
$
\langle P \lambda, \mu \rangle_{r-2} = \langle \lambda , Q \mu \rangle_{r}  = 0
$
for all $\mu \in \mS^{r-2}$. In particular, for $\mu = P \lambda$ we find
$
\langle P \lambda, P \lambda \rangle_{r-2} = 0
$
whence $P \lambda= 0$ or $\lambda \in \mS^r_r$.
\eop

\noindent
In a similar way the following complementary result is obtained.
\begin{proposition}
For $r=0,\ldots,p$, the space $\mS^{2p-r}$ may be decomposed as
\begin{eqnarray}
\label{step1bis}
\mS^{2p-r} = \mS^{2p-r}_r \oplus P \, \mS^{2p-r+2}
\end{eqnarray}
\end{proposition}

\noindent
Consecutive application of the decompositions (\ref{step1}) and (\ref{step1bis}) leads to the following result.

\begin{proposition}
\label{propositionS1}
One has, for all $r=0,\ldots,p$:
$$
\mS^r = \bigoplus_{j=0}^{\lfloor \frac{r}{2} \rfloor} \mS_{r-2j}^r \qquad \mbox{and} \qquad
\mS^{2p-r} = \bigoplus_{j=0}^{\lfloor \frac{r}{2} \rfloor} \mS_{r-2j}^{2p-r}
$$
\end{proposition}

\begin{corollary}
\label{corollaryS1}
One has, for $r=0,\ldots,p$ and $k=0,\ldots,\lfloor \frac{p}{2} \rfloor$:
\begin{enumerate}
\item[(i)] $\displaystyle \mbox{\em dim} \, \bigl ( \mS_r^{r+2k} \bigr ) = \binom{2p}{r} - \binom{2p}{r-2}$;
\item[(ii)] $\displaystyle \mbox{\em dim} \, \bigl ( \mS_r^{2p-r-2k} \bigr ) = \mbox{\em dim} \, \bigl( \mS_r^{r+2k} \bigr )$.
\end{enumerate}
\end{corollary}

\pf
(i) One has
$$
\mbox{dim} \, \bigl ( \mS_r^{r+2k} \bigr ) = \mbox{dim} \, \bigl ( \mS_r^r \bigr ) = \mbox{dim} \, \bigl ( \mS^r \bigr ) - \mbox{dim} \, \bigl ( \mS^{r-2} \bigr ) =  \binom{2p}{r} - \binom{2p}{r-2}
$$
and, in particular, $\mbox{dim} \, \bigl ( \mS_0^{2k} \bigr ) = \mbox{dim} \, \bigl ( \mS_0^0 \bigr ) = 1$ and
$\mbox{dim} \, \bigl ( \mS_1^{2k+1} \bigr ) = \mbox{dim} \, \bigl ( \mS_1^1 \bigr ) = 2p$.\\[2mm]
(ii) One has
\begin{eqnarray*}
\mbox{dim} \, \bigl ( \mS_r^{2p-r-2k} \bigr ) &=& \mbox{dim} \, \bigl ( \mS_r^{2p-r} \bigr ) = \mbox{dim} \, \bigl ( \mS^{2p-r} \bigr ) - \mbox{dim} \, \bigl ( \mS^{2p-r+2} \bigr ) \\
& = & \binom{2p}{2p-r} - \binom{2p}{2p-r+2} = \binom{2p}{r} - \binom{2p}{r-2} = \mbox{dim} \, \bigl ( \mS_r^{r+2k} \bigr )
\end{eqnarray*}
\qed

Now we claim that each of the cells $\mS_s^r$ in the above decompositions, is an irreducible $\mbox{Sp}(p)$--representation. We know that the homogeneous spinor spaces $\mS^{r}, r=0,\ldots,2p$ are irreducible modules for $\gsl_{2p}(\mC)$ and also for $\gl_{2p}(\mC)$; hence they are also invariant for $\gsp_{2p}(\mC) \subset \gsl_{2p}(\mC) \subset \gl_{2p}(\mC)$. Moreover the operators $P$ and $Q$ are invariant under $\gsp_{2p}(\mC)$; hence $\mS_r^{r} = \mS^{r} \cap \, \mbox{Ker} \, P$ and $\mS_r^{2p-r} = \mS^{2p-r} \cap \, \mbox{Ker} \, Q$ also are $\gsp_{2p}(\mC)$--invariant. They are however, a priori, not necessarily irreducible. This irreducibility will now be proven through branching when restricting $\gl_{2p}(\mC)$ to $\gsp_{2p}(\mC)$.
The corresponding branching rules could be found in full generality in \cite{HTW}, the branching multiplicities being expressed in terms of Littlewood--Richardson coefficients. However, due to the simple highest weight to start with, the actual branching is rather straightforward, and one obtains
\begin{equation}
\label{branching}
 \mS^r \Biggr\rvert_{\gsp_{2p}(\mC)}^{\gl_{2p}(\mC)} = (1_r)_s \oplus  (1_{r-2})_s \oplus \cdots \oplus  (1_{r-2\lfloor \frac{r}{2} \rfloor})_s
\end{equation}
where the shorthand notation $(1_r)_s$ refers to an irreducible representation for $\gsp_{2p}(\mC)$, and stands for the {\em symplectic} highest weight $(\underbrace{1,\ldots,1}_r, \underbrace{0,\ldots,0}_{p-r})$.

\begin{theorem}
\label{symplecticcells}
For $r=0,\ldots,p$ one has $\mS^r_r \cong (1_r)_s$ and
$$
\mS^r = \mS^r_r \oplus \mS^r_{r-2} \oplus \cdots \oplus \mS^r_{r-2\lfloor \frac{r}{2}\rfloor}
$$
is an $\gsp_{2p}(\mC)$--irreducible decomposition.
\end{theorem}

\pf
We proceed by induction. For $r=0,1$ the statement is trivial.\\
Assume that $\mS^{r-2}_{r-2} \cong (1_{r-2})_s$ and that
$$
\mS^{r-2} = \mS^{r-2}_{r-2} \oplus \mS^{r-2}_{r-4} \oplus \mS^{r-2}_{r-6} \oplus \cdots = 
\mS^{r-2}_{r-2} \oplus Q \, \mS^{r-4}_{r-4} \oplus Q^2 \, \mS^{r-6}_{r-6} \oplus \cdots 
$$
is an $\gsp_{2p}(\mC)$--irreducible decomposition. Then also
$$
Q \, \mS^{r-2} = Q \, \mS^{r-2}_{r-2} \oplus Q^2 \, \mS^{r-4}_{r-4} \oplus Q^3 \, \mS^{r-6}_{r-6} \oplus \cdots
$$
is an $\gsp_{2p}(\mC)$--irreducible decomposition, which also reads
$$
Q \, \mS^{r-2} = (1_{r-2})_s \oplus (1_{r-4})_s \oplus (1_{r-6})_s \oplus \cdots
$$
In view of the decomposition (\ref{step1}) and the branching (\ref{branching}) it follows that
$\mS^r_r \cong (1_r)_s$, which finishes the proof.
\eop\\

\noindent
In a similar way we obtain the following complementary result.

\begin{theorem}
\label{symplecticcellsbis}
For $r=0,\ldots,p$ one has $\mS^{2p-r}_r = (1_r)_s$ and
$$
\mS^{2p-r} = \mS^{2p-r}_r \oplus \mS^{2p-r}_{r-2} \oplus \cdots \oplus \mS^{2p-r}_{r-2\lfloor \frac{r}{2}\rfloor}
$$
is an $\gsp_{2p}(\mC)$--irreducible decomposition.
\end{theorem}

\begin{remark}
Recall that the real Lie algebra $\gsp(p)$ of skew--symplectic $M_p(\mH)$--matrices is isomorphic with the compact form $\gsp_{2p}(\mC) \cap \gu(2p)$ of the complex symplectic Lie algebra $\gsp_{2p}(\mC)$. So Theorems \ref{symplecticcells} and \ref{symplecticcellsbis} implie that all spaces $\mS_{r-2j}^{r}$, $r=0,\ldots,2p$, $j = 0,\ldots, \lfloor \frac{r}{2} \rfloor$, are irreducible $\gsp(p)$--modules as well.
\end{remark}

\noindent
We illustrate the above Fischer decomposition of the homogeneous spinor subspaces into symplectic cells by the following general triangular scheme and by some examples in low dimension.
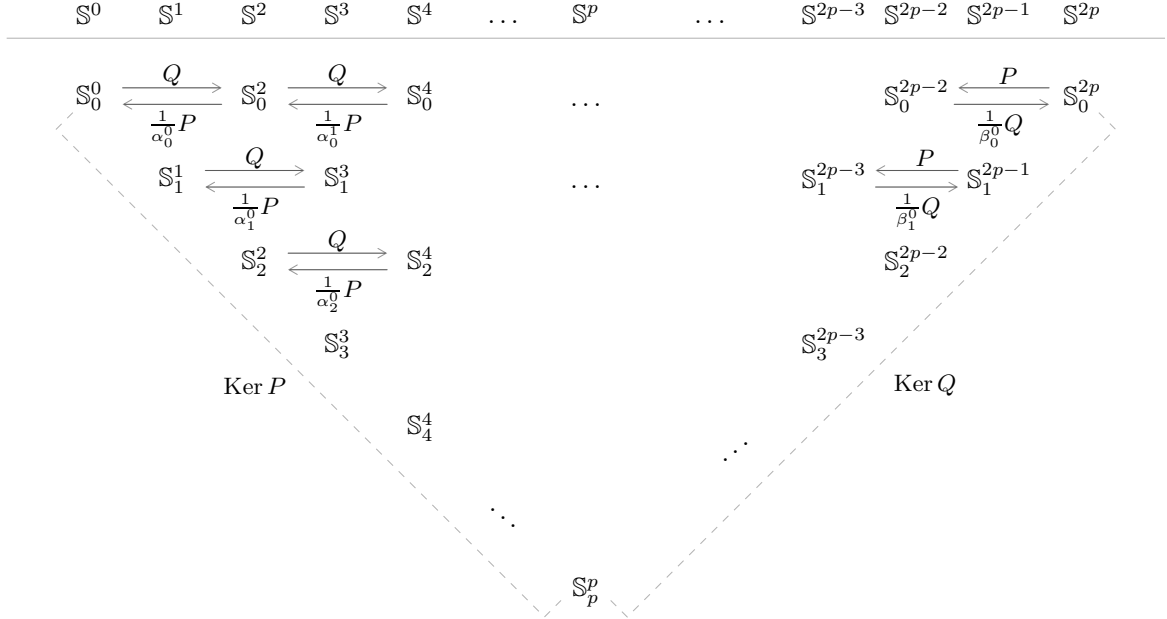
\begin{center}
\begin{tikzpicture}[scale=1.1]
\node at (0,0) {$\mS^0$};
\node at (1,0) {$\mS^1$};
\node at (2,0) {$\mS^2$};
\node at (3,0) {$\mS^3$};
\node at (4,0) {$\mS^4$};
\node at (5,-0.1) {$\ldots$};
\node at (6,0) {$\mS^p$};
\node at (7.5,-0.1) {$\ldots$};
\node at (9,0) {$\mS^{2p-3}$};
\node at (10,0) {$\mS^{2p-2}$};
\node at (11,0)  {$\mS^{2p-1}$};
\node at (12,0) {$\mS^{2p}$};
\draw[-, draw opacity=0.3] (-1,-0.3) -- (13,-0.3);
\node at (0,-1) {$\mS^0_0$};
\node at (2,-1) {$\mS^2_0$};
\node at (4,-1) {$\mS^4_0$};
\node at (6,-1.1) {$\ldots$};
\node at (6,-2.1) {$\ldots$};
\node at (10,-1) {$\mS^{2p-2}_0$};
\node at (12,-1) {$\mS^{2p}_0$};
\node at (1,-2) {$\mS^1_1$};
\node at (3,-2) {$\mS^3_1$};
\node at (11,-2) {$\mS^{2p-1}_1$};
\node at (9,-2) {$\mS^{2p-3}_1$};
\node at (2,-3) {$\mS_2^2$};
\node at (4,-3) {$\mS_2^4$};
\node at (10,-3) {$\mS_2^{2p-2}$};
\node at (9,-4) {$\mS_3^{2p-3}$};
\node at (3,-4) {$\mS_3^3$};
\node at (4,-5) {$\mS_4^4$};
\node at (5,-6) {$\ddots$};
\node at (6,-7) {$\mS_p^p$};
\node at (7.8,-5.2) {$\iddots$};
\draw[->,>=angle 60,draw opacity=0.4] (0.4,-0.9) -- (1.6,-0.9);
\node at (1,-0.75) {\small $Q$};
\draw[<-,>=angle 60,draw opacity=0.4] (0.4,-1.1) -- (1.6,-1.1);
\draw[<-,>=angle 60,draw opacity=0.4] (10.5,-0.9) -- (11.6,-0.9);
\draw[->,>=angle 60,draw opacity=0.4] (10.45,-1.1) -- (11.6,-1.1);
\node at (11.1,-0.75) {\small $P$};
\node at (3,-0.75) {\small $Q$};
\draw[->,>=angle 60,draw opacity=0.4] (2.4,-0.9) -- (3.6,-0.9);
\draw[<-,>=angle 60,draw opacity=0.4] (2.4,-1.1) -- (3.6,-1.1);
\node at (2,-1.75) {\small $Q$};
\draw[->,>=angle 60,draw opacity=0.4] (1.4,-1.9) -- (2.6,-1.9);
\draw[<-,>=angle 60,draw opacity=0.4] (1.4,-2.1) -- (2.6,-2.1);
\draw[<-,>=angle 60,draw opacity=0.4] (9.5,-1.9) -- (10.5,-1.9);
\draw[->,>=angle 60,draw opacity=0.4] (9.5,-2.1) -- (10.5,-2.1);
\node at (10.1,-1.75) {\small $P$};
\node at (3,-2.75) {\small $Q$};
\draw[->,>=angle 60,draw opacity=0.4] (2.4,-2.9) -- (3.6,-2.9);
\draw[<-,>=angle 60,draw opacity=0.4] (2.4,-3.1) -- (3.6,-3.1);
\draw[draw opacity = 0.3,dashed] (-0.2,-1.2) -- (-0.4,-1.4) -- (5.5,-7.3) -- (5.7,-7.1) ;
\node at (2,-4.5) {\small $\mbox{Ker} \, P$};
\draw[draw opacity = 0.3,dashed] (6.3,-7.1) -- (6.5,-7.3) -- (12.4,-1.4) -- (12.2,-1.2);
\node at (10.1,-4.5) {\small $\mbox{Ker} \, Q$};
\node at (1,-1.4) {\small $\frac{1}{\alpha_0^0} P$};
\node at (3,-1.4) {\small $\frac{1}{\alpha_0^1} P$};
\node at (2,-2.4) {\small $\frac{1}{\alpha_1^0} P$};
\node at (3,-3.4) {\small $\frac{1}{\alpha_2^0} P$};
\node at (11,-1.4) {\small $\frac{1}{\beta_0^0} Q$};
\node at (10,-2.4) {\small $\frac{1}{\beta_1^0} Q$};
\end{tikzpicture}
\end{center}

\noindent {\bf Example} ($p=1$)\\
For $p=1$ we have that $P=\gf_2 \gf_1$ and $Q=\gfd_1\gfd_2$. The complex dimension of spinor space is $\mbox{dim}(\mS)=2^2=4$ and it decomposes into $\mS^0 = \mS_0^0$, $\mS^1 = \mS_1^1$ and $\mS^2=\mS_0^2$. The triangular scheme, including basis elements, reduces to:
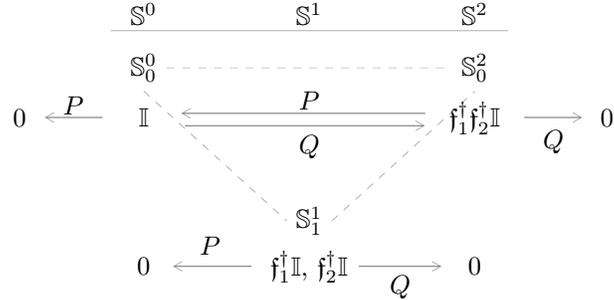
\begin{center}
\begin{tikzpicture}[scale=1.1]
\node at (0,0) {$\mS^0$};
\node at (2,0) {$\mS^1$}; 
\node at (4,0) {$\mS^2$};
\draw[draw opacity = 0.3] (-0.4,-0.2) -- (4.4,-0.2);
\node at (0,-0.65) (S00) {$\mS_0^0$};
\node at (0,-1.25) {$\mI$};
\node at (4,-0.65) (S02) {$\mS_0^2$};
\node at (4,-1.25) {$\gfd_1 \gfd_2 \mI$};
\draw[<-,>=angle 60,draw opacity = 0.3] (0.45,-1.2) -- (3.4,-1.2);
\draw[->,>=angle 60,draw opacity = 0.3] (0.5,-1.35) -- (3.4,-1.35);
\node at (2,-1.05) {$P$};
\node at (2,-1.6) {$Q$};
\draw[->,>=angle 60,draw opacity = 0.3] (4.6,-1.25) -- (5.3,-1.25);
\node at (5.6,-1.25) {$0$};
\node at (4.95,-1.55) {$Q$};
\draw[<-,>=angle 60, draw opacity = 0.3] (-1.2,-1.25) -- (-0.5,-1.25);
\node at (-1.5,-1.25) {$0$};
\node at (-0.85,-1.1) {$P$};
\node at (2,-2.5) (S11) {$\mS^1_1$}; 
\node at (0,-3.05) {$0$};
\node at (4,-3.05) {$0$};
\node at (2,-3.05) {$\gfd_1 \mI$, $\gfd_2 \mI$};
\draw[<-,>=angle 60, draw opacity = 0.3] (0.35,-3.05) -- (1.3,-3.05);
\node at (0.8,-2.8) {$P$};
\draw[->,>=angle 60, draw opacity = 0.3] (2.6,-3.05) -- (3.6,-3.05);
\node at (3.1,-3.3) {$Q$};
\draw[dashed,draw opacity = 0.2] (S00.east) to (S02.west);
\draw[dashed,draw opacity = 0.3] (S00.south) to (S11.west);
\draw[dashed,draw opacity = 0.3] (S11.east) to (S02.south);
\end{tikzpicture}
\end{center}

\noindent {\bf Example} ($p=2$)\\
Here we have that $P=\gf_2 \gf_1 + \gf_4 \gf_3$ and $Q = \gfd_1 \gfd_2 + \gfd_3 \gfd_4$. The complex dimension of spinor space is $\mbox{dim}(\mS)=2^3=8$ and it decomposes into $\mS^0 = \mS_0^0$, $\mS^1 = \mS_1^1$, $\mS^2=\mS_0^2 \oplus \mS_2^2$, $\mS^3 = \mS_1^3$ and $\mS^4 = \mS_0^4$. The respective complex dimensions of the homogeneous parts and the symplectic cells are given by
\begin{eqnarray*}
\mbox{dim}(\mS^0) & = & \mbox{dim}(\mS_0^0) \ =\ 1 \ = \  \mbox{dim}(\mS_0^4) \ = \ \mbox{dim}(\mS^4)\\
\mbox{dim}(\mS^1) & = & \mbox{dim}(\mS_1^1) \ = \ 4 \ = \ \mbox{dim}(\mS_1^3) \ = \ \mbox{dim}(\mS^3)
\end{eqnarray*}
and
$$
\mbox{dim}(\mS^2) = 6, \quad \mbox{with } \mbox{dim}(\mS_0^2) = 1 \mbox{ and } \mbox{dim}(\mS_2^2) = 5
$$
The triangular scheme now reads:
\begin{center}
\begin{tikzpicture}[scale=1.1]
\node at (0,0) {$\mS^0$};
\node at (2,0) {$\mS^1$}; 
\node at (4,0) {$\mS^2$};
\node at (6,0) {$\mS^3$}; 
\node at (8,0) {$\mS^4$};
\draw[draw opacity = 0.3] (-0.4,-0.2) -- (8.4,-0.2);
\node at (0,-0.65) (S00) {$\mS_0^0$};
\node at (0,-1.25) (I) {$\mI$};
\node at (4,-0.65) (S02) {$\mS_0^2$};
\node at (4,-1.25) (QI) {$(\gfd_1 \gfd_2 + \gfd_3 \gfd_4) \mI$};
\node at (8,-0.65) (S04) {$\mS_0^4$};
\node at (8,-1.25) (QQI) {$\gfd_1 \gfd_2 \gfd_3 \gfd_4 \mI$};
\draw[dashed,draw opacity = 0.2] (S00.east) to (S02.west);
\draw[dashed,draw opacity = 0.2] (S02.east) to (S04.west);
\draw[->,>=angle 60,draw opacity=0.3] (0.4,-1.2) to (3,-1.2);
\draw[<-,>=angle 60,draw opacity=0.3] (0.4,-1.35) to (3,-1.35);
\node at (1.7,-1.1) {$Q$};
\node at (1.7,-1.5) {$P$};
\draw[->,>=angle 60,draw opacity=0.3] (8.7,-1.25) to (9.3,-1.25);
\node at (9,-1.15) {$Q$};
\node at (9.5,-1.25) {$0$};
\draw[<-,>=angle 60,draw opacity=0.3] (-0.8,-1.25) to (-0.2,-1.25);
\node at (-0.5,-1.4) {$P$};
\node at (-1,-1.25) {$0$};
\draw[->,>=angle 60,draw opacity=0.3] (5,-1.2) to (7.3,-1.2);
\draw[<-,>=angle 60,draw opacity=0.3] (5,-1.35) to (7.3,-1.35);
\node at (6.2,-1.1) {$Q$};
\node at (6.2,-1.5) {$P$};
\node at (2,-2.65) (S11) {$\mS^1_1$}; 
\node at (6,-2.65) (S13) {$\mS^3_1$}; 
\draw[dashed,draw opacity = 0.3] (S00.south) to (S11.west);
\draw[dashed,draw opacity = 0.3] (S04.south) to (S13.east);
\node at (2,-3.25) {$\gfd_1 \mI$, $\gfd_2 \mI$, $\gfd_3 \mI$, $\gfd_4 \mI$};
\node at (6,-3.05) {$\gfd_1 \gfd_3 \gfd_4 \mI$, $\gfd_2 \gfd_3 \gfd_4 \mI$};
\node at (6,-3.5) {$\gfd_1 \gfd_2 \gfd_3 \mI$, $\gfd_1 \gfd_2 \gfd_4 \mI$};
\draw[->,>=angle 60,draw opacity=0.3] (3.4,-3.2) to (4.6,-3.2);
\draw[<-,>=angle 60,draw opacity=0.3] (3.4,-3.35) to (4.6,-3.35);
\node at (4,-3.1) {$Q$};
\node at (4,-3.5) {$P$};
\node at (4,-4.65) (S22) {$\mS^2_2$}; 
\draw[dashed,draw opacity = 0.3] (S11.south) to (S22.west);
\draw[dashed,draw opacity = 0.3] (S13.south) to (S22.east);
\node at (4,-5.05)  {$\gfd_1 \gfd_3 \mI$, $\gfd_1 \gfd_4 \mI$, $\gfd_2 \gfd_3 \mI$, $\gfd_2 \gfd_4 \mI$};
\node at (4,-5.5) {$(\gfd_1 \gfd_2 - \gfd_3 \gfd_4) \mI$};
\draw[<-,>=angle 60,draw opacity=0.3] (1.6,-5.35) -- (2.2,-5.35);
\node at (1.9,-5.5) {$P$};
\draw[->,>=angle 60,draw opacity=0.3] (5.8,-5.35) -- (6.4,-5.35);
\node at (6.1,-5.25) {$Q$};
\node at (1.3,-5.35) {$0$};
\node at (6.7,-5.35) {$0$};
\node at (-0.5,-3.25) {$0$};
\node at (8.5,-3.25) {$0$};
\draw[<-,>=angle 60, draw opacity = 0.3] (-0.15,-3.25) -- (0.8,-3.25);
\node at (0.4,-3.4) {$P$};
\draw[->,>=angle 60, draw opacity = 0.3] (7.1,-3.25) -- (8.1,-3.25);
\node at (7.6,-3.1) {$Q$};
\end{tikzpicture}
\end{center}

\noindent {\bf Example} ($p=3$)\\
Here we have that $P=\gf_2 \gf_1 + \gf_4 \gf_3 + \gf_6 \gf_5$ and $Q = \gfd_1 \gfd_2 + \gfd_3 \gfd_4 + \gfd_5 \gfd_6$. The complex dimension of spinor space is $\mbox{dim}(\mS)=2^6=64$ and it decomposes into $\mS^0 = \mS_0^0$, $\mS^1 = \mS_1^1$, $\mS^2=\mS_0^2 \oplus \mS_2^2$, $\mS^3 = \mS_1^3 \oplus \mS_3^3$, $\mS^4 = \mS_0^4 \oplus \mS_2^4$, $\mS^5 = \mS_1^5$ and $\mS^6 =\mS_0^6$. The respective complex dimensions of the symplectic cells are given by
\begin{eqnarray*}
\mbox{dim}(\mS_0^0) \ = \ \mbox{dim}(\mS_0^2) \ = \ \mbox{dim}(\mS_0^4) \ = \ \mbox{dim}(\mS_0^6) \ &=& \ \phantom{1}1 \\
\mbox{dim}(\mS_1^1) \ = \ \mbox{dim}(\mS_1^3)  \ = \ \mbox{dim}(\mS_1^5) \ &=& \ \phantom{1}6\\
\mbox{dim}(\mS_2^2) \ = \ \mbox{dim}(\mS_2^4) \ &=& \ 14\\
\mbox{dim}(\mS_3^3) \ &=& \ 14
\end{eqnarray*}
Respective bases are given by
\begin{eqnarray*}
\mS_0^0 &:& \mI\\
\mS_0^2 &:& (\gfd_1\gfd_2 + \gfd_3\gfd_4 + \gfd_5\gfd_6) \, \mI\\
\mS_0^4 &:& (\gfd_1\gfd_2  \gfd_3\gfd_4 + \gfd_1\gfd_2\gfd_5\gfd_6 + \gfd_3\gfd_4 \gfd_5\gfd_6) \, \mI \\
\mS_0^6 &:& \gfd_1\gfd_2 \gfd_3\gfd_4 \gfd_5\gfd_6 \, \mI\\ \\
\mS_1^1 &:& \gfd_1 \, \mI, \gfd_2 \, \mI, \gfd_3 \, \mI, \gfd_4 \, \mI, \gfd_5 \, \mI, \gfd_6 \, \mI \\
\mS_1^3 &:& (\gfd_1\gfd_3\gfd_4 + \gfd_1\gfd_5\gfd_6) \, \mI, 
(\gfd_2\gfd_3\gfd_4 + \gfd_2\gfd_5\gfd_6) \, \mI,
(\gfd_1\gfd_2\gfd_3 + \gfd_3\gfd_5\gfd_6) \, \mI,\\
&& (\gfd_1\gfd_2\gfd_4 + \gfd_4\gfd_5\gfd_6) \, \mI,
(\gfd_1\gfd_2\gfd_5 + \gfd_3\gfd_4\gfd_5) \, \mI,
(\gfd_1\gfd_2\gfd_6 + \gfd_3\gfd_4\gfd_6) \, \mI\\
\mS_1^5 &:& \gfd_1\gfd_3\gfd_4 \gfd_5\gfd_6 \, \mI, 
\gfd_2\gfd_3\gfd_4 \gfd_5\gfd_6 \, \mI,
\gfd_1\gfd_2\gfd_3 \gfd_5\gfd_6 \, \mI,
\gfd_1\gfd_2\gfd_4 \gfd_5\gfd_6 \, \mI,
\gfd_1\gfd_2\gfd_3 \gfd_4\gfd_5 \, \mI,
\gfd_1\gfd_2\gfd_3 \gfd_4\gfd_6 \, \mI \\ \\
\mS_2^2 &:& \gfd_1\gfd_3 \, \mI, \gfd_1\gfd_4 \, \mI, \gfd_1\gfd_5 \, \mI, \gfd_1\gfd_6 \, \mI, \gfd_2\gfd_3 \, \mI, \gfd_2\gfd_4 \, \mI, \gfd_2\gfd_5 \, \mI, \gfd_2\gfd_6 \, \mI, \gfd_3\gfd_5 \, \mI,  \gfd_3\gfd_6 \, \mI,  \gfd_4\gfd_5 \, \mI,  \gfd_4\gfd_6 \, \mI, \\
&& (2 \gfd_1\gfd_2 - \gfd_3\gfd_4 - \gfd_5\gfd_6) \, \mI , 
(- \gfd_1\gfd_2 + 2 \gfd_3\gfd_4 - \gfd_5\gfd_6) \, \mI \\
\mS_2^4 &:& \gfd_1\gfd_3\gfd_5\gfd_6 \, \mI, \gfd_1\gfd_4\gfd_5\gfd_6 \, \mI, \gfd_1\gfd_3\gfd_4\gfd_5 \, \mI, \gfd_1\gfd_3\gfd_4\gfd_6 \, \mI, \gfd_2\gfd_3\gfd_5\gfd_6 \, \mI, \gfd_2\gfd_4\gfd_5\gfd_6 \, \mI, \\
&&  \gfd_2\gfd_3\gfd_4\gfd_5 \, \mI, \gfd_2\gfd_3\gfd_4\gfd_6 \, \mI, \gfd_1\gfd_2\gfd_3\gfd_5 \, \mI, \gfd_1\gfd_2\gfd_3\gfd_6 \, \mI, \gfd_1\gfd_2\gfd_4\gfd_5 \, \mI, \gfd_1\gfd_2\gfd_4\gfd_6 \, \mI,\\
&& (\gfd_1\gfd_2\gfd_3\gfd_4 + \gfd_1\gfd_2\gfd_5\gfd_6 - 2 \gfd_3\gfd_4\gfd_5\gfd_6� \, \mI, (\gfd_1\gfd_2\gfd_3\gfd_4 -2 \gfd_1\gfd_2\gfd_5\gfd_6 + \gfd_3\gfd_4\gfd_5\gfd_6� \, \mI
\\ \\
\mS_3^3 &:& \gfd_1\gfd_3\gfd_5 \, \mI, \gfd_1\gfd_3\gfd_6 \, \mI, \gfd_1\gfd_4\gfd_5 \, \mI, \gfd_1\gfd_4\gfd_6 \, \mI, \gfd_2\gfd_3\gfd_5 \, \mI, \gfd_2\gfd_3\gfd_6 \, \mI, \gfd_2\gfd_4\gfd_5 \, \mI, \gfd_2\gfd_4\gfd_6 \, \mI\\
&& ( \gfd_1\gfd_3\gfd_4 -  \gfd_1\gfd_5\gfd_6) \, \mI, ( \gfd_2\gfd_3\gfd_4 -  \gfd_2\gfd_5\gfd_6) \, \mI, ( \gfd_1\gfd_2\gfd_3 -  \gfd_3\gfd_5\gfd_6) \, \mI, \\
&&( \gfd_1\gfd_2\gfd_4 -  \gfd_4\gfd_5\gfd_6) \, \mI, ( \gfd_1\gfd_2\gfd_5 -  \gfd_3\gfd_4\gfd_5) \, \mI, ( \gfd_1\gfd_2\gfd_6 -  \gfd_3\gfd_4\gfd_6) \, \mI
\end{eqnarray*}


\subsection{Projection on the symplectic cells}


We now aim at determining the projection operators $\Pi_s^r: \mS^r \longrightarrow \mS_s^r$,  $r=0,\ldots,p$, $s=r,r-2,\ldots$, the projection operators $\Pi_s^{2p-r}: \mS^{2p-r} \longrightarrow \mS_s^{2p-r}$,  $r=0,\ldots,p$, $s=r,r-2,\ldots$ being completely similar by symmetry.\\[-2mm]

We first proceed by direct calculation. Taking $r$ to be even, say $r=2j$, and turning our attention to the first row in the triangular scheme above, i.e.\ taking $s=0$, we thus aim at determining the projection operators $\Pi_0^{2j} : \mS^{2j} \rightarrow \mS_0^{2j}$, $j=0,\ldots,p$. Since $\mS^0 \equiv \mS^0_0$, it is clear that 
$$
\Pi_0^0 = \mathbf{1}
$$
Next we consider $\mS^2$, which decomposes as $\mS^2_0 \oplus \mS^2_2$ according to Proposition \ref{propositionS1}. We thus have that $P \mS^2 = P \mS_0^2 = \mS_0^0$ on account of Lemma \ref{lemmaS4} and Lemma \ref{lemmaS5}(ii), whence $QP \mS^2 = Q \mS_0^0 = \mS_0^2$ on account of Lemma \ref{lemmaS5}(i). In view of this observation, we claim that the projection operator on $\mS_0^2$ will be given by $\Pi_0^2 = \alpha \, QP$ where the constant $\alpha$ can be determined by expressing the fact that the element $Q \mI \in \mS_0^2$ should be mapped onto itself. Invoking Lemma \ref{lemmaS1}(i) we have that 
$$
\alpha \, QP \bigl ( Q \mI \bigr ) = \alpha Q (QP + p - \beta) \mI = \alpha  p  \bigl (Q \mI \bigr )
$$
since $\mI \in \mS_0^0 \subset \mbox{Ker} \, P$, which finally yields
$$
\Pi_0^2 = \frac{1}{p} \, QP
$$
Proceeding in the same way, we postulate the general form $\Pi_0^{2j} = \alpha \, Q^j P^j$. Indeed, according to Proposition \ref{propositionS1}, the space $\mS^{2j}$ decomposes as
$$
\mS^{2j} = \mS_0^{2j} \oplus \mS^{2j}_2 \oplus \ldots \oplus \mS^{2j}_{2j-2} \oplus \mS_{2j}^{2j}
$$    
By combination of Lemma \ref{lemmaS4} and Lemma \ref{lemmaS5}(ii), we may conclude that the action of $P^j$ makes us shift $j$ positions to the left in the same row of the triangular scheme, resulting in $P^j (\mS^{2j}) = P^j (\mS_0^{2j}) = \mS_0^0$. The subsequent action of $Q^j$ compensates for this shift, see Lemma \ref{lemmaS5}(i); more precisely $Q^j P^j (\mS^{2j}) = Q^j (\mS_0^0)= \mS_0^{2j}$, confirming the proposed form of $\Pi_0^{2j}$. Finally, the constant $\alpha$ is determined by letting act this proposed projection operator on the element $Q^j \mI \in \mS_0^{2j}$, which should be mapped onto itself. Direct calculation, involving the repeated application of Lemma \ref{lemmaS1}(i),  yields the auxiliary result
$$
Q^j P^j \bigl ( Q^j \mI \bigr ) = j! \, p(p-1)(p-2) \ldots (p-j+1) \, Q^j \mI
$$
whence eventually
\begin{equation}
\Pi_0^{2j} = \frac{1}{j! \, p(p-1)(p-2) \ldots (p-j+1)} \, Q^j P^j, \qquad j=1,\ldots,p
\label{Pi02j}
\end{equation}
We may now also determine all projection operators $\Pi_2^{2j}$. Indeed, the above decomposition of $\mS^{2j}$ allows us to write
$$
\mathbf{1}  = \Pi_0^{2j} + \Pi_2^{2j} + \Pi_4^{2j} + \ldots + \Pi_{2j}^{2j}, \quad \mbox{or} \quad \mathbf{1}  - \Pi_0^{2j} = \Pi_2^{2j} + \Pi_4^{2j} + \ldots + \Pi_{2j}^{2j}
$$
So the shifting property to the left of the operator $P$ implies that $P^{j-1} \bigl ( \mathbf{1} - \Pi_0^{2j} \bigr ) (\mS^{2j}) = P^{j-1} \bigl ( \mS_2^{2j} \bigr ) = \mS_2^2$, which is again neutralized by the shifting property to the right of the operator $Q$, yielding $Q^{j-1} P^{j-1} (\mS^{2j}) = \mS_2^{2j}$. This enables us to claim that 
$$
\Pi_2^{2j} = \alpha \, Q^{j-1} P^{j-1} \bigl ( \mathbf{1} - \Pi_0^{2j} \bigr )
$$
where, by considering the action of $\Pi_2^{2j}$ on a typical basis element of $\mS_2^{2j}$, the constant $\alpha$ is found to be
$$
\alpha = \frac{1}{(j-1)! \, (p-2)(p-3) \ldots (p-j)}
$$
 \\[-2mm]

Next, taking $r$ to be odd, say $r=2j+1$, and considering the corresponding first row in the triangular scheme, i.e.\ the row where $s=1$, an analogous procedure leads to
\begin{equation}
\Pi_1^{2j+1} =  \frac{1}{j! \, (p-1)(p-2) \ldots (p-j)} \, Q^j P^j, \qquad j=0,\ldots,p-1
\label{Pi12j+1}
\end{equation}
and formulae of similar structure for the projections $\Pi_3^{2j+1}$ and higher.\\[-2mm]

Clearly, this approach eventually will lead to explicit forms for all projection operators, the involved calculations becoming more and more lengthy, though.
However, these first results are also found to fit neatly into a more abstract approach. Indeed, in Lemma \ref{lemmaS1} it was stated that the operators $P$ and $Q$ generate an $\gsl_2(\mC)$--structure, so we may consider the so--called Casimir operator $\mathcal{C}$, given by
$$
\mathcal{C} = QP + \frac{1}{4} H(H+2)
$$
with $H = [P,Q] = p-\beta$.
It may be readily checked that all elements of the symplectic spinor cell $\mS^s_s$ are eigenvectors of the Casimir operator, with eigenvalue $c_s = \frac{1}{4} (p-s)(p+2-s)$, $s=0,\ldots,p$. However, there is more. Application of Lemma \ref{lemmaS5} reveals that, in fact, all symplectic cells on the same row in the triangular scheme belong to the same eigenspace corresponding to the eigenvalue $c_s$. This allows us to write (in accordance with an abstract result, see e.g.\ \cite{fulhar}) the projection operators on the individual cells as

\begin{equation}
\Pi_{2s}^r = \prod_{\stackrel{k=0}{2k \ne 2s}}^{\lfloor \frac{r}{2} \rfloor} \frac{\mathcal{C} - c_{2k}}{c_{2s}-c_{2k}} \quad {\rm and} \quad 
\Pi_{2s+1}^r = \prod_{\stackrel{k=0}{2k+1 \ne 2s+1}}^{\lfloor \frac{r}{2} \rfloor}   \frac{\mathcal{C} - c_{2k+1}}{c_{2s+1}-c_{2k+1}}
\label{casimir}
\end{equation}
Indeed, letting act the above operator on $\mS^r$, which is decomposed as in Proposition \ref{propositionS1}, it is now directly seen that it annihilates all symplectic components except $\mS_s^r$, on which it acts as the identity operator.\\
Let us illustrate this general formula by means of some examples. Let $r=2$, then
$$
\Pi_0^2 = \frac{\mathcal{C} - c_2}{c_0-c_2} \quad \mbox{and} \quad \Pi_2^2 = \frac{\mathcal{C} - c_0}{c_2-c_0}
$$
confirming that $\Pi_2^2 = \mathbf{1} - \Pi_0^2$, and explicitly yielding
$$
\Pi_0^2 = \frac{QP + \frac{1}{4} (p-2)p - \frac{1}{4} p(p-2)}{\frac{1}{4} (p+2)p - \frac{1}{4} p(p-2)} = \frac{1}{p} \, QP
$$
and
$$
\Pi_2^2 = \frac{QP + \frac{1}{4} (p-2)p - \frac{1}{4} p(p+2)}{\frac{1}{4} (p-2)p - \frac{1}{4} p(p+2)} = \mathbf{1} - \frac{1}{p} \, QP
$$
which is in accordance with the expressions obtained by direct calculation. Next, let $r=3$ and $s=1$, then
$$
\Pi_1^3 = \frac{\mathcal{C} - c_3}{c_1 - c_3} = \frac{QP + \frac{1}{4} (p-3)(p-1) - \frac{1}{4} (p-3)(p-1)}{\frac{1}{4} (p+1)(p-1) - \frac{1}{4} (p-1)(p-3)} = \frac{1}{p-1} QP 
$$
confirming the direct calculations leading to (\ref{Pi12j+1}). However, it is clear that the more interesting examples are those where more than one factor appears in formula (\ref{casimir}). To this end let us take $r=4$ and $s=0$, then (\ref{casimir}) reads
$$
\Pi_0^4 = \frac{(\mathcal{C} - c_2)(\mathcal{C} - c_4)}{(c_0-c_2)(c_0-c_4)} = \frac{( QP - (p-2)) QP}{2p(p-1)} = \frac{QPQP - (p-2)QP}{2p(p-1)}
$$
Invoking Lemma \ref{lemmaS4} we obtain 
$$
QPQP = Q^2 P^2 + (p-2)QP
$$
whence
$$
\Pi_0^4 = \frac{1}{2p(p-1)} \, Q^2 P^2
$$
in accordance with the general expression (\ref{Pi02j}). As a final example we take $r=5$ and $s=1$; here we have
$$
\Pi_1^5 = \frac{(\mathcal{C} - c_3)(\mathcal{C} - c_5)}{(c_1-c_3)(c_1-c_5)} =  \frac{( QP - (p-3)) QP}{2(p-1)(p-2)} = \frac{QPQP - (p-3)QP}{2(p-1)(p-2)}
$$
leading, in a similar way as above, to
$$
\Pi_1^5 = \frac{1}{2(p-1)(p-2)} \, Q^2 P^2
$$
which once more confirms (\ref{Pi12j+1}).


\subsection{Action of the Witt basis vectors on the symplectic cells}


We conclude this section by investigating the action of the Witt basis vectors $\gf_j$, $\gfd_j$, $j=1,\ldots,p$, by left multiplication on the symplectic cells $\mS_s^r$. To that end the following lemma is readily proven.
\begin{lemma}
\label{lemmaS7}
For $j=1,\ldots,p$ one has that
\begin{enumerate}
\item[(i)] $[ P,\gf_j]=0 = [Q,\gfd_j]$;
\item[(ii)] $[Q^2,\gf_j] = 0 = [P^2,\gfd_j]$.
\end{enumerate}
\end{lemma}

We then obtain the following result regarding the action of the Witt basis vectors on the symplectic spinor cells.

\begin{proposition}
By left multiplication the Witt basis vectors act as follows:
\begin{enumerate}
\item[(i)] $\gfd_j : \mS_s^r \longrightarrow \mS_{s-1}^{r+1} \oplus \mS_{s+1}^{r+1}$;
\item[(ii)] $\gf_j : \mS_s^r \longrightarrow \mS_{s-1}^{r-1} \oplus \mS_{s+1}^{r-1}$.
\end{enumerate}
\end{proposition}

\pf
(i) We know that $\gfd_j: \mS^r \rightarrow \mS^{r+1}$, and so the statement is trivial for $s=0,1,2$. Take $k$ such that $s+2k=r$, then by definition
$$
\mS_s^r = \mS_s^{s+2k} = Q^k \mS_s^s
$$
with $\mS_s^s \subset \mbox{Ker} P$. Since $Q$ and $\gfd_j$ are commuting, we have
$$
\gfd_j \mS_s^r = \gfd_j Q^k \mS_s^s = Q^k \gfd_j \mS_s^s
$$
Since $\gfd_j \mS_s^s \subset \mS^{s+1}$, we can decompose $\gfd_j \mS_s^s$ into a direct sum
$$
\gfd_j \mS_s^s = \mT_{s+1}^{s+1} \oplus \mT_{s-1}^{s+1} \oplus \mT^{s+1}
$$ 
where $\mT_{s+1}^{s+1} \subset \mS_{s+1}^{s+1}$, $\mT_{s-1}^{s+1} \subset \mS_{s-1}^{s+1}$ and $\mT^{s+1} \subset \mS^{s+1} \setminus \bigl ( \mS_{s+1}^{s+1} \oplus \mS_{s-1}^{s+1} \bigr )$.
Since $P^2$ and $\gfd_j$ are commuting, we have on the one hand
$$
P^2 \left ( \gfd_j \mS_s^s \right ) = \gfd_j P^2 \mS^s_s = 0
$$
while on the other
$$
P^2 \left ( \gfd_j \mS_s^s \right ) = P^2 \mT_{s+1}^{s+1} \oplus P^2 \mT_{s-1}^{s+1} \oplus P^2 \mT^{s+1}
$$
But $P^2 \mT_{s+1}^{s+1} =0$ since $\mS_{s+1}^{s+1} \subset \mbox{Ker} P$, and $P^2 \mT_{s-1}^{s+1} = 0$ since $P^2 \mS_{s-1}^{s+1} = P \mS_{s-1}^{s-1} = 0$ in view of $\mS_{s-1}^{s-1} \subset \mbox{Ker} P$. It follows that $P^2 \mT^{s+1}=0$ and hence that $\mT^{s+1} = 0$. In this way we have obtained that
$$
\gfd_j \mS_s^r = Q^k \gfd_j \mS_s^r = Q^k \mT_{s+1}^{s+1} \oplus Q^k \mT_{s-1}^{s+1} 
$$
or
$$
\gfd_j \mS_s^r \subset \mS_{s+1}^{r+1} \oplus \mS_{s-1}^{r+1}
$$
(ii) The proof is similar to the one of (i), switching the roles of the operators $P$ and $Q$.
\qed

The above proposition allows for the following decomposition of the multiplicative operators $\gf_j$ and $\gfd_j$, $j=1,\ldots,p$:
$$
\gf_j \Bigr\rvert_{\mS_s^r} = (\gf_j)^r_{s-} +  (\gf_j)^r_{s+} \quad , \quad (\gf_j)^r_{s\mp} : \mS_s^r \longrightarrow \mS_{s\mp1}^{r-1}
$$
and
$$
\gfd_j \Bigr\rvert_{\mS_s^r} = (\gfd_j)^r_{s-} +  (\gfd_j)^r_{s+} \quad , \quad (\gfd_j)^r_{s\mp} : \mS_s^r \longrightarrow \mS_{s\mp1}^{r+1}
$$

The following properties of those Witt basis components are immediate; they originate from the well--known properties of the Witt basis vectors themselves (see Section 4.1). In order to make the formulae more transparent we switch to shorthand notations $\gf_{j-}$ and $\gfd_{j+}$, assuming the symplectic cells, on which they act, to be clear from the context.

\begin{lemma}
\label{lemmaS8}
For $j,k=1,\ldots,p$ one has
\begin{itemize}
\item[(i)] $\gf_{j-} \, \gf_{k-} + \gf_{k-} \, \gf_{j-} = 0$
\item[(ii)] $\gf_{j+} \, \gf_{k+} + \gf_{k+} \, \gf_{j+} = 0$
\item[(iii)] $\gf_{j-} \, \gf_{k+}  + \gf_{j+} \, \gf_{k-} + \gf_{k-} \, \gf_{j+}  + \gf_{k+} \, \gf_{j-}= 0$
\end{itemize}
and in particular
\begin{itemize}
\item[(iv)] $\gf_{j-} \, \gf_{j-} = 0$
\item[(v)] $\gf_{j+} \, \gf_{j+} = 0$
\item[(vi)] $\gf_{j+} \, \gf_{j-}  + \gf_{j-} \, \gf_{j+} = 0$
\end{itemize}
\end{lemma}

\begin{lemma}
\label{lemmaS9}
For $j,k=1,\ldots,p$ one has
\begin{itemize}
\item[(i)] $\gfd_{j-} \, \gfd_{k-} + \gfd_{k-} \, \gfd_{j-} = 0$
\item[(ii)] $\gfd_{j+} \, \gfd_{k+} + \gfd_{k+} \, \gfd_{j+} = 0$
\item[(iii)] $\gfd_{j-} \, \gfd_{k+}  + \gfd_{j+} \, \gfd_{k-} + \gfd_{k-} \, \gfd_{j+}  + \gfd_{k+} \, \gfd_{j-}= 0$
\end{itemize}
and in particular
\begin{itemize}
\item[(iv)] $\gfd_{j-} \, \gfd_{j-} = 0$
\item[(v)] $\gfd_{j+} \, \gfd_{j+} = 0$
\item[(vi)] $\gfd_{j+} \, \gfd_{j-}  + \gfd_{j-} \, \gfd_{j+} = 0$
\end{itemize}
\end{lemma}

\begin{lemma}
\label{lemmaS10}
For $j \neq k=1,\ldots,p$ one has
\begin{itemize}
\item[(i)] $\gf_{j-} \, \gfd_{k-} + \gfd_{k-} \, \gf_{j-} = 0$
\item[(ii)] $\gf_{j+} \, \gfd_{k+} + \gfd_{k+} \, \gf_{j+} = 0$
\item[(iii)] $\gf_{j-} \, \gfd_{k+}  + \gf_{j+} \, \gfd_{k-} + \gfd_{k-} \, \gf_{j+}  + \gfd_{k+} \, \gf_{j-}= 0$
\end{itemize}
\end{lemma}

\begin{lemma}
\label{lemmaS11}
For $j=1,\ldots,p$ one has
\begin{itemize}
\item[(i)] $\gf_{j-} \, \gfd_{j-} + \gfd_{j-} \, \gf_{j-}= 0$
\item[(ii)] $\gf_{j+} \, \gfd_{j+} + \gfd_{j+} \, \gf_{j+}= 0$
\item[(iii)] $\gf_{j+} \, \gfd_{j-}  + \gf_{j-} \, \gfd_{j+}  + \gfd_{j+} \, \gf_{j-}  + \gfd_{j-} \, \gf_{j+}= {\bf 1}$
\end{itemize}
\end{lemma}
\vspace*{3mm}

Now we formulate explicitly the action of the Witt vector components.

\begin{proposition}
For $j=1,\ldots,p$, $r=0,1,\ldots,p$ and $k=0,1,\ldots,\lfloor \frac{r}{2} \rfloor$ one has
\begin{eqnarray*}
(\gfd_j)^r_{(r-2k)-} &=& \frac{1}{\gamma^r_{r-2k}} \, Q^{k+1} \, P^{k+1} \, \gfd_j\\
(\gfd_j)^{2p-r}_{(r-2k)-} &=& \frac{1}{\gamma^{2p-r}_{r-2k}} \, P^{k} \, Q^{k} \, \gfd_j\\
(\gf_j)^r_{(r-2k)-} &=&  \frac{1}{\gamma^{2p-r}_{r-2k}} \, Q^{k} \, P^{k} \, \gf_j\\
(\gf_j)^{2p-r}_{(r-2k)-} &=& \frac{1}{\gamma^{r}_{r-2k}} \, P^{k+1} \, Q^{k+1} \, \gf_j\\
\end{eqnarray*}
where
\begin{eqnarray*}
\gamma^r_{r-2k} &=& \alpha^{0}_{r-2k-1} \, \alpha^{1}_{r-2k-1}  \cdots  \alpha^{k}_{r-2k-1}\\[2mm]
&=& (k+1)!(p-r+k+1)(p-r+k+2)\cdots(p-r+2k+1)
\end{eqnarray*}
and
\begin{eqnarray*}
\gamma^{2p-r}_{r-2k} &=& \alpha^{0}_{r-2k-1} \, \alpha^{1}_{r-2k-1} \cdots  \alpha^{k-1}_{r-2k-1}  \\[2mm]
&=&(k)!(p-r+k+2)\cdots (p-r+2k)(p-r+2k+1)
\end{eqnarray*}
\end{proposition}

\pf
For the action of $\gfd_j$ on the symplectic cell $\mS_{r-2k}^{r}$ we consecutively have, in view of Lemma \ref{lemmaS1} and Lemma \ref{lemmaS5}, and making use of the shorthand notation $\gfd_{j-} = (\gfd_j)_{r-2k}^{r}$,
\begin{eqnarray*}
 Q^{k+1} \, P^{k+1} \, \gfd_j &=&  Q^{k+1} \, P^{k+1} \, ( \gfd_{j-} + \gfd_{j+}) =  Q^{k+1} \, P^{k+1} \, \gfd_{j-}\\
&=&  Q^{k} \, (Q \, P) \, P^{k} \, \gfd_{j-}\\ 
&=& Q^{k} \, (\alpha_{r-2k-1}^{0} ) \, P^{k} \, \gfd_{j-}\\
&=&  (\alpha_{r-2k-1}^{0} ) \, \, Q^{k-1} \, (\alpha_{r-2k-1}^{1})  \,  P^{k-1} \, \gfd_{j-}\\
&=& \ldots\\
&=&  \alpha^{0}_{r-2k-1} \, \alpha^{1}_{r-2k-1}  \cdots  \alpha^{k}_{r-2k-1}\\
&=& \gamma^r_{r-2k} \, \gfd_{j-}
\end{eqnarray*}
The expressions for the other three actions are proven in a similar way.
\eop

\begin{corollary}
For $j=1,\ldots,p$, $r=0,1,\ldots,p$ and $k=0,1,\ldots,\lfloor \frac{r}{2} \rfloor$ one has
\begin{eqnarray*}
(\gfd_j)^r_{(r-2k)+} &=& \left( 1  - \frac{1}{\gamma^r_{r-2k}} \, Q^{k+1} \, P^{k+1} \right) \, \gfd_j\\
(\gfd_j)^{2p-r}_{(r-2k)+} &=&  \left( 1 -  \frac{1}{\gamma^{2p-r}_{r-2k}} \, P^{k} \, Q^{k} \right) \, \gfd_j\\
(\gf_j)^r_{(r-2k)+} &=&  \left( 1 -    \frac{1}{\gamma^{2p-r}_{r-2k}} \, Q^{k} \, P^{k} \right) \, \gf_j\\
(\gf_j)^{2p-r}_{(r-2k)+} &=&  \left( 1 -  \frac{1}{\gamma^{r}_{r-2k}} \, P^{k+1} \, Q^{k+1} \right) \, \gf_j\\
\end{eqnarray*}
\end{corollary}


\subsection{Alternative Proof of Theorem \ref{symplecticcells}}


We conclude this section on spinor space by proving Theorem \ref{symplecticcells} in an alternative way; let us first recall this theorem:
it states that the homogeneous spinor spaces $\mS^{r}$ and $\mS^{2p-r}$, $r=0,\ldots,p$, may be decomposed into Sp$(p)$--irreducibles as
$$
\mS^{r} = \bigoplus_{j=0}^{\lfloor \frac{r}{2} \rfloor} \mS_{r-2j}^{r} \quad \mbox{ and} \quad
\mS^{2p-r} = \bigoplus_{j=0}^{\lfloor \frac{r}{2} \rfloor} \mS_{r-2j}^{2p-r} 
$$
We will prove the decomposition of $\mS^{2p-r}$, $r=0,\ldots,p$ through a series of lemmata, the one of  $\mS^{r}$ being completely similar by symmetry.

\begin{lemma}
\label{lemmaS12}
For $r=0,\ldots,2p$ one has
$$
\mS^{2p-r} \simeq \mbox{$\bigwedge^r$} \, \mC^{2p}
$$
\end{lemma}

\pf
We know that $\mS^{2p-r}$ is an irreducible $\gsl_{2p}(\mC)$--module. Recall from Section 4.1 that the Cartan subalgebra of $\gsl_{2p}(\mC)$ is given, in terms of bivectors, by
$$
\gh = \bigl \{ \gf_{j+1} \gfd_{j+1} - \gf_j \gfd_j : j=1,\ldots,2p-1 \bigr \}
$$
or, equivalently, by
$$
\gh = \bigl \{ H_j = I_j - I_{2p} : j=1,\ldots,2p-1 \bigr \}
$$
For $r=0$ we have $\mS^{2p} = \mbox{span}_{\mC} \bigl \{ \gfd_1 \gfd_2 \cdots \gfd_{2p} I \bigr \}$, and as
$$
H_j [ \gfd_1 \gfd_2 \cdots \gfd_{2p} I  ] = 0, \quad j=1,\ldots,2p-1
$$
its highest weight is $(0,\ldots,0)$ with length $2p-1$, which corresponds to $\mC \simeq \bigwedge^0 \mC^{2p}$. Similarly, for $r=1$, we have $\mS^{2p-1} = \mbox{span}_{\mC} \bigl \{ \gfd_1 \gfd_2 \cdots \widehat{\gfd_j} \cdots \gfd_{2p} I , \, j=1,\ldots,2p \bigr \}$ where the notation $\gfd_1 \gfd_2 \cdots \widehat{\gfd_j} \cdots \gfd_{2p} I$ quite naturally expresses the fact that $\gfd_j$ has been omitted from this product of Witt basis vectors. As
$$
H_1 [ \gfd_2 \gfd_3 \cdots \gfd_{2p} I  ] = \gfd_2 \gfd_3 \cdots \gfd_{2p} I 
$$
while
$$
H_j [ \gfd_2 \gfd_3 \cdots \gfd_{2p} I  ] = 0, \quad j=2,3,\ldots,2p
$$
its highest weight is $(1,0,\ldots,0)$ with length $2p-1$, which corresponds to $\mC^{2p} \simeq \bigwedge^1 \mC^{2p}$. Continuing in the same way we arrive at the case where $r=2p-1$ and $\mS^1 = \mbox{span}_{\mC} \bigl \{ \gfd_j I , \, j=1,\ldots,2p \bigr \}$. Seen the fact that 
$$
H_1[ \gfd_{2p} I] = H_2 [ \gfd_{2p} I ] = \ldots = H_{2p-1}[\gfd_{2p} I] = \gfd_{2p} I
$$
its highest weight is $(1,1,\ldots,1)$ with length $2p-1$, corresponding to $\bigwedge^{2p-1} \mC^{2p}$. Finally, for $r=2p$, there holds that $\mS^0 = \mbox{span}_{\mC} \bigl \{ I \bigr \}$ and
$$
H_j [I] = 0, \qquad j=1,\ldots,2p
$$
leading to the highest weight $(0,0,\ldots,0)$ of length $2p-1$, corresponding to $\mC \simeq \bigwedge^{2p} \mC^{2p}$.
\qed

\begin{lemma}
\label{lemmaS13}
Under the action of $\gsl_p(\mC)$ one has
$$
\mC^{2p} = V \oplus \overline{V}
$$
where $V$ is a fundamental representation for $\gsl_p(\mC)$ with highest weight $(1,0,\ldots,0)$ of length $p-1$, and $\overline{V}$ is its dual with highest weight $(1,1,\ldots,1)$ of length $p-1$.
\end{lemma}

\pf
Recall from the proof of Lemma \ref{lemmaS12} that
$$
\mC^{2p} \simeq \mbox{$\bigwedge^1$} \mC^{2p} \simeq \mS^{2p-1} = \mbox{span}_{\mC} \bigl \{ \gfd_1 \gfd_2 \cdots \widehat{\gfd_j} \cdots \gfd_{2p} I, \, j=1,\ldots,2p \bigr \}
$$
and thus also
$$
\mC^{2p} \simeq \mbox{span}_{\mC} \bigl \{ \gfd_1  \cdots \widehat{\gfd_{2k}} \cdots \gfd_{2p} I, \, k=1,\ldots,p \bigr \} \oplus \mbox{span}_{\mC} \bigl \{ \gfd_1  \cdots \widehat{\gfd_{2k-1}} \cdots \gfd_{2p} I, \, k=1,\ldots,p \bigr \}
$$
We put
$$
V = \mbox{span}_{\mC} \bigl \{ \gfd_1  \cdots \widehat{\gfd_{2k}} \cdots \gfd_{2p} I, \, k=1,\ldots,p \bigr \}
$$
and
$$
\overline{V} = \mbox{span}_{\mC} \bigl \{ \gfd_1  \cdots \widehat{\gfd_{2k-1}} \cdots \gfd_{2p} I, \, k=1,\ldots,p \bigr \}
$$
Recall from Section 4.2 that the Cartan subalgebra of $\gsl_p(\mC)$ is given, in terms of bivectors, by
$$
\gh^{{\rm sl}} = \{ H_j^{{\rm sl}} = H_j^{{\rm sympl}} - H_p^{{\rm sympl}} : j = 1,\ldots,p-1 \}
$$
In $V$ we have the highest weight vector $\gfd_1 \gfd_3 \cdots \gfd_{2p} I$ for which indeed
\begin{eqnarray*}
H_1^{{\rm sl}} \bigl [ \gfd_1 \gfd_3 \cdots \gfd_{2p} I \bigr ] & = & \bigr ( H_1^{{\rm sympl}} - H_p^{{\rm sympl}} \bigr ) \bigr [ \gfd_1 \gfd_3 \cdots \gfd_{2p} I \bigr ] \\
& = & \bigl ( \gf_2 \gfd_2 - \gf_1 \gfd_1 - \gf_{2p} \gfd_{2p} + \gf_{2p-1} \gfd_{2p-1} \bigr ) \bigl [ \gfd_1 \gfd_3 \cdots \gfd_{2p} I \bigr ] \\ 
& = & \gfd_1 \gfd_3 \cdots \gfd_{2p} I
\end{eqnarray*}
while for $j=2,\ldots,p-1$ there holds
\begin{eqnarray*}
H_j^{{\rm sl}} \bigl [ \gfd_1 \gfd_3 \cdots \gfd_{2p} I \bigr ] & = & \bigl ( H_j^{{\rm sympl}} - H_p^{{\rm sympl}} \bigr ) \bigl [ \gfd_1 \gfd_3 \cdots \gfd_{2p} I \bigr ] \\
& = & \bigl ( \gf_{2j} \gfd_{2j} - \gf_{2j-1} \gfd_{2j-1} - \gf_{2p} \gfd_{2p} + \gf_{2p-1} \gfd_{2p-1} \bigr ) \bigl [ \gfd_1 \gfd_3 \cdots \gfd_{2p} I \bigr ] \\ 
& = & 0
\end{eqnarray*}
which means that its highest weight is $(1,0,\ldots,0)$ of length $p-1$, which corresponds to $\mC^p$.\\
In $\overline{V}$ we have the highest weight vector $\gfd_1 \gfd_2 \cdots \gfd_{2p-2} \gfd_{2p} I$ for which indeed, for all $j=1,\ldots,p-1$:
\begin{eqnarray*}
H_j^{{\rm sl}} \bigl [ \gfd_1 \gfd_2 \cdots \gfd_{2p-2} \gfd_{2p} I \bigr ] & = & \bigl ( H_j^{{\rm sympl}} - H_p^{{\rm sympl}} \bigr ) \bigl [ \gfd_1 \gfd_2 \cdots \gfd_{2p-2} \gfd_{2p} I \bigr ] \\
& = & \bigl ( \gf_{2j} \gfd_{2j} - \gf_{2j-1} \gfd_{2j-1} - \gf_{2p} \gfd_{2p} + \gf_{2p-1} \gfd_{2p-1} \bigr ) \bigl [ \gfd_1 \gfd_2 \cdots \gfd_{2p-2} \gfd_{2p} I \bigr ] \\ 
& = & \gfd_1 \gfd_2 \cdots \gfd_{2p-2} \gfd_{2p} I
\end{eqnarray*}
This means that its highest weight is $(1,1,\ldots,1)$ of length $p-1$, corresponding with $\bigwedge^1 \mC^p \simeq \mC^p$.
\qed

Combining the results of the lemmata \ref{lemmaS12} and \ref{lemmaS13} we obtain,  for $r = 0, \ldots,2p$, the isomorphism 
$$
\mS^{2p-r} \simeq \mbox{$\bigwedge^r$} \mC^{2p} = \bigoplus_{a+b=r} \bigl ( \mbox{$\bigwedge^a$} V \otimes \mbox{$\bigwedge^b$} \overline{V} \bigr )
$$
clearly showing that $\mS^{2p-r}$ is not $\gsl_p(\mC)$--irreducible. Instead, by Pieri's formula (see  \cite{fulhar}) on the tensor product of $\gsl_p(\mC)$--modules, it can be decomposed into the direct sum
$$
\mS^{2p-r} \simeq \bigoplus_{j=0}^{\lfloor \frac{r}{2} \rfloor} \; \bigoplus_{a+b=r-2j} \bigl ( \mbox{$\bigwedge^a$} V \boxtimes \mbox{$\bigwedge^b$} \overline{V} \bigr )
$$
where $\boxtimes$ denotes the Cartan product. Each of the terms $\bigwedge^a V \boxtimes \bigwedge^b \overline{V}$ appearing in this decomposition is now characterized as an $\gsl_p(\mC)$--module.

\begin{lemma}
\label{lemmaS14}
One has
\begin{enumerate}
\item[(i)] the space
$$
\mbox{$ \bigwedge^a$} V \simeq ( \underbrace{1,\ldots,1}_a, \underbrace{0,\ldots,0}_{p-a-1})
$$ 
is an irreducible $\gsl_p(\mC)$--module generated by the highest weight vector $\gfd_{2p} \gfd_{2p-1} \cdots \gfd_{2a+2} I$;
\item[(ii)] the space
$$
\mbox{$ \bigwedge^b$} \overline{V} \simeq ( \underbrace{1,\ldots,1}_{p-b}, \underbrace{0,\ldots,0}_{b-1})
$$ 
is an irreducible $\gsl_p(\mC)$--module generated by the highest weight vector $\gfd_{1} \gfd_{3} \cdots \gfd_{2(p-b)+1} I$.
\end{enumerate}
\end{lemma}

\pf
First note that these realizations of $\bigwedge^a V$ and $\bigwedge^b \overline{V}$ differ from those in the proof of Lemma \ref{lemmaS12}, as we will now use only even, respectively odd labeled Witt basis vectors.\\[-2mm]

\noindent (i) We have
\begin{eqnarray*}
H_1^{{\rm sl}}\bigl [ \gfd_{2p} \gfd_{2p-2} \cdots \gfd_{2a+2} I \bigr ] & = & \bigl ( \gf_2 \gfd_2 - \gf_1 \gfd_1 - \gf_{2p} \gfd_{2p} + \gf_{2p-1} \gfd_{2p-1} \bigr ) \bigl [ \gfd_{2p} \gfd_{2p-2} \cdots \gfd_{2a+2} I \bigr ] \\
& = & \gfd_{2p} \gfd_{2p-2} \cdots \gfd_{2a+2} I - \gfd_{2p} \gfd_{2p-2} \cdots \gfd_{2a+2} I - 0 + \gfd_{2p} \gfd_{2p-2} \cdots \gfd_{2a+2} I \\
& = & \gfd_{2p} \gfd_{2p-2} \cdots \gfd_{2a+2} I 
\end{eqnarray*}
and similarly, for all $j \leq a$:
$$
H_j^{{\rm sl}}\bigl [ \gfd_{2p} \gfd_{2p-2} \cdots \gfd_{2a+2} I \bigr ] = \gfd_{2p} \gfd_{2p-2} \cdots \gfd_{2a+2} I
$$
On the contrary, for $j>a$ we obtain
$$
H_j^{{\rm sl}}\bigl [ \gfd_{2p} \gfd_{2p-2} \cdots \gfd_{2a+2} I \bigr ]=0
$$

\noindent (ii) We also have
\begin{eqnarray*}
H_1^{{\rm sl}}\bigl [ \gfd_1 \gfd_3 \cdots \gfd_{2p-2b+1} I \bigr ] & = & \bigl ( \gf_2 \gfd_2 - \gf_1 \gfd_1 - \gf_{2p} \gfd_{2p} + \gf_{2p-1} \gfd_{2p-1} \bigr ) \bigl [ \gfd_1 \gfd_3 \cdots \gfd_{2p-2b+1} I \bigr ] \\
& = & \gfd_1 \gfd_3 \cdots \gfd_{2p-2b+1} I  - \gfd_1 \gfd_3 \cdots \gfd_{2p-2b+1} I   + \gfd_1 \gfd_3 \cdots \gfd_{2p-2b+1} I  \\
& = & \gfd_1 \gfd_3 \cdots \gfd_{2p-2b+1} I 
\end{eqnarray*}
and similarly, for all $j \leq p-b$:
$$
H_j^{{\rm sl}}\bigl [ \gfd_1 \gfd_3 \cdots \gfd_{2p-2b+1} I \bigr ] = \gfd_1 \gfd_3 \cdots \gfd_{2p-2b+1} I
$$
while for $j > p-b$ there holds
$$
H_j^{{\rm sl}}\bigl [ \gfd_1 \gfd_3 \cdots \gfd_{2p-2b+1} I \bigr ] = 0
$$
\qed

\begin{lemma}
\label{lemmaS15}
The Cartan product $\bigwedge^a V \boxtimes \bigwedge^b \overline{V}$ is an irreducible $\gsl_p(\mC)$--module generated by the highest weight vector $\alpha = \gfd_{2p} \gfd_{2p-2} \cdots \gfd_{2a+2} \gfd_1 \gfd_3 \cdots \gfd_{2p-2b+1} I$.
\end{lemma}

\pf
We distinguish three cases. If $a+b<p$ then
\begin{eqnarray*}
H_j^{{\rm sl}}\bigl [ \alpha \bigr ] & = & 2 \alpha, \qquad \mbox{for } j=1,\ldots,a \\
H_j^{{\rm sl}}\bigl [ \alpha \bigr ] & = & \phantom{2} \alpha,  \qquad \mbox{for } j=a+1,\ldots,p-b \\
H_j^{{\rm sl}}\bigl [ \alpha \bigr ] & = & \phantom{2} 0,  \qquad \mbox{for } j=p-b+1,\ldots,p-1
\end{eqnarray*}
If $a+b=p$ then
\begin{eqnarray*}
H_j^{{\rm sl}}\bigl [ \alpha \bigr ] & = & 2 \alpha, \qquad \mbox{for } j=1,\ldots,a=p-b \\
H_j^{{\rm sl}}\bigl [ \alpha \bigr ] & = & \phantom{2} 0,  \qquad \mbox{for } j=a+1=p-b+1,\ldots,p-1
\end{eqnarray*}
If $a+b > p$ then
\begin{eqnarray*}
H_j^{{\rm sl}}\bigl [ \alpha \bigr ] & = & 2 \alpha, \qquad \mbox{for } j=1,\ldots,p-b \\
H_j^{{\rm sl}}\bigl [ \alpha \bigr ] & = & \phantom{2} \alpha,  \qquad \mbox{for } j=p-b+1,\ldots,a \\
H_j^{{\rm sl}}\bigl [ \alpha \bigr ] & = & \phantom{2} 0,  \qquad \mbox{for } j=a+1,\ldots,p-1
\end{eqnarray*}
In all three cases we notice that the highest weight of $\alpha$ is the sum of the highest weights for $\bigwedge^a V$ and $\bigwedge^b \overline{V}$. By the definition of the Cartan product this means that the highest weight of $\alpha$ is precisely the highest weight of $\bigwedge^a V \boxtimes \bigwedge^b \overline{V}$.
\qed

\begin{lemma}
\label{lemmaS16}
For $r=0,\ldots,p$ and whenever $a+b=r$ one has
$$
\bigwedge\nolimits^a V \ \boxtimes \ \bigwedge\nolimits^b \overline{V} \ \subset \ \mS^{2p-r} \cap \mbox{Ker} \, Q \, = \, \mS_r^{2p-r}
$$
\end{lemma}

\pf
From Lemma \ref{lemmaS15} we know that $\bigwedge^a V \boxtimes \bigwedge^b \overline{V}$ is generated by the highest weight vector 
$$
\gfd_{2p} \gfd_{2p-2} \cdots \gfd_{2a+2} \gfd_1 \gfd_3 \cdots \gfd_{2p-2b+1} I
$$
where the number of Witt basis vectors is $(p-a)+(p-b) = 2p-(a+b)$. Whenever $a+b=r$ this number equals $2p-r$, so, in those cases, this highest weight vector definitely belongs to $\mS^{2p-r}$. Seen the fact that the operator $Q$ is given by $\sum_{j=1}^p \gfd_{2j} \gfd_{2j-1}$, the remaining question to be answered is: does the highest weight vector contain, for each $j=1,\ldots,p$, either the vector $\gfd_{2j}$ or the vector $\gfd_{2j-1}$, so that it will indeed be annihilated by $Q$. When rewriting the highest weight vector as
$$
\gfd_{2p} \gfd_{2p-2} \cdots \gfd_{2r-2b+2} \gfd_1 \gfd_3 \cdots \gfd_{2p-2b+1} I
$$ 
the answer to this question easily is seen to be affirmative. Indeed, looking at the Witt basis vectors with odd indices, it only contains $\gfd_1$ up to $\gfd_{2p-2b-1}$, so it does not contain $\gfd_{2p-2b+1}$, $\gfd_{2p-2b+3}$ up to $\gfd_{2p-1}$. However, at the same time, looking at the even indexed ones, it does contain $\gfd_{2p-2b+2}$ up to $\gfd_{2p}$, which do exactly complement the lacking odd indexed ones in the definition of $Q$.
\qed

Now, taking $a=r$ and $b=0$ in the result of Lemma \ref{lemmaS16}, we know that $\bigwedge^r V \boxtimes \bigwedge^0 \overline{V}$ belongs to $\mS_r^{2p-r}$ and contains the highest weight vector $\gfd_{2p} \gfd_{2p-2} \ldots \gfd_{2r+2} \gfd_1 \gfd_3 \ldots \gfd_{2p-1} I$. This highest weight vector generates an $\gsp_{2p}(\mC)$--module $M$ with highest weight 
$$
(\underbrace{1,\ldots,1}_{r}, \underbrace{0,\ldots,0}_{p-r})
$$
Moreover $M \subset \mS_r^{2p-r}$ since $\mS_r^{2p-r}$ is $\gsp_{2p}(\mC)$--invariant. We will now prove that $M=\mS_r^{2p-r}$, by calculating its dimension.

\begin{lemma}
\label{lemmaS17}
For $r=0,\ldots,p$ it holds that ${\rm dim} (M) = {\rm \dim} (\mS_r^{2p-r})$.
\end{lemma}

\pf
We will calculate the dimension of the $\gsp_{2p}(\mC)$--module $M$ with highest weight 
$$
(\underbrace{1,\ldots,1}_{r}, \underbrace{0,\ldots,0}_{p-r})
$$
by means of the Weyl dimension formula (see \cite{fulhar}), which states that the dimension of the irreducible representation $\pi$ of the Lie algebra $\geg$ with highest weight $\mu$ is given by
$$
\mbox{dim} (\pi) = \prod_{\alpha \in R^+} \frac{\langle \alpha, \mu+\delta \rangle}{\langle \alpha,\delta \rangle}
$$
where $\delta$ is half the sum of the positive roots $\alpha \in R^+$. For $\geg = \gsp_{2p}(\mC)$ the positive root system is given by
$$
R^+ = \bigl \{ 2 L_i, L_i+L_j,L_i-L_j : i<j=1,\ldots,p \bigr \}
$$
and the inner product is the Euclidean one: $\langle L_i,L_j \rangle = \delta_{ij}$. We find
\begin{eqnarray*}
\delta & = & \frac{1}{2} \bigl ( \sum_{i=1}^p 2 L_i + \sum_{i<j} (L_i+L_j) + \sum_{i<j} (L_i-L_j) \bigr )\\
& = & \frac{1}{2} \bigl ( 2 \sum_{i=1}^p L_i + 2 \sum_{i<j} L_i \bigr ) \\
& = & p L_1 + (p-1) L_2 + \ldots + 2 L_{p-1} + L_p \ = \ (p,p-1,p-2,\ldots,2,1)
\end{eqnarray*}
For the denominators in the Weyl dimension formula, we find that
\begin{enumerate}
\item[(i)] for $\alpha=2L_i$, $i=1,\ldots,p$: $\langle \alpha,\delta\rangle = 2(p-i+1)$;
\item[(ii)]  for $\alpha=L_i+L_j$, $i<j=1,\ldots,p$: $\langle \alpha,\delta\rangle = 2p+2 -(i+j)$;
\item[(iii)] for $\alpha = L_i - L_j$, $i<j=1,\ldots,p$: $\langle \alpha,\delta\rangle = j-i$.
\end{enumerate}
For the numerators we find, with
\begin{eqnarray*}
\mu + \delta &=& (p+1,p,p-1,\ldots,p-r+2,p-r,\ldots,2,1)\\[-2mm]
& & \hspace*{36mm} \downarrow \hspace*{10mm} \downarrow \\[-3mm]
& & \hspace*{30mm} \mbox{\tiny $r$th place} \hspace*{2mm} \mbox{\tiny (r+1)th place}
\end{eqnarray*}

(i) for $\alpha=2L_i$, $i=1,\ldots,p$:
\begin{eqnarray*}
\langle \alpha,\mu+\delta \rangle &=& 2(p-i+2), \quad i=1,\ldots,r \\
\langle \alpha,\mu+\delta \rangle &=& 2(p-i+1), \quad i=r+1,\ldots,p 
\end{eqnarray*}
(ii) for $\alpha = L_i+L_j$, $i<j=1,\ldots,p$:
\begin{eqnarray*}
\langle \alpha,\mu+\delta \rangle &=& 2p+4-(i+j), \quad i<j=1,\ldots,r \\
\langle \alpha,\mu+\delta \rangle &=& 2p+3-(i+j), \quad i=1,\ldots,r \mbox{\ and\ } j=r+1,\ldots,p \\
\langle \alpha,\mu+\delta \rangle &=& 2p+2-(i+j), \quad i<j=r+1,\ldots,p 
\end{eqnarray*}
(iii) for  $\alpha = L_i-L_j$, $i<j=1,\ldots,p$:
\begin{eqnarray*}
\langle \alpha,\mu+\delta \rangle &=& j-i, \quad i<j=1,\ldots,r \\
\langle \alpha,\mu+\delta \rangle &=& j-i+1, \quad i=1,\ldots,r \mbox{\ and\ } j=r+1,\ldots,p \\
\langle \alpha,\mu+\delta \rangle &=& j-i, \quad i<j=r+1,\ldots,p 
\end{eqnarray*}
For the positive roots $\alpha=2L_i$, $i=1,\ldots,p$ we thus find
$$
\prod_{\alpha =2L_i} \frac{\langle \alpha, \mu+\delta \rangle}{\langle \alpha,\delta \rangle} = 
\frac{2^p (p+1) p (p-1) \ldots (p-r+2) (p-r) (p-r-1) \ldots 2}{2^p p (p-1)(p-2) \ldots (p-r+1) (p-r) (p-r-1) \ldots 2} = \frac{p+1}{p-r+1}
$$
For the positive roots $\alpha=L_i+L_j$, $i<j=1,\ldots,p$ we find, after some calculation
$$
\prod_{\alpha =L_i+L_j} \frac{\langle \alpha, \mu+\delta \rangle}{\langle \alpha,\delta \rangle} =  
\frac{(2p+1) (2p) (2p-1) \ldots (2p-r+3) (2p-2r+2)}{(p+1) p (p-1) \ldots (p-r+3) (p-r+2)}
$$
For the positive roots $\alpha=L_i-L_j$, $i<j=1,\ldots,p$ we find, again after some calculation
$$
\prod_{\alpha =L_i-L_j} \frac{\langle \alpha, \mu+\delta \rangle}{\langle \alpha,\delta \rangle} =  
\frac{p (p-1)\ldots (p-r+1)}{r (r-1) \ldots 2}
$$
This leads to
\begin{eqnarray*}
\mbox{dim}(M) &=& \frac{p+1}{p-r+1} \, \frac{(2p+1) (2p) \ldots (2p-r+3) (2p-2r+2)}{(p+1) p (p-1) \ldots (p-r+3) (p-r+2)} \, \frac{p (p-1)\ldots (p-r+1)}{r (r-1) \ldots 2}\\
&=& \frac{2}{r!} \, (2p+1) (2p) (2p-1)\ldots (2p+3-r) (p-r+1) 
\end{eqnarray*}
On the other hand, we know from Corollary \ref{corollaryS1} that
\begin{eqnarray*}
\mbox{dim}(\mS_r^{2p-r} ) & = & \binom{2p}{r} - \binom{2p}{r-2} \ = \ \frac{(2p)!}{r! \, (2p-r+2)!} \, \bigl [ (2p-r+1)(2p-r+2) - (r-1)r \bigr ] \\
& = & 2 \, \frac{(2p)!}{r! \, (2p-r+2)!} \, (2p+1)(p-r+1) \\
& = & \frac{2}{r!} \, (2p+1) (2p) (2p-1)\ldots (2p+3-r) (p-r+1)  \ = \ \mbox{dim}(M) 
\end{eqnarray*}
\qed

\begin{lemma}
\label{lemmaS18}
For $r=0,\ldots,p$ and all $j = 0,\ldots, \lfloor \frac{r}{2} \rfloor$, $\mS_{r-2j}^{2p-r}$ is an irreducible $\gsp_{2p}(\mC)$--module.
\end{lemma}

\pf
We know by Lemma \ref{lemmaS5} that the (left) multiplication operators $P$ and $Q$ are vector space isomorphisms between the spaces $\mS_{r-2j}^{2p-r}$, $j=0,\ldots, \lfloor \frac{r}{2} \rfloor$. Seen the fact that $P$ and $Q$ are $\gsp_{2p}(\mC)$--invariant operators, see Lemma \ref{lemmaS2}, it follows that $P$ and $Q$ then also are $\gsp_{2p}(\mC)$--representation isomorphisms between the spaces $\mS_{r-2j}^{2p-r}$. But $\mS_r^{2p-r}$ is an irreducible $\gsp_{2p}(\mC)$--module since it coincides with the irreducible $\gsp_{2p}(\mC)$--module $M$ constructed above. It follows that all spaces $\mS_{r-2j}^{2p-r}$, $j = 0,\ldots, \lfloor \frac{r}{2} \rfloor$, are irreducible $\gsp_{2p}(\mC)$--modules.
\qed


\section{Euclidean, hermitian and quaternionic Clifford analysis}


As already mentioned in the introduction (Section 1), the central notion in standard Clifford analysis  is that of a {\em monogenic} function. This is a continuously differentiable function defined in an open region of Euclidean space $\mR^m$, taking its values in the Clifford algebra $\mR_{0,m}$, or subspaces thereof, and vanishing under the action of the Dirac operator $\dirac = \sum_{\alpha=1}^m e_{\alpha}\, \p_{X_{\alpha}}$. This Dirac operator thus is a vector valued first order differential operator, which can be seen as the Fourier or Fischer dual of the Clifford variable $\uX$. The notion of monogenicity is the higher dimensional counterpart of holomorphy in the complex plane. As the Dirac operator factorizes the Laplacian: $\Delta_m = - \dirac^2$, Clifford analysis can be regarded as a refinement of harmonic analysis. It is important to note that the Dirac operator is invariant under the action of the SO$(m)$--group, and also of the O$(m)$--group and the conformal group, whence this framework is usually referred to as Euclidean (or orthogonal) Clifford analysis.\\[-2mm]

Taking the dimension of the underlying Euclidean vector space $\mR^m$ to be even: $m=2n$, renaming the variables:
$$
(X_1,\ldots,X_{2n}) = (x_1,y_1,x_2,y_2,\ldots,x_n,y_n)
$$
and considering the complex structure $\mI_{2n}$ (see Section 3), we define the twisted vector variable
\begin{eqnarray*}
\uX| = \mI_{2n}[\uX] & = & \mI_{2n} \left[ \sum_{k=1}^n (e_{2k-1} x_k + e_{2k} y_k ) \right]\\
& = & \sum_{k=1}^n \mI_{2n} [e_{2k-1}] x_k + \mI_{2n} [e_{2k}] y_k \ = \ \sum_{k=1}^n (-y_k e_{2k-1} + x_k e_{2k})
\end{eqnarray*}
and, correspondingly, the twisted Dirac operator
$$
\dirac| = \mI_{2n} [ \dirac] = \sum_{k=1}^n ( - \p_{y_k} e_{2k-1} + \p_{x_k} e_{2k})
$$

A differentiable function $F$ then is called {\em hermitian monogenic} in some region $\Omega$ of $\mR^{2n}$, if and only if in that region $F$ is a solution of the system
\begin{equation}
\dirac F = 0 = \dirac| F
\label{hmon}
\end{equation}
Observe that this notion of hermitian monogenicity does not involve the use of complex numbers, but instead, could be developed as a real function theory. There is however an alternative approach to the concept of hermitian monogenicity, making use of the projection operators $\frac{1}{2} ({\bf 1} \pm i \, \mI_{2n})$ (see Section 4) and thus involving a complexification.\\[-2mm]

In this approach we define the vector variables
\begin{eqnarray*}
\uz & = & - \frac{1}{2} ({\bf 1} - i \, \mI_{2n} ) [\uX] \ = \ - \sum_{k=1}^n x_k \frac{1}{2} ({\bf 1} - i \, \mI_{2n}) [e_{2k-1}] - \sum_{k=1}^n y_k \frac{1}{2} ({\bf 1} - i \, \mI_{2n} ) [ e_{2k} ] \\
&=& \sum_{k=1}^n ( x_k \gf_k + y_k (i \gf_k) ) \ = \ \sum_{k=1}^n (x_k + i y_k) \gf_k = \sum_{k=1}^n z_k \gf_k
\end{eqnarray*}
and
$$
\uzd = \frac{1}{2} ({\bf 1} + i \, \mI_{2n} ) [\uX] = \sum_{k=1}^n ( x_k \gfd_k + y_k (-i \gfd_k) ) = \sum_{k=1}^n \overline{z}_k \gfd_k
$$
and, correspondingly, the hermitian Dirac operators
$$
2 \upzd =  - \frac{1}{2} ( {\bf 1} - i \, \mI_{2n})[\dirac] \ = \ \sum_{k=1}^n ( \gf_k \p_{x_k} + i \gf_k \p_{y_k}) = \sum_{k=1}^n \gf_k ( \p_{x_k} + i \p_{y_k} ) = 2 \sum_{k=1}^n \p_{\overline{z}_k} \gf_k
$$
and
$$
2 \upz =  \frac{1}{2} ( {\bf 1} + i \, \mI_{2n})[\dirac] \ = \ \sum_{k=1}^n ( \gfd_k \p_{x_k} - i \gfd_k \p_{y_k}) = \sum_{k=1}^n \gfd_k ( \p_{x_k} - i \p_{y_k} ) = 2 \sum_{k=1}^n \p_{z_k} \gfd_k
$$
As $2(\upz - \upzd) = \dirac$ and $2(\upz+\upzd) = i \, \mI_{2n}[\dirac] = i \, \dirac|$, it follows that the system (\ref{hmon}) is equivalent with the system
\begin{equation}
\upz F = 0 = \upzd F
\label{hmoneq}
\end{equation}

When decomposing  a spinor valued function $F: \mC^{2n} \longrightarrow \mS$ in its components according to the homogeneous subspaces of spinor space:
$$
F = \sum_{r=0}^n F^r, \quad F^r: \mC^{2n} \longrightarrow \mS^r, \quad r=0,\ldots,n
$$
the monogenicity of $F$ does {\em not} imply the monogenicity of the components $F^r$, $r=0,\ldots,n$, but the hermitian monogenicity of $F$ does imply their hermitian monogenicity, and vice versa. This is due to the nature of the action of the Witt basis vectors as (left) multiplication operators, implying that
$$
\upz F^r: \mC^{2n} \longrightarrow \mS^{r+1} \quad \mbox{and} \quad \upzd F^r: \mC^{2n} \longrightarrow \mS^{r-1} 
$$
Moreover for each of the components $F^r$ the notions of monogenicity and hermitian monogenicity coincide, since
$$
\dirac F^r = 0 \Longleftrightarrow (\upz - \upzd) F^r = 0 \Longleftrightarrow \left \{ \begin{array}{lcl}
\upz F^r & = & 0 \\[1mm] \upzd F^r & = & 0 \end{array} \right .
$$

In conclusion we have the following result.
\begin{proposition}
The function $F = \sum_{r=0}^n F^r$, $F^r : \mC^{2n} \longrightarrow \mS^r$, $r=0,\ldots,n$, is hermitian monogenic in a certain region $\Omega$ of $\mC^{2n}$ if and only if each of the components $F^r$, $r=0,\ldots,n$, is monogenic in that region.
\end{proposition}

The hermitian Dirac operators $\upz$ and $\upzd$ are invariant under the action of the group $\mbox{SO}_{\mI}(2n)$, or its double cover $\mbox{Spin}_{\mI}(2n)$, which is isomorphic  to the unitary group $\mbox{U}(n)$  (see Section 3). We may thus say that $\mbox{U}(n)$ is the fundamental group underlying the function theory of hermitian monogenic functions.\\[-2mm]

We will now move on to a further refinement of hermitian Clifford analysis, by considering the quaternionic structure $\mQ = (\mI_{4p},\mJ_{4p},\mK_{4p})$ on $\mR^{4p} \simeq \mH^p$, the dimension $m=2n=4p$ assumed to be a 4-fold. This will lead, in a first step, to a function theory for so-called {\em quaternionic Clifford analysis}, and in a second step to so-called {\em symplectic Clifford analysis}, where the fundamental invariance will be the one of the symplectic group $\mbox{Sp}(p)$. The most genuine way to introduce the new concept of quaternionic monogenicity is to directly generalize the system (\ref{hmon}), expressing hermitian monogenicity, now making use of the quaternionic structure on $\mR^{4p}$, whence the following definition.
\begin{definition}
\label{defqmon}
A differentiable function $F: \mR^{4p} \longrightarrow \mS$ is called quaternionic monogenic in some region $\Omega$ of $\mR^{4p}$, if and only if in that region $F$ is a solution of the system
\begin{equation}
\dirac F = 0, \quad \mI_{4p}[\dirac] F = 0, \quad \mJ_{4p}[\dirac] F = 0, \quad \mK_{4p}[\dirac] F = 0
\label{qmon}
\end{equation}
\end{definition}
Observe that, in a similar way as it was possible to introduce the notion of hermitian monogenicity without involving complex numbers, the above Definition \ref{defqmon} expresses the notion of quaternionic monogenicity without having to resort to quaternionic numbers.\\[-2mm]

There is, quite naturally, an alternative characterization of quaternionic monogenicity possible in terms of the hermitian Dirac operators, yet still not involving quaternions. We recall these hermitian Dirac operators in the actual dimension:
\begin{eqnarray*}
\upz & = & \sum_{k=1}^{2p} \p_{z_k} \gfd_k \ = \ \sum_{j=1}^p ( \p_{z_{2j-1}} \gfd_{2j-1} + \p_{z_{2j}} \gfd_{2j} ) \\
\upzd & = & \sum_{k=1}^{2p} \p_{\overline{z}_k} \gf_k \ = \ \sum_{j=1}^p ( \p_{\overline{z}_{2j-1}} \gf_{2j-1} + \p_{\overline{z}_{2j}} \gf_{2j} )
\end{eqnarray*}
and compute their images under the action of $\mJ_{4p}$:
\begin{eqnarray*}
\upzJ \ = \ \mJ_{4p} [\upz] & = & \sum_{j=1}^p ( \p_{z_{2j}} \gf_{2j-1} - \p_{z_{2j-1}} \gf_{2j} ) \\
\upzJd \ = \ \mJ_{4p} [\upzd] & = & \sum_{j=1}^p ( \p_{\overline{z}_{2j}} \gfd_{2j-1} - \p_{\overline{z}_{2j-1}} \gfd_{2j} )
\end{eqnarray*}
Here use is made of the formulae
$$
\mJ_{4p} [ \gf_{2j-1} ] = - \gfd_{2j}, \quad \mJ_{4p} [\gf_{2j}] = \gfd_{2j-1}, \quad \mJ_{4p}[\gfd_{2j-1}] = - \gf_{2j}, \quad \mbox{and} \quad \mJ_{4p}[\gfd_{2j}] = \gf_{2j-1}
$$
Now the original Dirac operator $\dirac$ and its twisted versions $\mI_{4p}[\dirac]$, $\mJ_{4p}[\dirac]$ and $\mK_{4p}[\dirac]$ may be expressed in terms of the hermitian Dirac operators $( \upz,\upzd )$ and their twisted versions $( \upzJ,\upzJd )$. We indeed have
$$
\upz = \frac{1}{4} ( {\bf 1} + i \, \mI_{4p} ) [\dirac], \quad \upzd = -\frac{1}{4} ( {\bf 1} - i \, \mI_{4p} ) [\dirac] 
$$
and
$$
\upzJ = \frac{1}{4} ( \mJ_{4p} + i \, \mK_{4p}) [\dirac], \quad \upzJd = -\frac{1}{4} ( \mJ_{4p} - i \, \mK_{4p}) [\dirac]
$$
whence conversely
\begin{eqnarray*}
\dirac & = & 2 ( \upz - \upzd) \\
i \, \mI_{4p} [\dirac] & = & 2 ( \upz + \upzd ) \\
\mJ_{4p} [ \dirac] &=& 2 ( \upzJ - \upzJd ) \\
i \, \mK_{4p}[\dirac] &=& 2 ( \upzJ + \upzJd) 
\end{eqnarray*}
Note that we could also have used the projection operators $\frac{1}{2} ( {\bf 1} \pm j \, \mJ_{4p})$ involving the quaternionic generator $j$ instead, to define two other differential operators, which then are linear combinations of the hermitian Dirac operators and their $\mJ$-twisted versions:
$$
\frac{1}{2} ( {\bf 1} \pm j \, \mJ_{4p} ) [\dirac] = \upz - \upzd \pm j \, \upzJ \mp \upzJd
$$ 
Our aim to avoid explicit use of quaternions for our function theories explains the choice made above for the differential operators, which now leads to an alternative characterization of quaternionic monogenicity.
\begin{proposition}
A differentiable function $F: \mR^{4p} \simeq \mC^{2p} \longrightarrow \mS$ is quaternionic monogenic in the region $\Omega \subset \mR^{4p}$ if and only if $F$ is in $\Omega$ a simultaneous null solution of the operators $\upz$, $\upzd$, $\upzJ$ and $\upzJd$.
\end{proposition}

As the identification of an underlying symmetry group is necessary for the further development of an acceptable and powerful function theory, the next result is crucial.
\begin{proposition}
The operators $\upz$, $\upzd$, $\upzJ$ and $\upzJd$ are invariant under the action of the symplectic group $\mbox{\em Sp}(p)$.
\end{proposition}

\pf
The action of a $\mbox{Spin}(4p)$--element $s$ on a spinor valued function $F$ is the so-called $L$--action given by $L(s) [ F(\uX)] = s F ( s^{-1} \uX s)$. The Dirac operator $\dirac$ is invariant under $\mbox{Spin}(4p)$ (or, equivalently, under $\mbox{SO}(4p)$), i.e.
$$
[L(s),\dirac] = 0, \qquad \mbox{for all } s \in \mbox{Spin}(4p)
$$
which can be explained by
$$
L(s) \dirac_{\uX} F(\uX) = s \dirac_{s^{-1} \uX s} F (s^{-1} \uX s) = s ( s^{-1} \dirac_{\uX} s ) F ( s^{-1} \uX s)  = \dirac_{\uX} L(s) F(\uX)
$$
We have seen in Section 3 that $\mbox{Sp}(p)$ is isomorphic to a subgroup of $\mbox{Spin}(4p)$, whence the Dirac operator $\dirac$ is, quite trivially, also invariant under the action of $\mbox{Sp}(p)$. The invariance of the operators $\upz$, $\upzd$, $\upzJ$ and $\upzJd$ now follows from the fact that their respective definitions involve projection operators which are commuting with the $\mbox{Sp}(p)$--elements.
\qed
In this respect, also note that the hermitian Dirac operators $\upz$ and $\upzd$ are invariant under $\mbox{U}(2p)$, see e.g.\ \cite{partI}.


\section{Conclusion}


In this paper we have studied the fundaments of so--called quaternionic Clifford analysis, a recent refinement of hermitian Clifford analysis, in its turn a refinement of Euclidean Clifford analysis. These are function theories for functions in Euclidean space of general but appropriate dimension, and taking their values in a real or complex Clifford algebra or subspaces thereof. For the development of a function theory it is not only important to establish the traditional fundamental results such as a Cauchy formula, or a Taylor series expansion, but it is evenly crucial to show the invariance of the underlying operator(s) with respect to a fundamental group. For Euclidean Clifford analysis in $\mR^m$ this is the SO$(m)$ group, for hermitian Clifford analysis in $\mR^{2n}$ it is the U$(n)$ group and for quaternionic Clifford analysis in $\mR^{4p}$ it turns out to be the Sp$( p)$ group. Only when this group invariance has been unraveled, it is possible to obtain another important result, the so-called Fischer decomposition of spaces of homogeneous polynomials into irreducible representations for the corresponding group. To that end it is also necessary to decompose the value space into such irreducible representations. For hermitian Clifford analysis the spinor space had to be split into its homogeneous parts. For quaternion Clifford analysis the splitting had to be carried out still one step further and spinor space was split into the so--called symplectic cells; this is the core of the present paper. Splitting the value space inevitably entails the splitting of the system of differential equations defining quaternionic monogenicity into subsystems, the study of which contributes to a better insight in the concept of quarternionic monogenicity. This study is carried out in detail in the forthcoming paper \cite{paper2}. The Fischer decomposition in the framework of quaternionic Clifford analysis itself is the subject of the forthcoming papers \cite{paper3,paper4}. \\[-2mm]

One final remark should be made about the aim of this paper, which is double. Quite naturally it contains new research results concerning the foundations of the function theory under development, such as the construction of the symplectic cells and the projection operators on those cells. However it also contains some expository parts of a more educational nature, such as the introduction of the complex and the quaternion structure and the study of the Lie algebra $\gsp_{2p}(\mC)$, which were added for the reader's comfort, since also there the approach is not the traditional one but is embedded in the structure of the Clifford algebra.



\end{document}